\documentclass[11pt,a4paper]{article}
\usepackage[english]{babel}
\usepackage{amsmath}
\usepackage{amssymb}
\usepackage{amsfonts}
\usepackage{amsthm}
\usepackage{graphicx}
\usepackage{fancyhdr}

\newtheorem{theorem}{Theorem}
\newtheorem{corollary}[theorem]{Corollary}
\newtheorem{lemma}[theorem]{Lemma}
\newtheorem{proposition}[theorem]{Proposition}
\newtheorem{definition}[theorem]{Definition}

\newenvironment{example}{\vspace{1.5ex}\noindent\textit{Example \refstepcounter{theorem}\thetheorem}}{}

\numberwithin{equation}{section}
\newcommand{\seq}[1]{\textnormal{\boldmath$#1$}}
\newcommand{\Ac}[1]{\mathcal{A}(#1)}

\begin{document}
\begin{center}
{\Large\textbf{The Fine Structure of Dyadically \\[1ex] Badly Approximable Numbers }}

\vspace{2ex}
\textsc{Johan Nilsson}

{\small \texttt{johann@maths.lth.se}}
\end{center}

\begin{abstract} 
We consider badly approximable numbers in the case of dyadic diophantine approximation. For the unit circle $\mathbb{S}$ and the smallest distance to an integer $\|\cdot\|$ we give elementary proofs that the set $F(c) = \{ x \in \mathbb{S}: \|2^nx\| \geq c\,,  n\geq 0\}$ is a fractal set whose Hausdorff dimension depends continuously on $c$, is constant on intervals which form a set of Lebesgue measure 1 and is self-similar. Hence it has a fractal graph. Moreover, the dimension of $F(c)$ is zero if and only if $c\geq 1-2\tau$, where $\tau$ is the Thue-Morse constant. We completely characterise the intervals where the dimension remains unchanged. As a consequence we can completely describe the graph of $ c\mapsto \dim_H \{x\in[0,1]: \|x-\frac{m}{2^n}\|< \frac{c}{2^n} \textnormal{ finitely often}\}$.
\end{abstract}

\vspace{1.5ex}

\noindent AMS Mathematics Subject Classification; 11J70 Continued fractions and generalizations, 68R15 Combinatorics on words

\section{Introduction}
Let $(X,d)$ be a metric space. Given a sequence $\{x_n\}\subset X$, of maybe random numbers, and a sequence $\{l_n\}$ of positive real numbers we define the following two sets $ I = \left\{ y\in X: d(x_n,y)<l_n \ \textnormal{infinitely often} \right\}$ and $F=X\setminus I$. By the notion diophantine approximation we shall mean the study of the sets $I$ and $F$. Let us make the following remark: if the sequence $\{x_n\}$ is dense in $X$ then $I$ is a non-empty and hence a residual set in the sense of Baire.

Consider the sequence $\{x_{n,m}\}_{n\in\mathbb{N},\, 0\leq m<n}$ with $x_{n,m} = \frac{m}{n}$ and where $\gcd(m,n)=1$ and with the particular choice of the sequence $l_n = \frac{1}{n^{\alpha}}$. For this special choice of $\{x_n\}$ and $\{l_n\}$ we are in the case of the classical diophantine approximation by rational numbers. It is a well know fact that if $\alpha>2$ then $F$ is non-empty while it is empty when $\alpha<2$.

Inspired  by the above example, we continue in this direction and refine the definition of the set $F$ to be the following set
\begin{equation}
\label{eq: def of F(a)}
F(\alpha) = \left\{ y\in X: d(x_{n,m},y)<\frac{1}{n^{\alpha}} \ \textnormal{finitely often} \right\}.
\end{equation} 
An interesting question is to look at the critical exponent, $\alpha_0$, such that $F(\alpha)$ is empty if $\alpha<\alpha_0$ and is non-empty when $\alpha>\alpha_0$. For this special value $\alpha_0$ we say that the set $F(\alpha_0)$ is the set of \textit{Badly Approximable Numbers}, $\textit{BAN}$.

A second step in refinement of (\ref{eq: def of F(a)}) is to introduce the dependence on an extra parameter $c$,
\[
F_c(\alpha) = \left\{ y\in X: d(x_{n,m},y)<\frac{c}{n^{\alpha}} \ \textnormal{finitely often} \right\}.
\]
In the one-dimensional case this refinement leads to the area of continued fraction, which was first systematically studied by the Dutch astronomer Huygens in the 17-th century, motivated by technical problems while constructing a model of our solar system. Briefly, the continued fraction for a real $x$ is,
\[
x = a_0 + \cfrac{1}{ a_1 + \cfrac{1}{ a_2 + \cfrac{1}{a_3+\ldots}}}
\]
where $a_n\in \mathbb{N}$ are called \textit{partial denominators}. For brevity the continued fractions is often denoted by $[a_0, a_1, a_2, \ldots]$. The following theorem gives a neat connection between the badly approximable numbers and the continued fractions, for a proof see \cite{khintchine}.

\begin{theorem}
An irrational $x$ is  a $\textit{BAN}$ if and only if its partial denominators are bounded.
\end{theorem}

Yet another version, or refinement, of the $F$ set can be introduced via a condition on the partial denominators. We set
\[
F_N(2) = \left\{ x : x = [a_0,a_1,a_2,\ldots] \textnormal{ with } \ a_j < N  \right\}.
\]
The theory of iterated function system, \emph{IFS}-theory, gives an implicit formula for the Hausdorff dimension, $\dim_H$, of $F_N(2)$. The set $F_c(2)$ is finer as  $F_N(2)$ counts only the maximal $a_i$ while the $F_c(2)$ takes into account all $a_i$. In 1891 Hurwitz found that if $c<\frac{1}{\sqrt{5}}$ then $F_c(2)$ is empty and moreover the constant $\frac{1}{\sqrt{5}}$ is the best possible, but otherwise little is known about the set $F_c(2)$.

\vspace{1.5ex}
In this paper we are going to study a special case of diophantine approximation, approximation by dyadic rationals. Similar to the approximation by rationals we set the sequences $\{x_n\}$ and $\{l_n\}$ to be
\[
x_{n,m} = \frac{m}{2^n} \quad\textnormal{and}\quad l_{n} = \frac{c}{2^n}
\]
for $m$ odd.
We will turn our interest to the same type of questions as in the classical approximation case, and look at the set of badly approximable numbers in the dyadic case. We define $\hat{F}(c)$ to be the set 
\begin{equation}
\label{eq: def of F tilde}
	\hat{F}(c) = \left\{ x\in\mathbb{S}: \left\|x-\frac{m}{2^n}\right\| < \frac{c}{2^n}  \ \textnormal{finitely often} \right\},
\end{equation}
where $\|\cdot\|$ is the shortest distance to an integer. As we are going to study dimensional properties of $\hat{F(c)}$ we can restrict ourselves to the case when the inequality condition in (\ref{eq: def of F tilde}) is never fulfilled. So we define $F(c)$ by
\[
F(c) = \left\{ x\in\mathbb{S}:  \left\|2^n x  \right\| \geq c\, \textnormal{ for all $n\geq0$} \right\}
\] 
and we define the dimension function $\phi:(0,1)\to [0,1]$ such that $\phi(c) = \dim_H F(c)$. Then $\hat{F}(c)$ is the countable union of pre-images of $F(c)$ under multiplication by two. Hence the dimension does not change.

We prove that $\phi$ has derivative zero Lebesgue \mbox{a.e.}, that $\phi$ is continuous and is self-similar. Moreover we prove that the complementary zero-set, to where the derivative of $\phi$ is zero, has full Hausdorff dimension and we give the complete characterisation of the intervals where the derivative of $\phi$ is zero.

\subsection{Symbolic Dynamics}

Let $\Sigma_n^{} = \{0,1,\ldots,n-1\}^{\mathbb{N}} = \{ \,\seq{x}= x_1x_2\ldots: x_i\in \{0,1,\ldots,n-1 \} \} $ be the space of the one-sided infinite sequences on $n$ symbols, equipped with the product topology. Let similarly $\Sigma_n^* = \{ \,\seq{x}= x_1x_2\ldots x_m: x_i\in \{0,1,\ldots,n-1 \}, m\in \mathbb{N} \}$ be the set of all finite  sequences on $n$ symbols. There is a natural embedding of the finite sequences into the  set of infinite sequences, we can interpret a finite sequence as an infinite sequence ending with zeros. This gives that we can use the standard lexicographical order to compare sequences.

We are mainly going to consider sequences in $\Sigma_2^{}$ and $\Sigma_2^*$, the binary sequences. Therefore by the word \emph{sequence} we shall mean a binary sequence, finite or infinite, if not explicitly stated otherwise. The word sequence is used both for a finite sequence as well as for an infinite one. 

There is a correspondence between $\Sigma_2^{}$ and the real interval $[0,1]$, by simply considering the binary expansion of a real number. That is, for $x\in[0,1]$ we have
\begin{equation}
\label{eq: binary expansion}
	x = \sum_{i=1}^{\infty} \frac{x_i}{2^i}\quad \textnormal{with $x_i\in\{0,1\}$}
\end{equation}
and we let $\seq{x} = x_1x_2x_3\ldots$. This correspondence is one-to-one except for a countable set where it is two-to-one, but this will not cause us any trouble. We introduce here some notation that will be used.  

\begin{itemize}
\item By a \emph{concatenation} we mean that we append a sequence to a finite sequence, that is, the concatenation of $\seq{u}$ and $\seq{v}$ is $\seq{u}\seq{v}\,$, similarly we write $\seq{uu} = \seq{u}^2$.

\item We say that $\seq{x}$ is a \emph{prefix} of $\seq{s}$ if there exists a sequence $\seq{u}$ such that $\seq{s}= \seq{x}\seq{u}$ and similarly we then say that $\seq{u}$ is a \emph{suffix} of $\seq{s}$. If $\seq{u}$ is non-void then $\seq{x}$ is a proper prefix and similarly for a suffix.

\item By  $\seq{s}(k,n)$ we mean the sub-sequence $\seq{s}[k,n] = s_k s_{k+1}\ldots s_n$. And for a set $A$ of sequences the notation $A[k,n]$ is the set of sub-sequences, $A[k,n] = \{ \seq{s}[k,n]: \seq{s}\in A \}$.

\item The notation $|\cdot|$ will mean the length of a sequence, that is $|\seq{s}(k,n)| = n-k+1$. We will also use the $|\cdot|$-notation for the cardinality of a set. 

\item For a sequence $\seq{x}$, not necessarily binary, we define the left-shift $\sigma$ by $(\sigma(\seq{x}))_i = x_{i+1}$ and we let $\sigma^n = \sigma \circ \sigma^{n-1}$. If  $\seq{x}$ is a finite sequence then $\sigma^{|\seq{x}|}(\seq{x})$ is the empty sequence.  

\item By the notation $\seq{x}^*$ we mean the sequence $\seq{x}$ where we have changed zeros to ones and vice versa, the bit-wise inverse of $\seq{x}$. If $\seq{x}$ is finite then $\seq{x}^*$ can be seen as the inverse element of $\seq{x}$ in $\Sigma_2^*$. 

\item The notation $\seq{x}'$ will mean the inverse when seeing $\seq{x}$ as a real number, that is the inverse element of $\seq{x}$ in $\Sigma_2^{}$. If $\seq{x}$ is an infinite sequence then $\seq{x}^* = \seq{x}'$ but this equality does not hold in the finite case, as we then have to cast $\seq{x}$ to an element in $\Sigma_2^{}$, \mbox{i.e.} we have to append zeros at the end, that is $1' = (10^{\infty})' = 01^{\infty} = 1$ but $1^* =0$. We will always let $|\seq{x}| = |\seq{x}'|$.

\item For a finite sequence $\seq{x}$ the notation $\tilde{\seq{x}}$ means the sequence where the last symbol of $\seq{x}$ has been inverted.
\end{itemize}

Let $A$ be a square $\{0,1\}$ matrix with rows and columns indexed by the numbers $\{0,1,\ldots, n-1\}$. The matrix $A$ defines a closed, shift invariant subset of $\Sigma_n^{}$. It is defined by choosing the sequences as
\[
\left\{ \seq{x}\in\Sigma_n^{}: A_{x_ix_{i+1}} =1\,\textnormal{for all $i>0$}\right\}
\]
The dynamical system and the restriction of the shift transformation is the \emph{one-sided sub-shift of finite type defined by $A$.} We call such a matrix $A$ a transition matrix. The representation of a sub-shift via a transition matrix is not unique, two different matrices $A$ and $B$ may describe the same sub-shift. We say that a transition matrix $A$ is \emph{irreducible} if there for each pair of indices $i,j$ exists an $n$ such that $(A^n)_{ij}>0$.
Similarly, if there is an $m$ such that $(A^m)_{ij}>0$ for all pairs $i,j$ we say that the matrix is \emph{primitive}. Clearly primitivity implies irreducibility. A sub-shift of finite type is topological transitive if and only if it can be represented by an irreducible transition matrix and a sub-shift of finite type is topological mixing if and only if it can be represented by a primitive transition matrix.

For irreducible transition matrices we have the useful Perron-Frobenuis theorem, (see \cite{kitchens, perron}).

\begin{theorem}[Perron-Frobenius]	
\label{thm: perron-frobenius}
Suppose $A$ is a nonnegative, square matrix. If $A$ is irreducible there exists a real eigenvalue $\lambda>0$ such that 
\begin{enumerate}
\item $\lambda$ is a simple root of the characteristic polynomial;
\item $\lambda$ has strictly positive left and right eigenvectors;
\item the eigenvectors for $\lambda$ are unique up to constant multiple;
\item $\lambda>|\mu|$, where $\mu$ is any other eigenvalue;
\item if $0\leq B \leq A$ and $\beta$ is an eigenvalue for $B$ then 
		$|\beta|\leq\lambda$ and equality occurs if and only if $B=A$. 
\end{enumerate}
The special eigenvalue $\lambda$, is the Perron value of the matrix $A$.
A positive eigenvector corresponding to $\lambda$ is called a Perron eigenvector.
\end{theorem}

Note that the notion of Perron value coincides for non-negative irreducible matrices with the notion of spectral radius $\rho(A)$.

By coding each symbol in the alphabet $\{0,1,\ldots ,n-1\}$ with a finite word of zeros and ones, the transition matrix selects valid shifts in infinite binary sequences. That is, we index the rows and columns in $A$ by binary words of a fixed length. 

From (\ref{eq: binary expansion}) we see that multiplication by 2 of a real number $x\in[0,1]$ corresponds to shift the corresponding sequence $\seq{x}$ leftward once. Hence the investigation of the set $F(c)$ can now be turned to the investigation of the set  of sequences
\begin{equation}
\label{eq: def of F by seq}
	F(\seq{c}) = \left\{\seq{x}\in\Sigma_2^{}: \seq{c}'\geq \sigma^n(\seq{x})\geq \seq{c}\textnormal{ for all $n\geq0$} \right\}.
\end{equation}
From (\ref{eq: def of F by seq}) we have that if the sequence $\seq{c}$ is of finite length then $F(\seq{c})$ can be described by a transition matrix $A_{\seq{c}}$. Note that it is only a sufficient condition that $\seq{c}$ should be of finite length to be able to describe the set $F(\seq{c})$ by a transition matrix.

In \cite{allouche83:1, allouche83:2}, (see also \cite{allouche99}) Allouche and Cosnard consider iterations of unimodal functions. (A continuous function $f$ is said to be unimodal if for $a\in(0,1)$, $f(1)=0$ and $f(a) =1$, it is strictly increasing on $[0,a)$ and strictly decreasing on $(a,1]$). They give the result that the existence of unimodal functions is connected to elements in the set of binary sequences $\Gamma$, where 
\begin{equation}
\label{eq: allouche set}
\Gamma = \{ \seq{x}\in\Sigma_2: \seq{x}' \leq \sigma^n(\seq{x}) \leq \seq{x} \textnormal{ for all $n\geq0$} \}.
\end{equation}
Allouche and Cosnard presents some properties of the set $\Gamma$. They show that it is a self similar set and therefore a fractal set. In Corollary \ref{cor: dim Gamma = dim F} we show that the dimensional structure of $\Gamma$ is the same as the dimensional structure of $F(\seq{c})$. Furthermore in \cite{allouche83:1, allouche83:2}, Allouche and Cosnard consider also the more general set $\Gamma_{\seq{a}}$, where 
\begin{equation}
\label{eq: allouche set a}%
\Gamma_{\seq{a}} = \{ \seq{x}\in\Sigma_2: \seq{a}' \leq \sigma^n(\seq{x}) \leq \seq{a}  \textnormal{ for all $n\geq0$} \}.
\end{equation}
One of the main results achieved by Allouche and Cosnard on $\Gamma_{\seq{a}}$ is to present the threshold sequence $\seq{t}_2$ such that $\Gamma_{\seq{a}}$ is countable if and only if $\seq{a}<\seq{t}_2$. In \cite{moreira:2001}, Moreira improves this result and shows that $\dim_H \Gamma_{\seq{a}} = 0$ if and only if $\seq{a}\leq\seq{t}_2$. 

Moreira also turn his interest to how the dimension of sets like $\Gamma_{\seq{a}}$ depends on the parameter $\seq{a}$. In \cite{labarca}, Labarca and Moreira show that for $(\seq{a},\seq{b})\in \Sigma_2\times\Sigma_2$ the map 
\[	
	(\seq{a},\seq{b})\mapsto \dim_H\{\seq{x}\in \Sigma_2: \seq{a}\geq \sigma^n(\seq{x})\geq \seq{b}\textnormal{ for all $n\geq0$}\}
\]
is continuous in both $\seq{a}$ and $\seq{b}$. We simplify the proof given by Labarca and Moreira, and present an elementary proof that $\seq{c}\mapsto \dim_H F(\seq{c})$ is continuous. In Section \ref{sec: shift-bounded} we present in more detail some technical results by Allouche and Cosnard that we will make use of.

\subsection{Dimension}

Let us start with the notion of Hausdorff dimension.

\begin{definition}
\label{def: hausdorff}
Let $s\in[0,\infty]$. The $s$-dimensional Hausdorff measure $\mathcal{H}^s$(Y) of a subset of a metric space $X$ is defined by 
\[
\mathcal{H}^s(Y) = \lim_{\varepsilon\to0} \inf\left\{ \ 
\sum_{i=1}^{\infty} \textnormal{diam}(U_i)^s: Y\subset \bigcup_{i=1}^{\infty}U_i\,, \, \sup_i \textnormal{diam}(U_i) \leq \varepsilon \ \right\}.
 \]	
The unique $s_0$ such that 
\[
\mathcal{H}^s(Y) =\left\{
\begin{array}{ll}
	\infty & \textnormal{for } s< s_0\\	
	0 & \textnormal{for } s> s_0
\end{array}
\right.
\]
we call the Hausdorff dimension of the set $Y$ and it will be denoted by $\dim_H Y $.
\end{definition}

A way of estimating the Hausdorff dimension of a set is to use the connection between the H\"older exponent and the Hausdorff dimension. The following result is well known.

\begin{proposition}
\label{prop: holder-dim}
Let $X\subset\mathbb{R}^n$ and suppose that $f:X\to\mathbb{R}^m$ satisfies a H\"older condition
\[
|f(x)-f(y)|\leq C \, |x-y|^{\alpha}\quad (x,y\in X).
\]
Then $\dim_H f(X) \leq \frac{1}{\alpha} \dim_H X$.
\end{proposition}

For a deeper discussion of dimension theory and methods used therein see Falconer's book \cite{falconer}. 
Recall that by $F(\seq{c})[1,n]$ we denote the set of prefixes of length $n$ of sequences in $F(\seq{c})$, that is, $F(\seq{c})[1,n] = \{\seq{x}[1,n] : \seq{x}\in F(\seq{c}) \}$.

\begin{definition}
We define the topological entropy $h_\textnormal{top}$ of the set $F(\seq{c})$ as the growth rate of the number of sequences allowed as the length $n$ increases,
\[
\label{eq: def of top entropy}
h_\textnormal{top} (F(\seq{c})) \ = \ \lim_{n\to\infty} \frac{1}{n}\log \left(|F(\seq{c})[1,n]|\right),
\]
where $|\cdot|$ denotes the cardinality of a set.
\end{definition}

The existence of the above limit follows by simply noticing the sub-additivity property of the function $n\mapsto \log\left(|F(\seq{c})[1,n]|\right)$:
\[
\log \left( |F(\seq{c})[1,n+m]|\right)  \leq  \log\left(|F(\seq{c})[1,n]|\right)+\log\left(|F(\seq{c})[1,m]|\right).
\]
In the case when $F(\seq{c})$ is a sub-shift of finite type the existence of the limit implies that there exists constants $k_1$ and $k_2$ with $k_1 \lambda^n \leq |F(\seq{c})[1,n]| \leq k_2 \lambda^n$ for all sufficiently large $n$.

\begin{theorem}
\label{thm: pesin}
Let $F(\seq{c})$ be a sub-shift of finite type described by the transition matrix $A_{\seq{c}}$, with the spectral radius $\rho(A_{\seq{c}})$. Then 
\begin{enumerate}
\item $h_\textnormal{top}(F(\seq{c})) = \log \rho(A_{\seq{c}})$;
\item $\displaystyle{\dim_H F(\seq{c}) = \frac{\log \rho(A_{\seq{c}}) }{\log 2}}$.
\end{enumerate}
\end{theorem}

Theorem \ref{thm: pesin} gives a link between the topological entropy and the Hausdorff dimension via transition matrices for sub-shifts of finite type. For a proof of Theorem \ref{thm: pesin} see Pesin's book \cite{pesin} on dimension theory.

\section{Fundamental Properties}

Let us start with an example on the structure of $F(\seq{c})$ for a special choice of $\seq{c}$.  

\begin{example}
The set $F(0^k1)$ is the set of sequences from $\Sigma_2^{}$ containing at most $k$ consecutive zeros and $k$ consecutive ones. In particular $F(01)$ is the set containing only the two elements $(01)^{\infty}$ and $(10)^{\infty}$.
\end{example}

\begin{lemma}
\label{lemma: lambda}
For $k>0$ we have $\dim_H F(0^k1) = \frac{\log \lambda_k}{\log 2}$ where $\lambda_k $ is the largest real root of $\lambda^k = \lambda^{k-1} + \lambda^{k-2}+ \ldots + \lambda + 1$.
\end{lemma}

\begin{proof}
The polynomial equation is obtained by calculating the number of allowed words of a given length in $F(\seq{c})$. Combine this with Theorem \ref{thm: pesin}.
\end{proof}

\begin{lemma} 
\label{lemma: Fc = Fcoo}
Let $\seq{c}$ be a finite non-empty sequence. Then $F(\seq{c}) = F(\seq{c}^{\infty})$.
\end{lemma}
\begin{proof}
From (\ref{eq: def of F by seq}) it is clear that $F(\seq{c}^{\infty}) \subset F(\seq{c})$ as $\seq{c} < \seq{c}^{\infty}$. For the converse, let $\seq{x}\in F(\seq{c}) \setminus F(\seq{c}^{\infty})$. Then there is an $n$ such that either $\sigma^n(\seq{x})< \seq{c}^{\infty}$ or $(\seq{c}^{\infty})'<\sigma^n(\seq{x})$. Assume that $\sigma^n(\seq{x})< \seq{c}^{\infty}$ and let $k$ be the first position where $\sigma^n(\seq{x})$ differs from $\seq{c}^{\infty}$.  

\footnotesize
\[
\begin{picture}(230,55)
\multiput(40,27)( 0,0){1}{
\put( -35,  0){\makebox(30,12)[r]{$\sigma^n(\seq{x})=\ldots$}}
\multiput(0,0)( 0,0){1}{
\put(   0,  0){\line(1,0){193}}\put(  0,12){\line(1,0){193}}}
}
\multiput(40,5)( 0,0){1}{
\put(-35,  0){\makebox(30,12)[r]{$\seq{c}^{\infty}=\ldots$}}
\multiput(0,0)(0,0){1}{
\put(   0,  0){\line(1,0){ 30}}\put(   0,12){\line(1,0){ 30}}
\put(  30, 0){\line(0,1){ 12}}
\put(  13,  4){$\seq{c}$}}
\multiput(31,0)(46,0){3}{
\put(   0,  0){\line(1,0){ 45}}\put(   0,12){\line(1,0){ 45}}
\put(   0,  0){\line(0,1){ 12}}\put(  45, 0){\line(0,1){ 12}}
\put(  20,  4){$\seq{c}$}}
\multiput(169,0)(51,0){1}{
\put(   0,  0){\line(1,0){ 24}}\put(   0,12){\line(1,0){ 24}}
\put(   0,  0){\line(0,1){ 12}}}
}
\multiput(171,5)(0,4){9}{\line(0,1){2}}
\put(169,44){$k$}
\end{picture}
\]
\normalsize
We can write $k = m|\seq{c}|+ r$ for some non-negative integers $m, r$ with $r<|\seq{c}|$.  But then we must have $\sigma^{n+m|\seq{c}|}(\seq{x}) < \seq{c}$, and it follows that $\seq{x} \notin F(\seq{c})$. The case when $(\seq{c}^{\infty})'<\sigma^n(\seq{x})$ is treated in the same way.
\end{proof}

\begin{lemma}
\label{lemma: uuu}
Let $\seq{c}$ be a non-empty sequence of the form $\seq{c} = \tilde{\seq{u}}\,(\seq{u}^*)^k\seq{u}'\seq{v}$ for some $k\geq0$ and a finite non-empty sequence $\seq{u}$. If $\seq{x}\in F(\seq{c})$ contains the subsequence $\tilde{\seq{u}}$, (or symmetrically $\seq{u}'$), then $\seq{x}$ must be of the form 
\begin{equation}
\label{eq: uuu}
	\seq{w}\,\tilde{\seq{u}}\,(\seq{u}^*)^{k_1}\,
	\seq{u}' \, \seq{u}^{k_2}\,
	\tilde{\seq{u}}\,(\seq{u}^*)^{k_3}\,
	\seq{u}' \, \seq{u}^{k_4}\,	\tilde{\seq{u}}\ldots,
\end{equation}
with $0\leq k_i\leq k$ and where the sequence $\seq{w}$ does not contain the subsequence $\tilde{\seq{u}}$.
\end{lemma}

\begin{proof}
Let $n$ be the smallest integer such that $\sigma^n(\seq{x})= \tilde{\seq{u}}\ldots$. Let $\sigma^n(\seq{x}) = \tilde{\seq{u}}\, \seq{a}_1\,\seq{a}_2\ldots$, with $|\seq{a}_i| = |\seq{u}|$. Let $m$ be the smallest integer such that $\seq{a}_m\neq \seq{u}^*$. From the inequality $\sigma^n(\seq{x}) =\tilde{\seq{u}}\, \seq{a}_1\,\seq{a}_2\ldots \geq \tilde{\seq{u}}\,(\seq{u}^*)^k\seq{u}'\seq{v}$ we have that $1\leq m \leq k+1$. 

\footnotesize
\[
\begin{picture}(230,55)
\multiput(40,32)( 0,0){1}{
\put( -35,  0){\makebox(30,12)[r]{$\sigma^n(\seq{x})=$}}
\multiput(0,0)(41,0){4}{
\put(   0,  0){\line(1,0){ 40}}\put(   0,12){\line(1,0){ 40}}
\put(   0,  0){\line(0,1){ 12}}\put(  40, 0){\line(0,1){ 12}}}
\multiput(164,0)(41,0){1}{
\put(   0,  0){\line(1,0){ 29}}\put(   0,12){\line(1,0){ 29}}
\put(   0,  0){\line(0,1){ 12}}}
\put(  17,  4){$\tilde{\seq{u}}$}
\put(  58,  4){$\seq{a}_1$}
\put(  99,  4){$\seq{a}_2$}
\put( 140,  4){$\seq{a}_3$}
}
\multiput(40,10)( 0,0){1}{
\put(-35,  0){\makebox(30,12)[r]{$\seq{c}=$}}
\multiput(0,0)(41,0){4}{
\put(   0,  0){\line(1,0){ 40}}\put(   0,12){\line(1,0){ 40}}
\put(   0,  0){\line(0,1){ 12}}\put(  40, 0){\line(0,1){ 12}}}
\multiput(164,0)(41,0){1}{
\put(   0,  0){\line(1,0){ 29}}\put(   0,12){\line(1,0){ 29}}
\put(   0,  0){\line(0,1){ 12}}}
\put(  17,  4){$\tilde{\seq{u}}$}
\put(  58,  4){$\seq{u}^*$}
\put(  99,  4){$\seq{u}^*$}
\put( 140,  4){$\seq{u}^*$}
}
\end{picture}
\]
\normalsize
This implies that $\seq{a}_m \geq \seq{u}'$. By shifting the sequence $\seq{x}$ additionally $m|\seq{u}|$ times we obtain $\seq{u}'\geq \seq{a}_m$.
\footnotesize
\[
\begin{picture}(230,55)
\multiput(40,32)( 0,0){1}{
\put(-35,  0){\makebox(30,12)[r]{$\seq{c}'=$}}
\multiput(0,0)(41,0){4}{
\put(   0,  0){\line(1,0){ 40}}\put(   0,12){\line(1,0){ 40}}
\put(   0,  0){\line(0,1){ 12}}\put(  40, 0){\line(0,1){ 12}}}
\multiput(164,0)(41,0){1}{
\put(   0,  0){\line(1,0){ 29}}\put(   0,12){\line(1,0){ 29}}
\put(   0,  0){\line(0,1){ 12}}}
\put(  17,  4){$\seq{u}'$}
\put(  58,  4){$\seq{u}$}
\put(  99,  4){$\seq{u}$}
\put( 140,  4){$\seq{u}$}
}
\multiput(40,10)( 0,0){1}{
\put( -35,  0){\makebox(30,12)[r]{$\sigma^{n+m|\seq{u}|}(\seq{x})=$}}
\multiput(0,0)(41,0){4}{
\put(   0,  0){\line(1,0){ 40}}\put(   0,12){\line(1,0){ 40}}
\put(   0,  0){\line(0,1){ 12}}\put(  40, 0){\line(0,1){ 12}}}
\multiput(164,0)(41,0){1}{
\put(   0,  0){\line(1,0){ 29}}\put(   0,12){\line(1,0){ 29}}
\put(   0,  0){\line(0,1){ 12}}}
\put(  15,  4){$\seq{a}_{m  }$}
\put(  51,  4){$\seq{a}_{m+1}$}
\put(  92,  4){$\seq{a}_{m+2}$}
\put( 133,  4){$\seq{a}_{m+3}$}
}
\end{picture}
\]
\normalsize
Hence $\seq{u}'= \seq{a}_m$. The result now follows by symmetry. 
\end{proof}

\vspace{1.5ex}
For the special case when $k = 0$ in Lemma \ref{lemma: uuu} we have the following corollary, which also was given by Allouche in \cite{allouche83:1}.

\begin{corollary}[Allouche \cite{allouche83:1}]
\label{cor: aa*}
Let $\seq{c}$ be a finite non-empty sequence of the form $\seq{c} = \seq{u}\,\seq{u}^*$. If $\seq{x}\in F(\seq{c})$ contains the subsequence $\seq{u}$, (or symmetrically $\seq{u}^*$), then $\seq{x}$ must be 
of the form
\[
\seq{w}(\seq{u}\,\seq{u}^*)^{\infty}
\]
for some sequence $\seq{w}$ not containing the subsequence $\seq{u}$.
\end{corollary}

\section{Shift-Bounded Sequences}
\label{sec: shift-bounded}

\begin{definition}
\label{def: SB}
A finite sequence $\seq{s}$ fulfilling $\seq{s}'>\sigma^n(\seq{s})>\seq{s}$ for $0<n<|\seq{s}|$ is said to be a finite shift-bounded sequence. Similarly, an infinite sequence $\seq{s}$ fulfilling $\seq{s}'>\sigma^n(\seq{s})>\seq{s}$ for all $n>0$ is said to be an infinite shift-bounded sequence. For completeness we also say that the sequence $\seq{s}=1$ is shift-bounded, while the sequence $\seq{s}=0$ is not.
\end{definition}

Our definition of shift-bounded sequences coincides with and extends the definition of  \emph{admissible} sequences considered by Komornik and Loreti in \cite{komornik} and by Allouche and Cosnard in \cite{allouche01}. From the definition of a shift-bounded sequence we have directly the following important proposition

\begin{lemma}
\label{lemma: prefix vs suffix}
Let $\seq{s}$ be a finite shift-bounded sequence and let $\seq{\alpha}$ and $\seq{\gamma}$ be a prefix and a suffix respectively of $\seq{s}$ such that $0<|\seq{\alpha}|=|\seq{\gamma}|<|\seq{s}|$. Then $ \seq{\alpha}^* \geq \seq{\gamma} > \seq{\alpha}$ and $\seq{\alpha}^* > \tilde{\seq{\gamma}} \geq \seq{\alpha}$.
\end{lemma}

\begin{lemma}
\label{lemma: s = aa*}
Let $\seq{s}$ be a shift-bounded sequence. If $\seq{a}\,\seq{a}^*$ is a prefix of $\seq{s}$, where $\seq{a}$ is non-empty, then $\seq{s} = \seq{a}\,\seq{a}^*$.
\end{lemma}

\begin{proof}
Let $N$ be the maximal integer such that $\seq{s} = (\seq{a}\,\seq{a}^*)^N\seq{u}$ for some sequence $\seq{u}$. This number $N$ exists, since otherwise $\seq{s}$ would be periodic and hence not shift-bounded. By the shift-boundedness of $\seq{s}$ we have that $\seq{u}$ may not be empty. If $\seq{u}$ does not have $\seq{a}$ as a prefix we have, since $\seq{s}$ is shift-bounded, $\seq{a}^*\seq{a}\,\seq{u}'> \sigma^{2|\seq{a}|N-|\seq{a}|}(\seq{s}) = \seq{a}^*\seq{u}$. But also $\sigma^{2|\seq{a}|N}(\seq{s}) = \seq{u} > \seq{a}$. Hence $\seq{a}\,\seq{u}' > \seq{u} > \seq{a}$, a contradiction as $\seq{a}$ is not a prefix of $\seq{u}$. For the second case, if $\seq{u}$ has $\seq{a}$ as prefix we can write $\seq{s} = (\seq{a}\,\seq{a}^*)^N\seq{a}\,\seq{v}$ for some sequence $\seq{v}$ not having $\seq{a}^*$ as prefix. This gives $\seq{a}^*\seq{a}\,\seq{a}^*\seq{v}'> \sigma^{2|\seq{a}|N-|\seq{a}|}(\seq{s}) = \seq{a}^*\seq{a}\,\seq{v}$ and
$\sigma^{2|\seq{a}|N}(\seq{s}) = \seq{a}\,\seq{v} > \seq{a}\,\seq{a}^*$, that is, $\seq{a}^*\seq{v}'>\seq{v}>\seq{a}^*$, a contradiction.
\end{proof}

\begin{corollary}
\label{cor: b<a*}
Let $\seq{s}$ be a shift-bounded sequence and let $\seq{s} = \seq{abv}$ with $|\seq{a}| = |\seq{b}|$ and $|\seq{v}|>0$. Then $\seq{b}<\seq{a}^*$.
\end{corollary}

\begin{definition}
\label{def: f}
For a finite sequence $\seq{s}$ ending with a 1, we define the map $f$ by $f(\seq{s})= \tilde{\seq{s}}\,\seq{s}'$. We define the function $d$ as the function taking $\seq{s}$ to its limit point under self-composition of $f$,
\[
d(\seq{s})  = \lim_{k\to\infty}f^k(\seq{s}).
\]
\end{definition}

The limit in the definition of $d$ above exists as the function $f$ has a strictly decreasing and from below bounded orbit.  
In \cite{allouche83:1,allouche83:2, allouche01} Allouche and Cosnard consider a function $\varphi$ defined on periodic sequences by $\varphi( (\seq{a}0)^{\infty} ) = (\seq{a}0\seq{a}^* 1)^{\infty}$. Our function $f$ is $\varphi$ defined on finite sequences. The sequence $d(\seq{c})$ coincides with Allouche and Cosnard's notion of $q$-mirror sequences, where the $q$ is the lenght of $\seq{c}$.

\begin{lemma}[Allouche, Cosnard \cite{allouche83:1, allouche01}]
Let $\Gamma$ be the set defined in (\ref{eq: allouche set}) and let $\seq{x} = (\seq{a}0)^{\infty}$, where $(\seq{a}0)$ is the shortest period. Then $\seq{x}\in \Gamma$ if and only if $\varphi(\seq{x}) \in \Gamma$ and moreover $\seq{x}\in \Gamma$ if and only if $\lim_{n\to\infty}\varphi^n(\seq{x}) \in \Gamma$.
\end{lemma}

\begin{corollary}
Let $\seq{s}$ be a finite sequence. Then $\seq{s}$ is shift-bounded if and only if $f(\seq{s})$ is shift-bounded. Moreover $\seq{s}$ is shift-bounded if and only if $d(\seq{s})$ is shift-bounded.
\end{corollary}

The function $f$ could equally have been defined on the rational numbers. By a straight forward calculation we have

\begin{theorem}
\label{thm: transendental}
Let $\frac{a}{2^n}\in \mathbb{Q}^+\setminus \{0\}$. Then the limit 
\[	
\lim_{k\to\infty} \ f^k(\frac{a}{2^n}) = \frac{a}{2^n}\prod_{i=0}^{\infty}\left(1-\frac{1}{2^{2^in}}\right)
\]
is a well defined real number and moreover it is a transcendental number.
\end{theorem}

The second part of the theorem is a direct consequence of the following theorem by Mahler,

\begin{theorem}[Mahler \cite{mahler}]
\label{thm: mahler}
Let $0<|a|<1$ be an algebraic number. Then the product $\prod_{i=0}^{\infty}(1-a^{2^i})$ is transcendental.
\end{theorem}

The function $f$ is connected to the classical Thue-Morse sequence.

\begin{definition}[Thue-Morse sequence]
The sequence $\seq{t}$ recursively defined by $t_1=0$ and $t_{2n+1}^{}=t_{n+1}^{}$, $t_{2n+2}^{} = t^*_{n+1}$, is called the Thue-Morse sequence. 
\end{definition}

The first entries of the Thue- Morse sequence $\seq{t}$ and its inverse are easily seen to be 
\begin{eqnarray*}
\seq{t}
&=& 0\,1\,1\,0\,1\,0\,0\,1\,1\,0\,0\,1\,0\,1\,1\,0\,1\,0\,0\,1\,0\,1\,1\,0\,0\,1\,1\,0\,1\,0\,0\,1\ldots\\[1ex]
\seq{t}'
&=& 1\,0\,0\,1\,0\,1\,1\,0\,0\,1\,1\,0\,1\,0\,0\,1\,0\,1\,1\,0\,1\,0\,0\,1\,1\,0\,0\,1\,0\,1\,1\,0\ldots
\end{eqnarray*}
The Thue-Morse sequence $\seq{t}$ is widely studied and it appears in several different areas. We reefer to \cite{allouche99} for more about this and further references.  In \cite{allouche83:1,allouche83:2}, Allouche and Cosnard observed that we may obtain the Thue-Morse sequence via the limit under iteration of the function $f$. 

\begin{theorem}[Allouche, Cosnard \cite{allouche83:1,allouche83:2}]
The sequence $\seq{d} = d(1)$ is the shift\-ed inverse Thue-Morse sequence $\seq{t}$\,, \mbox{i.e.} $\seq{d} = \sigma(\seq{t}')$. In particular $\seq{d}$ is shift-bounded.
\end{theorem}

From Theorem \ref{thm: transendental} and Mahler's Theorem \ref{thm: mahler} we have directly the transcendence of the Thue-Morse constant. 

\begin{theorem}[Dekking \cite{dekking}, Mahler \cite{mahler}]
The Thue-Morse constant $\tau = \sum_{i=1}^{\infty}\frac{t_i}{2^{i+1}}=0.41245403\ldots$ is transcendental, where $\seq{t}$ is the Thue-Morse sequence.
\end{theorem}

\begin{lemma}
\label{lemma: dim Fc = dim Ffc}
Let $\seq{c}$ be a finite sequence. Then $\dim_H F(f^n(\seq{c})) = \dim_H F(\seq{c})$ for all $n\geq 0$.
\end{lemma}

\begin{proof}
Let $\seq{x}\in F(f(\seq{c}))\setminus F(\seq{c})$. Then $x$ must at least once contain the pattern $\tilde{\seq{c}}$ (or symmetrically $\seq{c}'$). But as $f(\seq{c}) = \tilde{\seq{c}}(\tilde{\seq{c}})^*$ Corollary \ref{cor: aa*} gives that $\seq{x}$ must end with  $(\tilde{\seq{c}}(\tilde{\seq{c}})^*)^{\infty}$. Hence $F(f(\seq{c}))$ contains only countably more elements than $F(\seq{c})$. To prove that we may extend to $f^2$ replace $\seq{c}$ by $\seq{c}_1 = f(\seq{c})$ and repeat the argumentation.
\end{proof}

\vspace{1.5ex}
In the same spirit as in the proof of Lemma \ref{lemma: dim Fc = dim Ffc} Allouche and Cosnard showed that the Thue-Morse sequence is related to the threshold for $F(\seq{c})$ being countable.

\begin{theorem}[Allouche, Cosnard  \cite{allouche83:1,allouche83:2}]
\label{thm: allouche Fc countable}
The set $F(\seq{c})$ is countable if and only if $\seq{c} > \sigma(\seq{t}')$, where $\seq{t}$ is the Thue-Morse sequence.
\end{theorem}

\vspace{1.5ex}
In \cite{moreira:2001} Moreira improved Theorem \ref{thm: allouche Fc countable} and showed that the sequence $\sigma(\seq{t}')$ also is the threshold for the dimension of $F(\seq{c})$.

\begin{theorem}[Moreira \cite{moreira:2001}]
\label{thm: moreira}
The Hausdorff dimension of $F(\seq{c})$ is zero if and only if $\seq{c}\geq \sigma(\seq{t}')$, where $\seq{t}$ is the Thue-Morse sequence.
\end{theorem}

Next, we give a lemma on the distribution of finite shift-bounded sequences.

\begin{lemma}
Let $\seq{s}$ be a finite shift-bounded sequence such that there exists no sequence $\seq{u}$ such that $f(\seq{u}) = \seq{s}$. Then the sequences $\{\seq{u}_k = f^k(\seq{s}): k\geq0 \}$ are the only shift-bounded sequences in the interval $\big(d(\seq{s}), \seq{s}^{\infty}\big]$.
\end{lemma}

\begin{proof}
Let $\seq{v}$ be a shift-bounded sequence in the interval $\big(d(\seq{s}), \seq{s}\big]$. Then there is a $k\geq0$ such that 
\begin{equation}
\label{eq: u_k+1 < v < u_k}
\seq{u}_{k+1} < \seq{v} \leq \seq{u}_k.
\end{equation}
Hence $\seq{u}_k[1,|\seq{u}_k|-1] = \seq{v}[1,|\seq{u}_k|-1]$. 
If $\seq{u}_k = \seq{v}[1,|\seq{u}_k|]$ then by (\ref{eq: u_k+1 < v < u_k})  we must have $\seq{u}_k = \seq{v}$. 
For the case $\seq{u}_{k+1}[1,|\seq{u}_k|] = \seq{v}[1,|\seq{u}_k|]$, we have by Lemma \ref{lemma: s = aa*} and (\ref{eq: u_k+1 < v < u_k}) that $\seq{u}_{k+1}$ can not be a prefix of $\seq{v}$. 
Hence there is a first position $|\seq{u}_k|< i \leq |\seq{u}_{k+1}|$ where $\seq{u}_{k+1}$ and $\seq{v}$ differ. But then $\sigma^{|\seq{u}_k|}(\seq{v}) > \seq{v}'$, contradicting $\seq{v}$ being shift-bounded. 

If $\seq{v}$ is a shift-bounded sequence in the interval $(\seq{s},\seq{s}^{\infty})$ then we must have $\seq{v} = \seq{s}^k\seq{b}$ where $\seq{b}>\seq{s}$. But then $\seq{v}> \seq{s}^{\infty}$, a contradiction.
\end{proof}

\vspace{1.5ex}

For the next definition recall that for a finite sequence $\seq{u}$ we use the notation $\tilde{\seq{u}}$ for the sequence $\seq{u}$ where the last symbol has been inverted.   

\begin{definition}
For a finite shift-bounded sequence $\seq{s}=\seq{u}\,\seq{v}\,\seq{u}^* <1$, where $\seq{u}$ is the longest possible we define the prefix-suffix reduction function $p$ by $p(\seq{s}) = \widetilde{\seq{uv}}$.    
\end{definition}

Note that $\seq{v}$ in the definition above may be empty while $\seq{u}$ is always non-empty as $\seq{s}$ is shift-bounded. The shift-boundedness of $\seq{s}$ in the definition also gives that $p(\seq{s})$ is well defined, that $\frac{1}{2}|\seq{s}|< |p(\seq{s})|<|\seq{s}|$ and that $\seq{s} < p(\seq{s})$.

\begin{lemma}
\label{lemma: p(s) SB}
Let $\seq{s}$ be a finite shift-bounded sequence such that $|\seq{s}|>1$. Then $p(\seq{s})$ is shift-bounded.
\end{lemma}

\begin{proof}
Let $\seq{s}=\seq{u}\seq{v}\seq{u}^*$ where $p(\seq{s}) = \widetilde{\seq{uv}}$. The inequality $\sigma^n(p(\seq{s}))>p(\seq{s})$ for $0<n<|p(\seq{s})|$ follows from the definition of $p$ and that $\seq{s}$ is shift-bounded. For the upper bounding inequality for shift-boundedness we consider first the case when $0<n < |\seq{u}|$. Let
$\seq{\alpha} = p(\seq{s})[1,n]$ and
$\seq{\beta} =  p(\seq{s})[n+1,2n]$.

\footnotesize
\[
\begin{picture}(180,65)
\multiput(40,37)( 0,0){1}{
\put(-35,  0){\makebox(30,12)[r]{$p(\seq{s})'=$}}
\multiput(0,0)(0,0){1}{
\put(   0,  0){\line(1,0){ 80}}\put(   0,12){\line(1,0){ 80}}
\put(   0,  0){\line(0,1){ 12}}\put(  80, 0){\line(0,1){ 12}}}
\multiput(81,0)(0,0){1}{
\put(   0,  0){\line(1,0){ 50}}\put(   0,12){\line(1,0){ 50}}
\put(   0,  0){\line(0,1){ 12}}\put(  50, 0){\line(0,1){ 12}}}
\put(  12,  4){$\seq{\alpha}^*$}
\put(  42,  4){$\seq{\beta}^*$}
\put(   0,  0){$\overbrace{\makebox(80,15)[c]{}}_{}$}
\put(  38, 24){$\seq{u}^*$}
\put( 102,  4){$\seq{v}$}
}
\multiput(40,15)( 0,0){1}{
\put( -35,  0){\makebox(30,12)[r]{$\sigma^n(p(\seq{s}))=$}}
\multiput(0,0)(0,0){1}{
\put(   0,  0){\line(1,0){ 50}}\put(   0,12){\line(1,0){ 50}}
\put(   0,  0){\line(0,1){ 12}}\put(  50, 0){\line(0,1){ 12}}}
\multiput(51,0)(0,0){1}{
\put(   0,  0){\line(1,0){ 50}}\put(   0,12){\line(1,0){ 50}}
\put(   0,  0){\line(0,1){ 12}}\put(  50, 0){\line(0,1){ 12}}}
\put(  12,  4){$\seq{\beta}$}
\put(   0, -2){$\underbrace{\makebox(50,1)[c]{}}_{}$}
\put(  12,-17){$\sigma^n(\seq{u})$}
\put(  72,  4){$\tilde{\seq{v}}$}
}
\multiput( 70,15)(0,4){9}{\line(0,1){2}}
\multiput(100,39)(0,4){3}{\line(0,1){2}}
\end{picture}
\]
\normalsize
Then as  $\seq{s}$ is shift-bounded we have by Corollary \ref{cor: b<a*} that $\seq{\alpha}^*>\seq{\beta}$ and therefore $p(\seq{s})'>\sigma^n(p(\seq{s}))$.

For $|\seq{u}|\leq n < |p(\seq{s})|$ note first that $\seq{u}$ is non-empty as a shift-bounded sequence must end with a 1. Let
$\seq{\alpha} = \seq{s}[1,|\seq{u}\seq{v}|-n]$,
$\seq{\beta}  = \seq{s}[n+1,|\seq{u}\seq{v}|]$ and
$\seq{\gamma} = \seq{s}[|\seq{u}\seq{v}|-n+1,|\seq{u}\seq{v}|-n+|\seq{u}|]$.

\footnotesize
\[
\begin{picture}(180,65)
\multiput(40,37)( 0,0){1}{
\put(-35,  0){\makebox(30,12)[r]{$\seq{s}'=$}}
\multiput(0,0)(0,0){1}{
\put(   0,  0){\line(1,0){ 80}}\put(   0,12){\line(1,0){ 80}}
\put(   0,  0){\line(0,1){ 12}}\put(  80, 0){\line(0,1){ 12}}}
\multiput(81,0)(0,0){1}{
\put(   0,  0){\line(1,0){ 50}}\put(   0,12){\line(1,0){ 50}}
\put(   0,  0){\line(0,1){ 12}}\put(  50, 0){\line(0,1){ 12}}}
\put(   0,  0){$\overbrace{\makebox(39,15)[c]{}}_{}$}\put(  16, 24){$\seq{\alpha}^*$}
\put(  28,  4){$\seq{u}^*$}
\put( 102,  4){$\seq{v}'$}
\put(  41,  0){$\overbrace{\makebox(80,15)[c]{}}_{}$}\put( 77, 24){$\seq{\gamma}^*$}
}
\multiput(40,15)( 0,0){1}{
\put( -35,  0){\makebox(30,12)[r]{$\sigma^n(\seq{s})=$}}
\multiput( 0,0)(0,0){1}{
\put(   0,  0){\line(1,0){ 40}}\put(   0,12){\line(1,0){ 40}}
\put(   0,  0){\line(0,1){ 12}}\put(  40, 0){\line(0,1){ 12}}}
\put(  41,  0){\dashbox(80,12)}
\put(   0, -2){$\underbrace{\makebox(39,1)[c]{}}_{}$}\put(  16,-17){$\seq{\beta}$}
\put(  77,  4){$\seq{u}^*$}
}
\multiput( 80,15)(0,4){9}{\line(0,1){2}}
\multiput(161,15)(0,4){9}{\line(0,1){2}}
\end{picture}
\]
\normalsize
By the definition of $\seq{u}$ we have $\seq{\alpha}^*\seq{\gamma}^*>\seq{\beta}\seq{u}^*$. But as $\seq{u}^* \geq \seq{\gamma}^*$ we must have $\seq{\alpha}^*>\seq{\beta}$ and hence $p(\seq{s})'>\sigma^n(p(\seq{s}))$.
\end{proof}

\vspace{1.5ex}
In \cite{nilsson} the following lemma was given

\begin{lemma}[Nilsson \cite{nilsson}]
\label{lemma: one-sided prefix minimal}
For a sequence $\seq{s}$ let $n_{\seq{s}} = \inf\{n\in\mathbb{N}: \seq{s}[1,n]^{\infty}\geq \seq{s} \}$. Then $\seq{s}[1,m]^{\infty}<\seq{s}[1,n_{\seq{s}}]$ for $m<n_{\seq{s}}$.
\end{lemma}

\begin{lemma}
\label{lemma: shift - minimal}
Let $\seq{s}$ be any non-empty sequence. Then $\sigma^n(\seq{s})>\seq{s}$ for $0<n<N$ if and only if $\seq{s}[1,n]^{\infty}<\seq{s}$ for $0<n<N$.
\end{lemma}

\begin{proof}
Assume that $\sigma^n(\seq{s})>\seq{s}$ for $0<n<N$ and that there is a smallest $0<m<N$ such that $\seq{s}[1,m]^{\infty} \geq \seq{s}$. We cannot have equality, as it would imply that $\seq{s}$ is periodic. Hence $\seq{s}[1,m]^{\infty} > \seq{s}$. Let $\seq{a} = \seq{s}[1,m]$. Then we can for some maximal $k$ write $\seq{s} = \seq{a}^k\seq{b}$, for some sequence $\seq{b}$ with $\seq{b}<\seq{a}$. This implies $\sigma^m(\seq{s}) = \seq{a}^{k-1}\seq{b} < \seq{a}^k\seq{b}=\seq{s}$, a contradiction.

Conversely, assume that $\seq{s}[1,n]^{\infty}<\seq{s}$ for $0<n<N$ and that there is a smallest $0<m<N$ such that $\sigma^m(\seq{s}) \leq \seq{s}$. Again we can out-rule the case of equality as it would imply periodicity. Hence $\sigma^m(\seq{s}) < \seq{s}$. Let $\seq{a} = \seq{s}[1,m]$. Then we can for some maximal $k$ write $\seq{s} = \seq{a}^k\seq{b}$, for some sequence $\seq{b}$ with $\seq{b}<\seq{a}$. This implies $\seq{s}[1,m]^{\infty}=\seq{a}^{\infty}> \seq{a}^k\seq{b} = \seq{s}$, a contradiction. 
\end{proof}

\begin{lemma}
\label{lemma: SB prefix of infinite SB}
Let $\seq{s}$ be an infinite shift-bounded sequence. Then there exists a strictly increasing infinite sequence of integers $\{n_k\}$ such that $\seq{s}[1,n_k]$ is a finite shift-bounded sequence for all $k\geq1$.
\end{lemma}

\begin{proof}
It is clear that for any integer $t$ we have that $(\seq{s}[1,t])' > \sigma^k(\seq{s}[1,t])$ for $0<k<t$. Hence we only have to consider the lower inequality in the definition of shift-boundedness.

There is a $k$ such that $0^k1$ is a prefix of $\seq{s}$. Hence we may put $n_1 = k+1$. Assume for contradiction that there are only $N$ finite shift-bounded prefixes of $\seq{s}$. The shift-boundedness of $\seq{s}$ and Lemma \ref{lemma: shift - minimal} gives $\seq{s}[1,k]^{\infty}<\seq{s}$ for all $k\geq1$. 
For any $m> n_N$ there exists a smallest $n \leq m$ such that $\seq{s}[1,n]^{\infty} \geq \seq{s}[1,m]$, (as this inequality holds for $n=m$). 
The shift-boundedness of $\seq{s}[1,n_N]$ and Lemma \ref{lemma: shift - minimal} implies that $\seq{s}[1,n]^{\infty}<\seq{s}[1,n_N]<\seq{s}[1,m]$ for all $0<n<n_N$. Hence we must have $n_N\leq n$. 
If $n = n_N$ for all $m$ then we obtain $\seq{s}[1,n_N]^{\infty}\geq \seq{s}$, which contradicts the shift-boundedness of $\seq{s}$. Hence there must be an $m$ and smallest $n$ such that $n_N<n\leq m$ with $\seq{s}[1,n]^{\infty} > \seq{s}[1,m]$. Hence Lemma \ref{lemma: one-sided prefix minimal} gives that $\seq{s}[1,k]^{\infty} < \seq{s}[1,n]$ for $0<k<n$, but then Lemma \ref{lemma: shift - minimal} gives that $\sigma^k(\seq{s}[1,n])>\seq{s}[1,n]$ for $0<k<n$, a contradiction to the maximality to $N$.
\end{proof}

\begin{lemma}
\label{lemma: S measure 0}
Let \textit{ISB} be the set of all infinite shift-bounded sequences. Then \textit{ISB} has Lebesque measure zero.
\end{lemma}

\begin{proof}
Let $\lambda$ be the Lebesgue measure. As $\lambda$ is invariant under $x\mapsto 2x$ on the unit circle $\lambda$-almost every $x$ has a dense orbit.
Hence as \textit{ISB} is a set of $x$'s with bounded orbit it must have Lebesgue measure 0.
\end{proof}

\section{Minimal Sequences}

\begin{definition}
For a finite sequence $\seq{s}$ ending with a 1, we define the function $e$ by
$e(\seq{s})  = \tilde{\seq{s}}\, (\seq{s}^*)^{\infty}$. For an infinite sequence $\seq{s}$ we let $e(\seq{s})= \seq{s}$. 
\end{definition}

The motivation for the definition of the function $e$ comes from the special kind of sequences given in (\ref{eq: uuu}) where the upper bounding $k$ has been removed.

\begin{definition}
We denote $\seq{e}_i= e(f^i(1))$ for $i\geq0$ and we say that a sequence $\seq{s}$ not containing only zeros,  finite or infinite, is an $e_i$-sequence if $\seq{e}_{i-1} \leq \seq{s} < \seq{e}_i$.
\end{definition}

Note that $\seq{e}_i$ grows monotonically to $d(1)$ as $i$ tends to infinity. Moreover, for any $i\geq 1$ we have $\seq{e}_{i-1}[1,2^i+2^{i-1}-1] = \seq{e}_{i}[1,2^i+2^{i-1}-1]$. Combining this equality and the fact that $(\seq{e}_{i})_{2^i} = 0$ we see that there are no shift-bounded $e_i$-sequences of length precisely $2^i$, (as any shift-bounded sequence must end with a 1). 
From Lemma \ref{thm: allouche Fc countable} we see that we only have to consider $F(\seq{c})$ for an $e_i$-sequence $\seq{c}$.

\begin{lemma}
\label{lemma: pc = f}
Let $\seq{c}$ be a finite shift-bounded $e_i$-sequence for $i\geq1$. Then there exists an $n>0$ such that $p^n(\seq{c}) = f^i(1)$. 
\end{lemma}

\begin{proof}
Let us use the notation $\seq{w}_k = f^{k}(1)$ for $k\geq0$. Then $|\seq{w}_k| = 2^k$. Assume for contradiction that there is an $n$ such that 
\begin{equation}
\label{eq: pn < w < pn+1}
p^n(\seq{c}) < \seq{w}_i < p^{n+1}(\seq{c}).
\end{equation}
We claim that the above assumption gives the following chain of inequalities
\begin{equation}
\label{eq: e<pn<e<pn+1}
\seq{e}_{i-1}< p^n(\seq{c}) <  \seq{e}_{i} < p^{n+1}(\seq{c}).
\end{equation}
The left-most inequality of (\ref{eq: e<pn<e<pn+1}) is clear as $\seq{c}$ is an $e_i$-sequence and therefore $\seq{e}_{i-1} < \seq{c} < p^n(\seq{c})$. The right-most inequality of (\ref{eq: e<pn<e<pn+1}) is given by our assumption, $\seq{e}_{i} < \seq{w}_i < p^{n+1}(\seq{c})$.
For the middle inequality of the claim (\ref{eq: e<pn<e<pn+1}), assume that $p^n(\seq{c})>\seq{e}_i$. Then $p^n(\seq{c}) = \tilde{\seq{w}}_i\,\seq{s}$ for some non-empty sequence $\seq{s}$ with $|\seq{s}|\leq |\seq{w}_i|$ and $\seq{s}>\seq{w}_i^*$. (If $\seq{s}$ were empty then $p^n(\seq{c}) < \seq{e}_{i-1}$ since $\seq{e}_{i-1}[1,2^i] = \seq{e}_i[1,2^i]$). Thus
\[
p^n(\seq{c})' = \seq{w}_i^*\seq{s}' < \sigma^{|\seq{w}_i|}(p^n(\seq{c})) = \seq{s},
\] 
contradicting the shift-boundedness of $p^n(\seq{c})$ and concludes the claim. 

The assumption (\ref{eq: pn < w < pn+1}) also gives that $|p^{n+1}(\seq{c})|<|\seq{w}_i|$, since otherwise $\seq{w}_i$ would be a prefix of $p^n(\seq{c})$. From the equality $\seq{e}_{i-1}[1,2^i] = \seq{e}_i[1,2^i]$ and (\ref{eq: e<pn<e<pn+1}) we also have $|\seq{w}_i|<|p^n(\seq{c})|$.

Let $p^n(\seq{c}) = \seq{u}\seq{v}\seq{u}^*$ where $p^{n+1}(\seq{c}) = \widetilde{\seq{u}\seq{v}}$. 
Put $\seq{\alpha} = p^n(\seq{c})[1,2^i-|\seq{u}\seq{v}|]$, 
$\seq{\beta}^*  = p^n(\seq{c})[|\seq{u}\seq{v}|+1,2^i]$ and 
$\seq{\gamma}^* = \seq{e}_{i-1}[|\seq{uv}|+1,2^i]$.

\footnotesize
\[
\begin{picture}(270,65)
\multiput(40,37)( 0,0){1}{
\put( -35,  0){\makebox(30,12)[r]{$\seq{e}_{i-1}=$}}
\multiput(0,0)( 81,0){2}{
\put(   0,  0){\line(1,0){ 80}}\put(  0,12){\line(1,0){ 80}}
\put(   0,  0){\line(0,1){ 12}}\put( 80, 0){\line(0,1){ 12}}}
\multiput(162,0)(0,0){1}{
\put(   0,  0){\line(1,0){ 71}}\put(  0,12){\line(1,0){ 71}}
\put(   0,  0){\line(0,1){ 12}}}
\put(  81,  0){$\overbrace{\makebox(78,15)[c]{}}_{}$}\put(116,24){$\seq{w}_{i-1}^*$}
\put(  34,  4){$\tilde{\seq{w}}_{i-1}$}
\put( 142,  4){$\seq{\gamma}^*$}
\put( 190,  4){$\tilde{\seq{w}}_{i-1}$}
}
\multiput(40,15)( 0,0){1}{
\put(-35,  0){\makebox(30,12)[r]{$p^n(\seq{c})=$}}
\multiput(0,0)(0,0){1}{
\put(   0,  0){\line(1,0){ 60}}\put(   0,12){\line(1,0){ 60}}
\put(   0,  0){\line(0,1){ 12}}\put(  60, 0){\line(0,1){ 12}}}
\multiput(61,0)(0,0){1}{
\put(   0,  0){\line(1,0){ 70}}\put(   0,12){\line(1,0){ 70}}
\put(   0,  0){\line(0,1){ 12}}\put(  70, 0){\line(0,1){ 12}}}
\multiput(132,0)(0,0){1}{
\put(   0,  0){\line(1,0){ 60}}\put(   0,12){\line(1,0){ 60}}
\put(   0,  0){\line(0,1){ 12}}\put(  60, 0){\line(0,1){ 12}}}
\put(   1, -2){$\underbrace{\makebox(58,1)[c]{}}_{}$}\put( 26,-17){$\seq{u}$}
\put( 133, -2){$\underbrace{\makebox(58,1)[c]{}}_{}$}\put(158,-17){$\seq{u}^*$}
\multiput(30,0)(0,4){3}{\line(0,1){2}}
\put(  12,  4){$\seq{\alpha}$}
\put(  93,  4){$\seq{v}$}
\put( 142,  4){$\seq{\beta}^*$}
}
\multiput(201,15)(0,4){9}{\line(0,1){2}}
\multiput(171,15)(0,4){9}{\line(0,1){2}}
\end{picture}
\]
\normalsize
By the definition of $p$ we have $\seq{\alpha} = \seq{\beta}$ and as $\seq{w}_i$ is shift-bounded we have $\seq{\alpha}^* > \seq{\gamma}^*$, that is, $\seq{\gamma} \neq \seq{\beta}$, a contradiction.
\end{proof}

\begin{definition}
For $n\in\mathbb{N}$ such that $n\leq|\seq{s}|$ we define the function $g$ by $g_n(\seq{s}) = \seq{s}[1,n-1]1$. If $n = \infty$ we let $g_n(\seq{s}) = \seq{s}$.
\end{definition}

\begin{definition}
\label{def: min-seq}
For an $e_i$-sequence $\seq{s}$ we define the integer $m_{\seq{s}}$ by,
\[	
m_{\seq{s}} = \inf \left\{ n\geq 2^i: e(g_n(\seq{s})) \leq \seq{s} \leq g_n(\seq{s})^{\infty}  \right\}.
\]
If $m_{\seq{s}}$ is undefined we set $m_{\seq{s}}=\infty$ and $g_{m_{\seq{s}}}(\seq{s}) = \seq{s}$. We say that $g_{m_{\seq{s}}}(\seq{s})$ is an $e_i$-minimal prefix of $\seq{s}$. An $e_i$-sequence $\seq{s}$ is a finite $e_i$-minimal sequence if $g_{m_{\seq{s}}}(\seq{s}) = \seq{s}[1,m_{\seq{s}}] = \seq{s}\,$ for $m_{\seq{s}}<\infty$ and $\seq{s}$ is an infinite $e_i$-minimal sequence if $m_{\seq{s}} = \infty$.
\end{definition}

\begin{lemma}
An $e_i$-minimal prefix is an $e_i$-sequence.
\end{lemma}

\begin{proof}
Let $\seq{s}$ be an $e_i$-sequence. The lemma is clear if the $e_i$-minimal prefix of $\seq{s}$ is an infinite sequence. Hence we assume that $\seq{s}$ has the finite $e_i$-minimal prefix $\seq{c}$, that is, $\seq{c} = g_{m_{\seq{s}}}(\seq{s})$.

Assume for contradiction that $\seq{c}>\seq{e}_i$. We must have that $\tilde{\seq{c}}$ is a prefix of $\seq{e}_i$.
Let $\seq{w} = f^i(1)$, 
$\seq{\alpha}^* = e(\seq{c})[|\seq{c}|+1, |\seq{w}|k]$ and 
$\seq{\gamma}^* = \seq{e}_i[|\seq{c}|+1, |\seq{w}|k]$, where $k$ is the smallest integer such that $|\seq{w}|k \geq |\seq{c}|$.

\footnotesize
\[
\begin{picture}(230,65)
\multiput(40,37)( 0,0){1}{
\put(-35,  0){\makebox(30,12)[r]{$\seq{e}_i=\ldots$}}
\multiput(0,0)(0,0){1}{
\put(   0,  0){\line(1,0){ 50}}\put(   0,12){\line(1,0){ 50}}
\put(  50, 0){\line(0,1){ 12}}}
\multiput(51,0)(61,0){2}{
\put(   0,  0){\line(1,0){ 60}}\put(   0,12){\line(1,0){ 60}}
\put(   0,  0){\line(0,1){ 12}}\put(  60, 0){\line(0,1){ 12}}}
\multiput(173,0)(0,0){1}{
\put(   0,  0){\line(1,0){ 20}}\put(   0,12){\line(1,0){ 20}}
\put(   0,  0){\line(0,1){ 12}}}
\put(  52,  0){$\overbrace{\makebox(58,15)[c]{}}_{}$}\put(  78, 24){$\seq{w}^*$}
\put( 113,  0){$\overbrace{\makebox(58,15)[c]{}}_{}$}\put( 139, 24){$\seq{w}^*$}
\put(  87,  4){$\seq{\gamma}^*$}
}
\multiput(40,15)( 0,0){1}{
\put( -35,  0){\makebox(30,12)[r]{$e(\seq{c})=\ldots$}}
\multiput(0,0)(0,0){1}{
\put(   0,  0){\line(1,0){ 70}}\put(   0,12){\line(1,0){ 70}}
\put(  70,  0){\line(0,1){ 12}}}
\multiput(71,0)(81,0){1}{
\put(   0,  0){\line(1,0){122}}\put(   0,12){\line(1,0){122}}
\put(   0,  0){\line(0,1){ 12}}}
\put(  72, -2){$\underbrace{\makebox(58,1)[c]{}}_{}$}\put( 93,-17){$(\tilde{\seq{w}})^*$}
\put(  87,  4){$\seq{\alpha}^*$}
}
\multiput(111,15)(0,4){9}{\line(0,1){2}}
\multiput(151,15)(0,4){9}{\line(0,1){2}}
\multiput(171,15)(0,4){3}{\line(0,1){2}}
\end{picture}
\]
\normalsize
If $|\seq{w}|k>|\seq{c}|$ then as $\seq{w}$ is shift-bounded we have $\seq{\alpha}^*>\seq{\gamma}^*$. Therefore $\seq{s} \geq e(\seq{c})>\seq{e}_i$, contradicting $\seq{s}$ being an $e_i$-sequence. 
For $|\seq{w}|k = |\seq{c}|$ then we reach $\seq{e}_i<e(\seq{c})\leq\seq{s}$ as  $\seq{w}^* < (\tilde{\seq{w}})^*$, which gives a contradiction to our assumption. 

For the second case, assume for contradiction that $\seq{c}<\seq{e}_{i-1}$. We must then have that $\seq{c}$ is a prefix of $\seq{e}_{i-1}$. 
Let $\seq{w} = f^{i-1}(1)$ and note that $|\seq{w}| = 2^{i-1}$. Furthermore let  
$\seq{v}^*            = \seq{e}_{i-1}[|\seq{c}|+1, |\seq{w}|k]$, 
$\seq{u}              = \seq{c}^{\infty}[|\seq{c}|+1, |\seq{w}|k]$, 
$\seq{\alpha}^*       = \seq{e}_i[2^{i-1}k+1, |\seq{c}|+ |\seq{w}|]$ and 
$\tilde{\seq{\gamma}} = \seq{c}^{\infty}[|\seq{w}|k+1, |\seq{c}|+|\seq{w}|]$
where $k$ is the smallest integer such that $|\seq{w}|k \geq |\seq{c}|$. 

\footnotesize
\[
\begin{picture}(230,65)
\multiput(40,37)( 0,0){1}{
\put( -35,  0){\makebox(30,12)[r]{$\seq{c}^{\infty}=\ldots$}}
\multiput(0,0)(0,0){1}{
\put(   0,  0){\line(1,0){110}}\put(   0,12){\line(1,0){110}}
\put( 110, 0){\line(0,1){ 12}}}
\multiput(111,0)(0,0){1}{
\put(   0,  0){\line(1,0){ 82}}\put(   0,12){\line(1,0){ 82}}
\put(   0,  0){\line(0,1){ 12}}}
\put( 111,  0){$\overbrace{\makebox(58,15)[c]{}}_{}$}\put( 137, 24){$\tilde{\seq{w}}$}
\put(  58,  4){$\seq{c}$}
\put( 129,  4){$\seq{u}$}
\put( 160,  4){$\tilde{\seq{\gamma}}$}
}
\multiput(40,15)( 0,0){1}{
\put(-35,  0){\makebox(30,12)[r]{$\seq{e}_{i-1}=\ldots$}}
\multiput(0,0)(0,0){1}{
\put(   0,  0){\line(1,0){ 30}}\put(   0,12){\line(1,0){ 30}}
\put(  30, 0){\line(0,1){ 12}}}
\multiput(31,0)(61,0){2}{
\put(   0,  0){\line(1,0){ 60}}\put(   0,12){\line(1,0){ 60}}
\put(   0,  0){\line(0,1){ 12}}\put(  60, 0){\line(0,1){ 12}}}
\multiput(153,0)(0,0){1}{
\put(   0,  0){\line(1,0){ 40}}\put(   0,12){\line(1,0){ 40}}
\put(   0,  0){\line(0,1){ 12}}}
\put(  58,  4){$\seq{w}^*$}
\put(  93, -2){$\underbrace{\makebox(58,1)[c]{}}_{}$}\put( 119,-17){$\seq{w}^*$}
\put( 129,  4){$\seq{v}^*$}
\put( 158,  4){$\seq{\alpha}^*$}
}
\multiput(151,15)(0,4){9}{\line(0,1){2}}
\multiput(192,15)(0,4){9}{\line(0,1){2}}
\multiput(212,15)(0,4){9}{\line(0,1){2}}
\end{picture}
\]
\normalsize
If $|\seq{w}|k > |\seq{c}|$ then as $\seq{w}$ is shift-bounded we have $\seq{v}^*\seq{\alpha}^*>\seq{u}\tilde{\seq{\gamma}} = \tilde{\seq{w}}$ and therefore 
$\seq{s} \leq  \seq{c}^{\infty} < \seq{e}_{i-1}$, a contradiction. If $|\seq{w}|k = |\seq{c}|$ the result follows as $\tilde{\seq{w}}<\seq{w}^*$, which concludes the proof.
\end{proof}

\begin{lemma}
\label{lemma: c'>s(c)>c n<2i}
Let $\seq{s}$ be a non-zero $e_i$-sequence. Then $\seq{s}'[1,|\seq{s}|-n]>\sigma^n(\seq{s})>\seq{s}$ for all $0< n < 2^i$.
\end{lemma}

\begin{proof}
The result is clear for any $e_1$-sequence. Hence we only have to consider the case with $i>1$. Let $\seq{u}_k = f^k(1)$ and put 
\[
\seq{s}_{i} :=\tilde{\seq{u}}_{i}^{}\,\seq{u}_{i-2}' 
=\tilde{\seq{u}}_{i-2}^{}\,\seq{u}_{i-2}^*\,\seq{u}_{i-2}'\,\tilde{\seq{u}}_{i-2}^{}\,\seq{u}_{i-2}'.
\]
Note that $|\seq{s}_i| = 5\cdot 2^{i-2}$. We have that $\seq{s}_i$ is a prefix of all $e_i$-sequences and moreover $\seq{s}_i$ is a prefix of $f^{i+1}(1)$. To prove the lemma it is enough to show that $\seq{s}_i'[1,|\seq{s}_i|-n]>\sigma^n(\seq{s}_i)>\seq{s}_i$ holds for $0< n < 2^i$, as $\seq{s}_i$ is a prefix of any $e_i$-sequence. Since $f^{i+1}(1)$ is a shift-bounded sequence we have that 
\begin{equation}
\label{eq: s'>s>s}
\seq{s}_i'[1,|\seq{s}_i|-n] \geq\sigma^n(\seq{s}_i) > \seq{s}_i[1,|\seq{s}_i|-n]
\end{equation}
holds for $0<n<|\seq{s}_i|$. Hence we have to show that these (\ref{eq: s'>s>s}) shift-inequalities are strict for $0<n<2^i$.
From the definition of $\seq{s}_i$ we have directly that $\seq{s}_2 = 00101$ and by a straight forward calculation we see that lemma holds in this case. Hence we may assume that $i\geq3$.
 
For the upper bounding inequality in (\ref{eq: s'>s>s}) let 
$\seq{\alpha}        = \seq{s}_i[|\seq{u}_i|+1,|\seq{u}_i|+n]$ and
$\tilde{\seq{\gamma}}= \seq{s}_i[|\seq{u}_i|-n+1,|\seq{u}_i|]$.

\footnotesize
\[
\begin{picture}(270,65)
\multiput(40,37)( 0,0){1}{
\put(-35,  0){\makebox(30,12)[r]{$\seq{s}_i'=$}}
\multiput(0,0)(46,0){1}{
\put(   0,  0){\line(1,0){183}}\put(   0,12){\line(1,0){183}}
\put(   0,  0){\line(0,1){ 12}}\put( 183, 0){\line(0,1){ 12}}}
\multiput(184,0)(46,0){1}{
\put(   0,  0){\line(1,0){ 45}}\put(   0,12){\line(1,0){ 45}}
\put(   0,  0){\line(0,1){ 12}}\put(  45, 0){\line(0,1){ 12}}}
\put(   1,  0){$\overbrace{\makebox(182,15)[c]{}}_{}$}\put(  88, 24){$\seq{u}'_i$}
\put(  52,  4){$\seq{\alpha}^*$}
\put( 201,  4){$\seq{u}'_{i-2}$}
}
\multiput(40,15)( 0,0){1}{
\put( -35,  0){\makebox(30,12)[r]{$\sigma^n(\seq{s}_i)=$}}
\multiput(0,0)(0,0){1}{
\put(   0,  0){\line(1,0){112}}\put(   0,12){\line(1,0){112}}
\put(   0,  0){\line(0,1){ 12}}\put( 112, 0){\line(0,1){ 12}}}
\multiput(113,0)(0,0){1}{
\put(   0,  0){\line(1,0){ 45}}\put(   0,12){\line(1,0){ 45}}
\put(   0,  0){\line(0,1){ 12}}\put(  45, 0){\line(0,1){ 12}}}
\put(   1,  -2){$\underbrace{\makebox(110,1)[c]{}}_{}$}\put( 45,-17){$\sigma^n(\tilde{\seq{u}}_i)$}
\put(  52,  4){$\tilde{\seq{\gamma}}$}
\put( 131,  4){$\seq{u}'_{i-2}$}
}
\multiput(152,15)(0,4){9}{\line(0,1){2}}
\end{picture}
\]
\normalsize
As $\seq{u}_i$ is shift-bounded we have $\seq{\alpha}^*>\tilde{\seq{\gamma}}$ and therefore $\seq{s}_i'[1,|\seq{s}_i|-n]>\sigma^n(\seq{s}_i)$. 

To prove the lower inequality of (\ref{eq: s'>s>s}) we consider first the case when $0 < n < |\seq{u}_{i-2}|$. Let 
$\seq{\alpha}^*      = \seq{s}_i[|\seq{u}_{i-2}|+1,|\seq{u}_{i-2}|+n]$ and
$\tilde{\seq{\gamma}}= \seq{s}_i[|\seq{u}_{i-2}|-n+1,|\seq{u}_{i-2}|]$.

\footnotesize
\[
\begin{picture}(270,65)
\multiput(40,37)( 0,0){1}{
\put( -35,  0){\makebox(30,12)[r]{$\sigma^n(\seq{s}_i)=$}}
\multiput(21,0)(46,0){4}{
\put(   0,  0){\line(1,0){ 45}}\put(   0,12){\line(1,0){ 45}}
\put(   0,  0){\line(0,1){ 12}}\put(  45, 0){\line(0,1){ 12}}}
\multiput(0,0)(0,0){1}{
\put(   0,  0){\line(1,0){ 20}}\put(   0,12){\line(1,0){ 20}}
\put(   0,  0){\line(0,1){ 12}}\put(  20, 0){\line(0,1){ 12}}}
\put(  22,  0){$\overbrace{\makebox(43,15)[c]{}}_{}$}\put(  41, 24){$\seq{u}^*_{i-2}$}
\put(  30,  4){$\seq{\alpha}^*$}
\put(  85,  4){$\seq{u}'_{i-2}$}
\put( 131,  4){$\tilde{\seq{u}}_{i-2}$}
\put( 177,  4){$\seq{u}'_{i-2}$}
}
\multiput(40,15)( 0,0){1}{
\put(-35,  0){\makebox(30,12)[r]{$\seq{s}_i=$}}
\multiput(0,0)(46,0){5}{
\put(   0,  0){\line(1,0){ 45}}\put(   0,12){\line(1,0){ 45}}
\put(   0,  0){\line(0,1){ 12}}\put(  45, 0){\line(0,1){ 12}}}
\put(  1,  -2){$\underbrace{\makebox(43,1)[c]{}}_{}$}\put( 20,-17){$\tilde{\seq{u}}_{i-2}$}
\put(  30,  4){$\tilde{\seq{\gamma}}$}
\put(  63,  4){$\seq{u}^*_{i-2}$}
\put( 109,  4){$\seq{u}'_{i-2}$}
\put( 155,  4){$\tilde{\seq{u}}_{i-2}$}
\put( 201,  4){$\seq{u}'_{i-2}$}
}
\multiput( 61,15)(0,4){9}{\line(0,1){2}}
\multiput( 85,15)(0,4){9}{\line(0,1){2}}
\end{picture}
\]
\normalsize
As $\seq{u}_{i-2}$ is shift-bounded we have $\seq{\alpha}^*>\tilde{\seq{\gamma}}$ and therefore 
$\sigma^n(\seq{s}_i) > \seq{s}_i$. 
The case $n = |\seq{u}_{i-2}|$ is clear as $\tilde{\seq{u}}_{i-2}< \seq{u}_{i-2}^*$. 

For $|\seq{u}_{i-2}| < n < 2|\seq{u}_{i-2}|$ let
$\seq{\alpha}^*      = \seq{s}_i[2|\seq{u}_{i-2}|+1,|\seq{u}_{i-2}|+n]$ and
$\tilde{\seq{\gamma}}= \seq{s}_i[|\seq{u}_{i-2}|-n+1,|\seq{u}_{i-2}|]$.

\footnotesize
\[
\begin{picture}(270,65)
\multiput(40,37)( 0,0){1}{
\put( -35,  0){\makebox(30,12)[r]{$\sigma^n(\seq{s}_i)=$}}
\multiput(21,0)(46,0){3}{
\put(   0,  0){\line(1,0){ 45}}\put(   0,12){\line(1,0){ 45}}
\put(   0,  0){\line(0,1){ 12}}\put(  45, 0){\line(0,1){ 12}}}
\multiput(0,0)(0,0){1}{
\put(   0,  0){\line(1,0){ 20}}\put(   0,12){\line(1,0){ 20}}
\put(   0,  0){\line(0,1){ 12}}\put(  20, 0){\line(0,1){ 12}}}
\put(  22,  0){$\overbrace{\makebox(43,15)[c]{}}_{}$}\put(  41, 24){$\seq{u}'_{i-2}$}
\put(  30,  4){$\seq{\alpha}^*$}
\put(  85,  4){$\tilde{\seq{u}}_{i-2}$}
\put( 131,  4){$\seq{u}'_{i-2}$}
}
\multiput(40,15)( 0,0){1}{
\put(-35,  0){\makebox(30,12)[r]{$\seq{s}_i=$}}
\multiput(0,0)(46,0){5}{
\put(   0,  0){\line(1,0){ 45}}\put(   0,12){\line(1,0){ 45}}
\put(   0,  0){\line(0,1){ 12}}\put(  45, 0){\line(0,1){ 12}}}
\put(  1,  -2){$\underbrace{\makebox(43,1)[c]{}}_{}$}\put( 20,-17){$\tilde{\seq{u}}_{i-2}$}
\put(  30,  4){$\tilde{\seq{\gamma}}$}
\put(  63,  4){$\seq{u}^*_{i-2}$}
\put( 109,  4){$\seq{u}'_{i-2}$}
\put( 155,  4){$\tilde{\seq{u}}_{i-2}$}
\put( 201,  4){$\seq{u}'_{i-2}$}
}
\multiput( 61,15)(0,4){9}{\line(0,1){2}}
\multiput( 85,15)(0,4){9}{\line(0,1){2}}
\end{picture}
\]
\normalsize
Again by the shift-boundedness of $\seq{u}_{i-2}$  we have $\seq{\alpha}^*>\tilde{\seq{\gamma}}$ and therefore 
$\sigma^n(\seq{s}_i) > \seq{s}_i$.
The case $n = 2|\seq{u}_{i-2}|$ follows as $\tilde{\seq{u}}_{i-2}<\seq{u}_{i-2}'$.

For $2|\seq{u}_{i-2}| < n < 3|\seq{u}_{i-2}|$ let
$\seq{\alpha}          = \seq{s}_i[1,3|\seq{u}_{i-2}|-n]$ and
$(\tilde{\seq{\gamma}})^*= \seq{s}_i[n+1,3|\seq{u}_{i-2}|]$.

\footnotesize
\[
\begin{picture}(270,70)
\multiput(40,37)( 0,0){1}{
\put( -35,  0){\makebox(30,12)[r]{$\sigma^n(\seq{s}_i)=$}}
\multiput(0,0)(0,0){1}{
\put(   0,  0){\line(1,0){ 30}}\put(   0,12){\line(1,0){ 30}}
\put(   0,  0){\line(0,1){ 12}}\put(  30, 0){\line(0,1){ 12}}}
\multiput(31,0)(46,0){2}{
\put(   0,  0){\line(1,0){ 45}}\put(   0,12){\line(1,0){ 45}}
\put(   0,  0){\line(0,1){ 12}}\put(  45, 0){\line(0,1){ 12}}}
\put(   1,  0){$\overbrace{\makebox(28,15)[c]{}}_{}$}\put(  -8, 24){$\sigma^{n-2|\seq{u}_{i-2}|}(\seq{u}'_{i-2})$}
\put(   8,  4){$(\tilde{\seq{\gamma}})^*$}
\put(  49,  4){$\tilde{\seq{u}}_{i-2}$}
\put(  95,  4){$\seq{u}'_{i-2}$}
}
\multiput(40,15)( 0,0){1}{
\put(-35,  0){\makebox(30,12)[r]{$\seq{s}_i=$}}
\multiput(0,0)(46,0){5}{
\put(   0,  0){\line(1,0){ 45}}\put(   0,12){\line(1,0){ 45}}
\put(   0,  0){\line(0,1){ 12}}\put(  45, 0){\line(0,1){ 12}}}
\put(  1,  -2){$\underbrace{\makebox(43,1)[c]{}}_{}$}\put( 20,-17){$\tilde{\seq{u}}_{i-2}$}
\put(  12,  4){$\seq{\alpha}$}
\put(  63,  4){$\seq{u}^*_{i-2}$}
\put( 109,  4){$\seq{u}'_{i-2}$}
\put( 155,  4){$\tilde{\seq{u}}_{i-2}$}
\put( 201,  4){$\seq{u}'_{i-2}$}
}
\multiput( 71,15)(0,4){9}{\line(0,1){2}}
\end{picture}
\]
\normalsize
The shift-boundedness of $\seq{u}_{i-2}$ gives again $\seq{\alpha}<(\tilde{\seq{\gamma}})^*$ and therefore 
$\sigma^n(\seq{s}_i) > \seq{s}_i$.
The case $n = 3|\seq{u}_{i-2}|$ is clear as $\seq{u}_{i-2}^*<\seq{u}_{i-2}'$.

For $3|\seq{u}_{i-2}| < n < 4|\seq{u}_{i-2}|$ let
$\seq{\alpha}^*      = \seq{s}_i[4|\seq{u}_{i-2}|+1,|\seq{u}_{i-2}|+n]$ and
$\tilde{\seq{\gamma}}= \seq{s}_i[4|\seq{u}_{i-2}|-n+1,|\seq{u}_{i-2}|]$.

\footnotesize
\[
\begin{picture}(270,65)
\multiput(40,37)( 0,0){1}{
\put( -35,  0){\makebox(30,12)[r]{$\sigma^n(\seq{s}_i)=$}}
\multiput(21,0)(46,0){1}{
\put(   0,  0){\line(1,0){ 45}}\put(   0,12){\line(1,0){ 45}}
\put(   0,  0){\line(0,1){ 12}}\put(  45, 0){\line(0,1){ 12}}}
\multiput(0,0)(0,0){1}{
\put(   0,  0){\line(1,0){ 20}}\put(   0,12){\line(1,0){ 20}}
\put(   0,  0){\line(0,1){ 12}}\put(  20, 0){\line(0,1){ 12}}}
\put(  22,  0){$\overbrace{\makebox(43,15)[c]{}}_{}$}\put(  41, 24){$\seq{u}'_{i-2}$}
\put(  30,  4){$\seq{\alpha}^*$}
}
\multiput(40,15)( 0,0){1}{
\put(-35,  0){\makebox(30,12)[r]{$\seq{s}_i=$}}
\multiput(0,0)(46,0){5}{
\put(   0,  0){\line(1,0){ 45}}\put(   0,12){\line(1,0){ 45}}
\put(   0,  0){\line(0,1){ 12}}\put(  45, 0){\line(0,1){ 12}}}
\put(  1,  -2){$\underbrace{\makebox(43,1)[c]{}}_{}$}\put( 20,-17){$\tilde{\seq{u}}_{i-2}$}
\put(  30,  4){$\tilde{\seq{\gamma}}$}
\put(  63,  4){$\seq{u}^*_{i-2}$}
\put( 109,  4){$\seq{u}'_{i-2}$}
\put( 155,  4){$\tilde{\seq{u}}_{i-2}$}
\put( 201,  4){$\seq{u}'_{i-2}$}
}
\multiput( 61,15)(0,4){9}{\line(0,1){2}}
\multiput( 85,15)(0,4){9}{\line(0,1){2}}
\end{picture}
\]
\normalsize
As $\seq{u}_{i-2}$ is shift-bounded we have $\seq{\alpha}<(\tilde{\seq{\gamma}})^*$ and therefore 
$\sigma^n(\seq{s}_i) > \seq{s}_i$. 
The case $n = 4|\seq{u}_{i-2}|$ is as before clear as $\tilde{\seq{u}}_{i-2}<\seq{u}_{i-2}'$, concluding the proof of the lower inequality of (\ref{eq: s'>s>s}). 
\end{proof}

\begin{lemma}
\label{lemma: prefix is minimal}
An $e_i$-minimal prefix is an $e_i$-minimal sequence.
\end{lemma}

\begin{proof}
It is clear that the statement holds in the case when the $e_i$-minimal prefix is an infinite sequence. Let $\seq{c}$ be the finite $e_i$-minimal prefix of the sequence $\seq{s}$, \mbox{i.e.} $\seq{c} = g_{m_{\seq{s}}}(\seq{s})$. We have to show that the $e_i$-minimal prefix of $\seq{c}$ is $\seq{c}$ itself, that is, $\seq{c} = g_{m_{\seq{c}}}(\seq{c})$. Assume for contradiction that $m_{\seq{c}}<m_{\seq{s}}$. 

If $c_{m_{\seq{c}}} = 0$ then by definition of an $e_i$-minimal prefix we have
\[
e(g_{m_{\seq{c}}}(\seq{c})) < \seq{c}[1,m_{\seq{c}}] < g_{m_{\seq{c}}}(\seq{c})^{\infty}, 
\]
but this is a contradiction, as $\seq{c}[1,m_{\seq{c}}]$ is a proper prefix of $e(g_{m_{\seq{c}}}(\seq{c}))$.

For the case $c_{m_{\seq{c}}} = 1$, consider first the case when $s_{m_{\seq{s}}} = 0$. Let 
$\seq{\gamma}_{\seq{s}} = \seq{s}[m_{\seq{c}}+1,m_{\seq{s}}]$,
$\seq{\gamma}_{\seq{c}} = \seq{c}[m_{\seq{c}}+1,m_{\seq{s}}]$ and 
$\seq{\alpha} = g_{m_{\seq{c}}}(\seq{c})^{\infty}[m_{\seq{c}}+1,m_{\seq{s}}]$. 

\footnotesize
\[
\begin{picture}(265,90)
\multiput(40,67)( 0,0){1}{
\put( -35,  0){\makebox(30,12)[r]{$\seq{s}=$}}
\multiput(0,0)(0,0){1}{
\put(   0,  0){\line(1,0){ 50}}\put(   0,12){\line(1,0){ 50}}
\put(   0,  0){\line(0,1){ 12}}\put(  50, 0){\line(0,1){ 12}}}
\multiput(51,0)(0,0){1}{
\put(   0,  0){\line(1,0){133}}\put(   0,12){\line(1,0){133}}
\put(   0,  0){\line(0,1){ 12}}\put( 133, 0){\line(0,1){ 12}}}
\multiput(185,0)(0,0){1}{
\put(   0,  0){\line(1,0){ 40}}\put(   0,12){\line(1,0){ 40}}
\put(   0,  0){\line(0,1){ 12}}}
\put( 112,  4){$\seq{\gamma}_{\seq{s}}$}
\put(  48, 20){$m_{\seq{c}}$}
\put( 182, 20){$m_{\seq{s}}$}
\put( 177,  4){$0$}
}
\multiput(40,41)( 0,0){1}{
\put(-35,  0){\makebox(30,12)[r]{$\seq{c}=$}}
\multiput(0,0)(0,0){1}{
\put(   0,  0){\line(1,0){ 50}}\put(   0,12){\line(1,0){ 50}}
\put(   0,  0){\line(0,1){ 12}}\put(  50, 0){\line(0,1){ 12}}}
\put(  43,  4){$1$}
\multiput(51,0)(0,0){1}{
\put(   0,  0){\line(1,0){133}}\put(   0,12){\line(1,0){133}}
\put(   0,  0){\line(0,1){ 12}}\put( 133, 0){\line(0,1){ 12}}}
\put( 112,  4){$\seq{\gamma}_{\seq{c}}$}
\put( 177,  4){$1$}
}
\multiput(40,15)( 0,0){1}{
\put(-35,  0){\makebox(30,12)[r]{$g_{m_{\seq{c}}}(\seq{c})^{\infty}=$}}
\multiput(0,0)(51,0){4}{
\put(   0,  0){\line(1,0){ 50}}\put(   0,12){\line(1,0){ 50}}
\put(   0,  0){\line(0,1){ 12}}\put(  50, 0){\line(0,1){ 12}}}
\multiput(204,0)(0,0){1}{
\put(   0,  0){\line(1,0){ 21}}\put(   0,12){\line(1,0){ 21}}
\put(   0,  0){\line(0,1){ 12}}}
\put(  52,  -2){$\underbrace{\makebox(131,1)[c]{}}_{}$}\put(112,-17){$\seq{\alpha}$}
}
\multiput( 90,15)(0,4){16}{\line(0,1){2}}
\multiput(224,15)(0,4){16}{\line(0,1){2}}
\end{picture}
\]
\normalsize
It is clear that $\seq{\gamma}_{\seq{s}}<\seq{\gamma}_{\seq{c}}$.  As $\seq{c}$ is the $e_i$-minimal prefix of $\seq{s}$ we have $\seq{\alpha}\leq \seq{\gamma}_{\seq{s}} $ and as $g_{m_{\seq{c}}}(\seq{c})$ is the $e_i$-minimal prefix of $\seq{c}$ we have $\seq{\alpha} \geq \seq{\gamma}_{c}$. Hence $\seq{\alpha} \leq \seq{\gamma}_{\seq{s}} <\seq{\gamma}_{\seq{c}} \leq \seq{\alpha}$, a contradiction.
 
Finally, let us turn to the case with $s_{m_{\seq{c}}} = 1$. From Lemma \ref{lemma: c'>s(c)>c n<2i} we have that $\sigma^n(\seq{s})>\seq{s}$ for all $0< n < 2^i$ and combining this with Lemma \ref{lemma: shift - minimal} we get that $\seq{s}[1,n]^{\infty}< \seq{s}$ for $0<n<2^i$ and our assumption extends this to that $\seq{s}[1,n]^{\infty}< \seq{s}$ for $0<n<m_{\seq{s}}$. 
Hence $m_{\seq{s}}$ coincides with the integer $n_{\seq{s}} = \inf \{n\in\mathbb{N} : \seq{s}[1,n]^{\infty}\geq \seq{s}\}$, that is $m_{\seq{s}} = n_{\seq{s}}$.
But then Lemma \ref{lemma: one-sided prefix minimal} gives
\[
g_{m_{\seq{c}}}(\seq{c})^{\infty} = \seq{c}[1,m_{\seq{c}}]^{\infty}<\seq{s}[1,m_{\seq{s}}] = \seq{c},
\]
a contradiction to that $g_{m_{\seq{c}}}(\seq{c})$ is the $e_i$-minimal prefix of $\seq{c}$.
\end{proof}

\begin{lemma}
\label{lemma: min is SB}
An $e_i$-minimal sequence is shift-bounded.
\end{lemma}

\begin{proof}
Let $\seq{s}$ be an $e_i$-minimal sequence. From Lemma \ref{lemma: shift - minimal} and Lemma \ref{lemma: c'>s(c)>c n<2i} we have that $n_{\seq{s}}\geq2^i$. But as $\seq{s}$ is an $e_i$-minimal sequence we have also that $\seq{s}> g_{n}(\seq{s})^{\infty}$ for $2^i\leq n <|\seq{s}|$. Hence $\seq{s}>\seq{s}[1,n]^{\infty}$ for $0<n<|\seq{s}|$, which by Lemma \ref{lemma: shift - minimal} implies $\sigma^n(\seq{s})>\seq{s}$ for $0<n<|\seq{s}|$.

For the upper bounding inequality in the definition of shift-bounded\-ness we have by Lemma \ref{lemma: c'>s(c)>c n<2i} that $\seq{s}' > \sigma^n(\seq{s})$ for $0<n<2^i$. Moreover, by the $e_i$-minimality of $\seq{s}$ we have that $e(g_{n}(\seq{s}))> \seq{s}$ for $2^i\leq n <|\seq{s}|$. For $2^i\leq n <|\seq{s}|$ let $\seq{a} = \seq{s}[1,n]$. 
Then $e(g_{n}(\seq{s})) = \seq{a} ((\tilde{\seq{a}})^*)^{\infty}$ and $\seq{s} = \seq{a}\seq{b}$ for some sequence $\seq{b}$ such that $(\tilde{\seq{a}})^*>\seq{b}$. This implies $\seq{s}' \geq \seq{a}^* > (\tilde{\seq{a}})^* > \seq{b} = \sigma^n(\seq{s})$.
\end{proof}

\begin{example}
There are shift-bounded $e_i$-sequences which are not $e_i$-minimal sequences. The sequence $\seq{s}=000111$ is shift-bounded but not $e_1$-minimal, it has the $e_1$-minimal prefix $001$. 
\end{example}

\begin{lemma}
\label{lemma: u->c}
Let $\seq{c}$ be a finite sequence such that $[\seq{c}]\cap F(\seq{c}) \neq \emptyset$ and let $\seq{u}$ be such that $[\seq{u}]\cap F(\seq{c})\neq \emptyset$ and $|\seq{c}|=|\seq{u}|$. Then there exists $1\leq k\leq |\seq{u}|$ such that $[\seq{u}[1,k]\seq{c}]\cap F(\seq{c})\neq \emptyset$ 
\end{lemma}

\begin{proof}
Let $\seq{w}$ be an infinite sequence such that $\seq{c}\seq{w}\in F(\seq{c})$. Assume there exists a smallest $k$ such that $\seq{u}[k+1,|\seq{u}|] = \seq{c}[1,|\seq{u}|-k+1]$. If we for some $n<k$ would have $\sigma^n(\seq{u}) = \seq{c}[1,|\seq{c}|-n]$ then we would have a contradiction to the choice of $k$. Hence $\sigma^n(\seq{u}) > \seq{c}[1,|\seq{c}|-n]$ for $n<k$ and therefore $\sigma^n(\seq{u}[1,k]\seq{c}\seq{w})\geq \seq{c}$ for $n\geq 0$. 

For any continuation $\seq{v}$ of $\seq{u}$ such that $\seq{u}\seq{v}\in F(\seq{c})$ we have 
\[
\seq{c}' > \sigma^n\big(\seq{u}\,\seq{v} ) \geq \sigma^n(\seq{u}[1,k]\,\, \seq{c}\seq{w}),
\]
for $n\geq0$. If $\seq{c}$ does not overlap $\seq{u}$ then clearly we must have both $\sigma^n(\seq{u}\seq{c}\seq{w}) > \seq{c}$ and $\seq{c}'>\sigma^n(\seq{u}\seq{c}\seq{w})$ for $n \geq 0$.
\end{proof}

\begin{lemma}
\label{lemma: c->w->01}
Let $\seq{c} $ be a finite $e_1$-minimal sequence. Then there exists a finite sequence $\seq{w}$ such that $\seq{c}\, \seq{w}\, (01)^{\infty}\in F(\seq{c})$.
\end{lemma}

\begin{proof}
Let $\seq{a}_k = p^k(\seq{c})$ for $0 \leq k\leq N$ where $N$ is such that $\seq{a}_N = 01$, which exists by Lemma \ref{lemma: pc = f}. Now let 
\[
\seq{b}_k = \seq{a}_k\, (\seq{a}_{k+1})^{\infty} = \seq{u}\,\seq{v}\,\seq{u}^* ( \widetilde{\seq{u}\seq{v}} )^{\infty}.
\]
We claim that $\seq{c}'> \sigma^n(\seq{b}_k) >\seq{c}$ for $0\leq n$. To prove the claim it is enough to prove that it holds for $0\leq n\leq |\seq{a}_k|$ as $\seq{c}' >\sigma^n(\seq{a}_r^{\infty})>\seq{c}$ for $n\geq 0$ and all $0\leq r\leq N$. The lower inequality, $\sigma^n(\seq{b}_k)>\seq{c}$, follows direct from the definition of $p$. For the upper inequality, $\seq{c}'> \sigma^n(\seq{b}_k)$, we start by notice that when $n = 0$ the result follows trivially as $\seq{b}_k$ starts with a 0 while $\seq{c}'$ starts with a 1.

For $0 < n < \frac{1}{2}|\seq{a}_k|$ let
$\seq{\alpha} = \seq{b}_k[1,n]$ and
$\seq{\beta}  = \seq{b}_k[n+1,2n]$.

\footnotesize
\[
\begin{picture}(230,65)
\multiput(40,37)( 0,0){1}{
\put( -35,  0){\makebox(30,12)[r]{$\seq{c}'=$}}
\multiput(0,0)(0,0){1}{
\put(   0,  0){\line(1,0){110}}\put(  0,12){\line(1,0){110}}
\put(   0,  0){\line(0,1){ 12}}\put(110, 0){\line(0,1){ 12}}}
\multiput(111,0)(0,0){1}{
\put(   0,  0){\line(1,0){82}}\put(  0,12){\line(1,0){82}}
\put(   0,  0){\line(0,1){ 12}}}
\put(   1,  0){$\overbrace{\makebox(110,15)[c]{}}_{}$}\put( 52,24){$\seq{a}_k'$}
\multiput(70,2)(0,4){3}{\line(0,1){2}}
\put(  14,  4){$\seq{\alpha}^*$}
\put(  48,  4){$\seq{\beta}^*$}
}
\multiput(40,15)( 0,0){1}{
\put(-35,  0){\makebox(30,12)[r]{$\sigma^n(\seq{b}_k)=$}}
\multiput(0,0)(0,0){1}{
\put(   0,  0){\line(1,0){ 80}}\put(   0,12){\line(1,0){ 80}}
\put(   0,  0){\line(0,1){ 12}}\put(  80, 0){\line(0,1){ 12}}}
\multiput(81,0)(71,0){1}{
\put(   0,  0){\line(1,0){ 70}}\put(   0,12){\line(1,0){ 70}}
\put(   0,  0){\line(0,1){ 12}}\put(  70, 0){\line(0,1){ 12}}}
\multiput(152,0)(0,0){1}{
\put(   0,  0){\line(1,0){ 41}}\put(   0,12){\line(1,0){ 41}}
\put(   0,  0){\line(0,1){ 12}}}
\put(   1, -2){$\underbrace{\makebox(78,1)[c]{}}_{}$}\put( 26,-17){$\sigma^n(\seq{a}_k)$}
\put(  14,  4){$\seq{\beta}$}
\put( 109,  4){$\seq{a}_{k+1}$}
}
\multiput( 75,15)(0,4){9}{\line(0,1){2}}
\end{picture}
\]
\normalsize
As $\seq{c}$ is shift-bounded we have $\seq{\alpha}^*>\seq{\beta}$ and therefore $\seq{c}'>\sigma^n(\seq{b}_k)$.

For $n = \frac{1}{2}|\seq{a}_k|$, and if $\seq{v}$ is void then since $|\seq{a}_k|\geq2$ the $e_1$-minimality of $\seq{c}$ gives $\seq{c}'> \seq{u}^* (\tilde{\seq{u}})^{\infty}= \sigma^n(\seq{b}_k)$. If $\seq{v}$ is non-void then the result follows by the definition of $\seq{a}_{k+1}$ via $p$. 

\footnotesize
\[
\begin{picture}(230,65)
\multiput(40,37)( 0,0){1}{
\put( -35,  0){\makebox(30,12)[r]{$\seq{c}'=$}}
\multiput(0,0)(0,0){1}{
\put(   0,  0){\line(1,0){110}}\put(  0,12){\line(1,0){110}}
\put(   0,  0){\line(0,1){ 12}}\put(110, 0){\line(0,1){ 12}}}
\multiput(111,0)(0,0){1}{
\put(   0,  0){\line(1,0){ 82}}\put(  0,12){\line(1,0){ 82}}
\put(   0,  0){\line(0,1){ 12}}}
\put(   1,  0){$\overbrace{\makebox(110,15)[c]{}}_{}$}\put( 52,24){$\seq{a}_k'$}
\multiput(45,2)(0,4){3}{\line(0,1){2}}
\put(  19,  4){$\seq{u}^*$}
}
\multiput(40,15)( 0,0){1}{
\put(-35,  0){\makebox(30,12)[r]{$\sigma^n(\seq{b}_k)=$}}
\multiput(0,0)(0,0){1}{
\put(   0,  0){\line(1,0){ 55}}\put(   0,12){\line(1,0){ 55}}
\put(   0,  0){\line(0,1){ 12}}\put(  55, 0){\line(0,1){ 12}}}
\multiput(56,0)(71,0){1}{
\put(   0,  0){\line(1,0){ 70}}\put(   0,12){\line(1,0){ 70}}
\put(   0,  0){\line(0,1){ 12}}\put(  70, 0){\line(0,1){ 12}}}
\multiput(127,0)(0,0){1}{
\put(   0,  0){\line(1,0){ 66}}\put(   0,12){\line(1,0){ 66}}
\put(   0,  0){\line(0,1){ 12}}}
\put(   1, -2){$\underbrace{\makebox(53,1)[c]{}}_{}$}\put( 13,-17){$\sigma^n(\seq{a}_k)$}
\multiput(10,0)(0,4){3}{\line(0,1){2}}
\put(  29,  4){$\seq{u}^*$}
\put(  83,  4){$\seq{a}_{k+1}$}
\put( 154,  4){$\seq{a}_{k+1}$}
}
\end{picture}
\]
\normalsize

For $\frac{1}{2}|\seq{a}_k|<n < |\seq{a}_k|$ let
$\seq{\alpha}^* = \seq{c}'[1,|\seq{a}_k|-n]$,
$\seq{\beta}^*  = \seq{c}'[|\seq{a}_k|-n+1, 2|\seq{a}_k|-2n]$ and
$\seq{\gamma}   =\seq{b}_k[n+1, |\seq{a}_k|]$.

\footnotesize
\[
\begin{picture}(230,65)
\multiput(40,37)( 0,0){1}{
\put( -35,  0){\makebox(30,12)[r]{$\seq{c}'=$}}
\multiput(0,0)(0,0){1}{
\put(   0,  0){\line(1,0){110}}\put(  0,12){\line(1,0){110}}
\put(   0,  0){\line(0,1){ 12}}\put(110, 0){\line(0,1){ 12}}}
\multiput(111,0)(0,0){1}{
\put(   0,  0){\line(1,0){ 82}}\put(  0,12){\line(1,0){ 82}}
\put(   0,  0){\line(0,1){ 12}}}
\put(   1,  0){$\overbrace{\makebox(110,15)[c]{}}_{}$}\put( 52,24){$\seq{a}_k'$}
\put(  18,  4){$\seq{\alpha}^*$}
\put(  58,  4){$\seq{\beta}^*$}
}
\multiput(40,15)( 0,0){1}{
\put(-35,  0){\makebox(30,12)[r]{$\sigma^n(\seq{b}_k)=$}}
\multiput(0,0)(0,0){1}{
\put(   0,  0){\line(1,0){ 40}}\put(   0,12){\line(1,0){ 40}}
\put(   0,  0){\line(0,1){ 12}}\put(  40, 0){\line(0,1){ 12}}}
\multiput(41,0)(71,0){2}{
\put(   0,  0){\line(1,0){ 70}}\put(   0,12){\line(1,0){ 70}}
\put(   0,  0){\line(0,1){ 12}}\put(  70, 0){\line(0,1){ 12}}}
\multiput(183,0)(0,0){1}{
\put(   0,  0){\line(1,0){ 10}}\put(   0,12){\line(1,0){ 10}}
\put(   0,  0){\line(0,1){ 12}}}
\put(   1, -2){$\underbrace{\makebox(38,1)[c]{}}_{}$}\put(  5,-17){$\sigma^n(\seq{a}_k)$}
\put(  41, -2){$\underbrace{\makebox(68,1)[c]{}}_{}$}\put( 68,-17){$\seq{a}_{k+1}$}
\put(  18,  4){$\seq{\gamma}$}
\put(  58,  4){$\seq{\alpha}$}
\put( 139,  4){$\seq{a}_{k+1}$}
}
\multiput( 80,15)(0,4){9}{\line(0,1){2}}
\multiput(120,15)(0,4){9}{\line(0,1){2}}
\end{picture}
\]
\normalsize
We have $\seq{\alpha}^*\geq \seq{\gamma}$ and $\seq{\beta}^* >\seq{\alpha}$. If $\seq{u}$ is void we have directly $\seq{\alpha}^*>\seq{\gamma}$. Therefore $\seq{c}'>\sigma^n(\seq{b}_k)$, which proves the claim. Put $\seq{w} = \seq{a}_1^{n_1}\seq{a}_2^{n_2}\ldots \seq{a}_N$ with $n_k = [ \frac{|\seq{c}|}{ |\seq{a}_k|}  ] +1$. By repeated use of the just proved claim we have $\seq{c}\,\seq{w}\,(01)^{\infty} \in F(\seq{c})$. 
\end{proof}

\begin{theorem}
\label{thm: mixing}
Let $\seq{c}$ be a finite $e_1$-minimal sequence. Then $\sigma:F(\seq{c})\to F(\seq{c})$ is topologically mixing.
\end{theorem}

\begin{proof}
Let $U = [\seq{u}]\cap F(\seq{c})$ and $V = [\seq{v}] \cap F(\seq{c})$ and assume they are both non-empty. By Lemma \ref{lemma: u->c} there is a $k$ such that $[\seq{u}[1,k] \,\seq{c}] \cap U$ is non-empty. Lemma \ref{lemma: c->w->01} gives that there is a finite sequence $\seq{w}$ such that $\seq{u}[1,k] \,\seq{c}\,\seq{w}\,(01)^{\infty} \in U$.  Let $\seq{a} = 0$ if $v_1 = 1$ and let $\seq{a}$ be void if $v_1=0$. Then there exists a positive integer $N_1$ such that
\begin{equation}	
\label{eq: u-v 1}
\big[\seq{u}[1,k] \,\seq{c}\,\seq{w}\,(01)^{n_1} \seq{a}\, \seq{v}\big] \cap U \neq \emptyset 
\end{equation}
for $n_1>N_1$. As $\seq{c}$ is a finite $e_1$-minimal sequence there exist $N_2$ and $N_3$ such that  
\begin{equation}	
\label{eq: u-v 2}
\big[\seq{u}[1,k] \,\seq{c}\,\seq{w}\,(01)^{n_2}0(01)^{n_3} \seq{a}\, \seq{v}\big] \cap U \neq \emptyset
\end{equation}
for $n_2>N_2$ and $n_3>N_3$. Combining (\ref{eq: u-v 1}) and (\ref{eq: u-v 2}) gives $\sigma^n(U)\cap V\neq \emptyset$ for all $n$ lager than some $N_0$.
\end{proof}

\begin{example}
Letting $\seq{c}$ be an $e_1$-sequence is crucial in Lemma \ref{lemma: c->w->01} and Theorem \ref{thm: mixing}. If we assume that $\seq{c}$ is a finite minimal $e_i$-minimal sequence for $i\geq2$ then $\seq{c}$ must have a prefix $\seq{p}$ of the form $\tilde{\seq{u}}(\seq{u}^*)^k\seq{u}'$ for some $k>0$. Lemma \ref{lemma: uuu} now gives that we can never find a sequence $\seq{w}$ such that $\seq{c}\seq{w}(01)^{\infty}$ is a sequence in $F(\seq{c})$, and therefore we do not have topologically mixing. 
\end{example}

\begin{corollary}
\label{cor: A primitiv}
Let $\seq{c}$ be a finite $e_1$-minimal sequence. Then the transition matrix $A_{\seq{c}}$ corresponding to $F(\seq{c})$ is primitive.
\end{corollary}

We end the section by proving two accumulation results on finite $e_i$-minimal sequences.

\begin{lemma}
\label{lemma: e-a_k}
Let $\seq{s}$ be a finite $e_i$-minimal sequence and let $\seq{s} = \seq{u} \seq{v} \seq{u}^*$ where $p(\seq{s}) = \widetilde{\seq{uv}}$. Put $\seq{a}_k(\seq{s}) = \tilde{\seq{s}}(\seq{s}^*)^k \seq{u}^*$ for $k \geq 1$. Then the $\seq{a}_k$'s are $e_i$-minimal and $\seq{a}_k \nearrow e(\seq{s})$ when $k$ tends to infinity.
\end{lemma}

\begin{proof}
We first have to show that $\seq{a}_k$ is an $e_i$-sequence. As $\seq{s}$ is a finite $e_i$-sequence we have that $\seq{e}_{i-1}[1,|\seq{s}|] < \seq{s}$ and since $\seq{s}$ is $e_i$-minimal we must have $2^i<|\seq{s}|$, (as there are no shift-bounded $e_i$-sequences of length $2^i$). We only have to consider the case when $\seq{e}_{i-1}[1,|\seq{s}|] = \tilde{\seq{s}}$. To do so, let 
$\seq{w} = f^{i-1}(1)$, 
$\seq{\alpha}^* = \seq{a}_k[|\seq{s}|+1, t|\seq{w}|]$ and
$\seq{\gamma}^* = \seq{e}_{i-1}[|\seq{s}|+1, t|\seq{w}|]$, where $t$ is the smallest integer such that $t|\seq{w}|>|\seq{s}|$.

\footnotesize
\[
\begin{picture}(230,65)
\multiput(40,37)( 0,0){1}{
\put( -35,  0){\makebox(30,12)[r]{$\seq{e}_{i-1}=\ldots$}}
\multiput(0,0)(0,0){1}{
\put(   0,  0){\line(1,0){ 60}}\put(   0,12){\line(1,0){ 60}}
\put(  60, 0){\line(0,1){ 12}}}
\multiput(61,0)(0,0){1}{
\put(   0,  0){\line(1,0){ 80}}\put(   0,12){\line(1,0){ 80}}
\put(   0,  0){\line(0,1){ 12}}\put(  80, 0){\line(0,1){ 12}}}
\multiput(142,0)(0,0){1}{
\put(   0,  0){\line(1,0){ 51}}\put(   0,12){\line(1,0){ 51}}
\put(   0,  0){\line(0,1){ 12}}}
\put(  62,  0){$\overbrace{\makebox(78,15)[c]{}}_{}$}\put( 99,24){$\seq{w}^*$}
\put(  27,  4){$\seq{w}^*$}
\put( 112,  4){$\seq{\gamma}^*$}
\put( 165,  4){$\seq{w}^*$}
}
\multiput(40,15)( 0,0){1}{
\put(-35,  0){\makebox(30,12)[r]{$\seq{a}_k=\ldots$}}
\multiput(0,0)(0,0){1}{
\put(   0,  0){\line(1,0){ 90}}\put(   0,12){\line(1,0){ 90}}
\put(  90,  0){\line(0,1){ 12}}}
\multiput(91,0)(0,0){1}{
\put(   0,  0){\line(1,0){102}}\put(   0,12){\line(1,0){102}}
\put(   0,  0){\line(0,1){ 12}}}
\put( 112,  4){$\seq{\alpha}^*$}
\put(  37,  4){$\tilde{\seq{s}}$}
\put( 165,  4){$\seq{s}^*$}
}
\multiput(131,15)(0,4){9}{\line(0,1){2}}
\multiput(181,15)(0,4){9}{\line(0,1){2}}
\end{picture}
\]
\normalsize
As $\seq{w}$ is shift-bounded and that $\tilde{\seq{s}}=\seq{e}_i[1,|\seq{s}|]$ we have $\seq{\alpha}^*>\seq{\gamma}^*$ and therefore $\seq{e}_{i-1}<\seq{a}_k$. Hence we have $\seq{e}_{i-1} < \seq{a}_k < \seq{s} < \seq{e}_i$, since $\seq{s}$ is an $e_i$-sequence, and therefore we see that $\seq{a}_k$ is an $e_i$-sequence.

For the minimality we have to show that at least one of the two inequalities
\begin{equation}
\label{eq: e(a)>a}
e(g_n(\seq{a}_k)) > \seq{a}_k
\end{equation}
and
\begin{equation}
\label{eq: a>g(a)}
\seq{a}_k > g_n(\seq{a}_k)^{\infty}
\end{equation}
hold for $2^i\leq n<|\seq{a}_k|$. Let us first turn to the inequality (\ref{eq: e(a)>a}). It is clear that (\ref{eq: e(a)>a}) fails whenever $n$ is such that $(\seq{a}_k)_n=1$, hence we may assume that $n$ is such that $(\seq{a}_k)_n=0$. For $2^{i}\leq n < |\seq{s}|$ we have that $e(g_n(\seq{a}_k)) = e(g_n(\seq{s}))>\seq{s}>\seq{a}_k$, which gives that (\ref{eq: e(a)>a}) holds in this case.

For $n = r|\seq{s}|$ with $1<r<k$ we have that (\ref{eq: e(a)>a}) holds as $(\tilde{\seq{s}})^*>\seq{s}^*\geq \seq{u}^*$.

For $j|\seq{s}|<n<(j+1)|\seq{s}|$ with $0<j<k$ let 
$\seq{\alpha}^* = e(g_n(\seq{a}_k))[n+1, j|\seq{s}|]$ and 
$\seq{\gamma}^* = \seq{a}_k[n+1, j|\seq{s}|]$.

\footnotesize
\[
\begin{picture}(230,65)
\multiput(40,37)( 0,0){1}{
\put(-35,  0){\makebox(30,12)[r]{$e(g_n(\seq{a}_k))=\ldots$}}
\multiput(0,0)(0,0){1}{
\put(   0,  0){\line(1,0){ 90}}\put(   0,12){\line(1,0){90}}
\put(  90,  0){\line(0,1){ 12}}}
\multiput( 91,0)(0,0){1}{
\put(   0,  0){\line(1,0){102}}\put(   0,12){\line(1,0){102}}
\put(   0,  0){\line(0,1){ 12}}}
\put(  92,  0){$\overbrace{\makebox(78,15)[c]{}}_{}$}\put( 129, 24){$\seq{s}^*$}
\put( 112,  4){$\seq{\alpha}^*$}
\multiput(170,2)(0,4){3}{\line(0,1){2}}
}
\multiput(40,15)( 0,0){1}{
\put( -35,  0){\makebox(30,12)[r]{$\seq{a}_k=\ldots$}}
\multiput(0,0)(0,0){1}{
\put(   0,  0){\line(1,0){ 60}}\put(   0,12){\line(1,0){ 60}}
\put(  60, 0){\line(0,1){ 12}}}
\multiput(61,0)(0,0){1}{
\put(   0,  0){\line(1,0){ 80}}\put(   0,12){\line(1,0){ 80}}
\put(   0,  0){\line(0,1){ 12}}\put(  80, 0){\line(0,1){ 12}}}
\multiput(142,0)(0,0){1}{
\put(   0,  0){\line(1,0){ 51}}\put(   0,12){\line(1,0){ 51}}
\put(   0,  0){\line(0,1){ 12}}\multiput( 30,0)(0,4){3}{\line(0,1){2}}}
\put(  62, -2){$\underbrace{\makebox(78,1)[c]{}}_{}$}\put( 97,-17){$\seq{s}^*$}
\put(  27,  4){$\seq{s}^*$}
\put( 112,  4){$\seq{\gamma}^*$}
\put( 155,  4){$\seq{u}^*$}
}
\multiput(131,15)(0,4){9}{\line(0,1){2}}
\multiput(181,15)(0,4){9}{\line(0,1){2}}
\end{picture}
\]
\normalsize
As $\seq{s}$ is shift-bounded we have $\seq{\alpha}^*>\seq{\gamma}^*$ and therefore (\ref{eq: e(a)>a}) holds. 

For $k|\seq{s}|< n < |\seq{a}_k|$ let 
$\seq{\alpha}^* = e(g_n(\seq{a}_k))[n+1, |\seq{a}_k|]$ and 
$\seq{\gamma}^* = \seq{a}_k[n+1, |\seq{a}_k|]$.

\footnotesize
\[
\begin{picture}(230,65)
\multiput(40,37)( 0,0){1}{
\put(-35,  0){\makebox(30,12)[r]{$e(g_n(\seq{a}_k))=\ldots$}}
\multiput(0,0)(0,0){1}{
\put(   0,  0){\line(1,0){150}}\put(   0,12){\line(1,0){150}}
\put( 150,  0){\line(0,1){ 12}}}
\multiput(151,0)(0,0){1}{
\put(   0,  0){\line(1,0){42}}\put(   0,12){\line(1,0){42}}
\put(   0,  0){\line(0,1){ 12}}}
\put( 157,  4){$\seq{\alpha}^*$}
}
\multiput(40,15)( 0,0){1}{
\put( -35,  0){\makebox(30,12)[r]{$\seq{a}_k=\ldots$}}
\multiput(0,0)(0,0){1}{
\put(   0,  0){\line(1,0){ 60}}\put(   0,12){\line(1,0){ 60}}
\put(  60,  0){\line(0,1){ 12}}}
\multiput(61,0)(0,0){1}{
\put(   0,  0){\line(1,0){ 80}}\put(   0,12){\line(1,0){ 80}}
\put(   0,  0){\line(0,1){ 12}}\put(  80, 0){\line(0,1){ 12}}}
\multiput(142,0)(0,0){1}{
\put(   0,  0){\line(1,0){ 30}}\put(   0,12){\line(1,0){ 30}}
\put(   0,  0){\line(0,1){ 12}}\put(  30, 0){\line(0,1){ 12}}}
\put( 142, -2){$\underbrace{\makebox(28,1)[c]{}}_{}$}\put( 152,-17){$\seq{u}^*$}
\put(  27,  4){$\seq{s}^*$}
\put(  98,  4){$\seq{s}^*$}
\put( 157,  4){$\seq{\gamma}^*$}
}
\multiput(191,15)(0,4){9}{\line(0,1){2}}
\multiput(212,15)(0,4){9}{\line(0,1){2}}
\end{picture}
\]
\normalsize
The shift-boundedness of $\seq{s}$ and the definition of $\seq{u}$ gives $\seq{\alpha}^*>\seq{\gamma}^*$ and hence (\ref{eq: e(a)>a}) holds.

Now let us turn to the inequality (\ref{eq: a>g(a)}). It is clear that (\ref{eq: a>g(a)}) fails whenever $n$ is such that $(\seq{a}_k)_n=0$, hence we may assume that $n$ is such that $(\seq{a}_k)_n=1$. For $2^i \leq n < |\seq{s}|$ we have as $\seq{s}$ is $e_i$-minimal that $\tilde{\seq{a}}_k[1,|\seq{s}|] \geq g_n(\seq{a}_k)^{\infty}[1,|\seq{s}|]$. If the inequality is strict we are done, hence we only have to consider the case when having equality, $\tilde{\seq{a}}_k[1,|\seq{s}|] = g_n(\seq{a}_k)^{\infty}[1,|\seq{s}|]$.
Let $r$ be the smallest integer such that $nr>|\seq{s}|$. 

If $nr - |\seq{s}| > \frac{1}{2}n$ then let
$\seq{z} = g_n(\seq{a}_k)$,
$\seq{\alpha}^* = \seq{a}_k[|\seq{s}|+1, 2|\seq{s}|-n(r-1)]$ and
$\seq{\beta}  = g_n(\seq{a}_k)^{\infty}[|\seq{s}|+1, 2|\seq{s}|-n(r-1)]$.

\footnotesize
\[
\begin{picture}(230,65)
\multiput(40,37)( 0,0){1}{
\put( -35,  0){\makebox(30,12)[r]{$g_n(\seq{a}_k)^{\infty}=\ldots$}}
\multiput(0,0)(0,0){1}{
\put(   0,  0){\line(1,0){ 60}}\put(   0,12){\line(1,0){ 60}}
\put(  60, 0){\line(0,1){ 12}}}
\multiput(61,0)(0,0){1}{
\put(   0,  0){\line(1,0){ 80}}\put(   0,12){\line(1,0){ 80}}
\put(   0,  0){\line(0,1){ 12}}\put(  80, 0){\line(0,1){ 12}}}
\multiput(142,0)(0,0){1}{
\put(   0,  0){\line(1,0){ 51}}\put(   0,12){\line(1,0){ 51}}
\put(   0,  0){\line(0,1){ 12}}}
\put(  62,  0){$\overbrace{\makebox(78,15)[c]{}}_{}$}\put( 99,24){$\seq{z}$}
\put(  27,  4){$\seq{z}$}
\put(  73,  4){$\seq{\alpha}$}
\put( 103,  4){$\seq{\beta}$}
\put( 165,  4){$\seq{z}$}
}
\multiput(40,15)( 0,0){1}{
\put(-35,  0){\makebox(30,12)[r]{$\seq{a}_k=\ldots$}}
\multiput(0,0)(0,0){1}{
\put(   0,  0){\line(1,0){ 90}}\put(   0,12){\line(1,0){ 90}}
\put(  90,  0){\line(0,1){ 12}}}
\multiput(91,0)(0,0){1}{
\put(   0,  0){\line(1,0){102}}\put(   0,12){\line(1,0){102}}
\put(   0,  0){\line(0,1){ 12}}}
\put(  92, -2){$\underbrace{\makebox(78,1)[c]{}}_{}$}\put( 129,-17){$\seq{z}^*$}
\put(  45,  4){$\tilde{\seq{s}}$}
\put( 103,  4){$\seq{\alpha}^*$}
\multiput(170,0)(0,4){3}{\line(0,1){2}}
}
\multiput(131,15)(0,4){9}{\line(0,1){2}}
\multiput(161,15)(0,4){9}{\line(0,1){2}}
\end{picture}
\]
\normalsize
If we assume that $n$ is the smallest integer such that (\ref{eq: a>g(a)}) does not hold then $\seq{z}$ is the $e_i$-minimal prefix of $\seq{a}_k$. But then $\seq{z}$ is shift-bounded and we must have $\seq{\alpha}^*>\seq{\beta}$, which contradicts that $\seq{z}$ is the $e_i$-minimal prefix of $\seq{a}_k$.

If $nr - |\seq{s}| \leq \frac{1}{2}n$ then let
$\seq{z} = g_n(\seq{a}_k)$,
$\seq{\alpha}^* = \seq{a}_k[|\seq{s}|+1, nr]$,
$\seq{\gamma}   = g_n(\seq{a}_k)^{\infty}[|\seq{s}|+1, nr]$ and
$\seq{\beta}  = \seq{a}_k[nr+1,2nr-|\seq{s}| ]$.

\footnotesize
\[
\begin{picture}(230,65)
\multiput(40,37)( 0,0){1}{
\put( -35,  0){\makebox(30,12)[r]{$g_n(\seq{a}_k)^{\infty}=\ldots$}}
\multiput(0,0)(0,0){1}{
\put(   0,  0){\line(1,0){ 60}}\put(   0,12){\line(1,0){ 60}}
\put(  60, 0){\line(0,1){ 12}}}
\multiput(61,0)(0,0){1}{
\put(   0,  0){\line(1,0){ 80}}\put(   0,12){\line(1,0){ 80}}
\put(   0,  0){\line(0,1){ 12}}\put(  80, 0){\line(0,1){ 12}}}
\multiput(142,0)(0,0){1}{
\put(   0,  0){\line(1,0){ 51}}\put(   0,12){\line(1,0){ 51}}
\put(   0,  0){\line(0,1){ 12}}}
\put(  62,  0){$\overbrace{\makebox(78,15)[c]{}}_{}$}\put(  99,24){$\seq{z}$}
\put(  27,  4){$\seq{z}$}
\put( 125,  4){$\seq{\gamma}$}
\put( 152,  4){$\seq{\alpha}$}
}
\multiput(40,15)( 0,0){1}{
\put(-35,  0){\makebox(30,12)[r]{$\seq{a}_k=\ldots$}}
\multiput(0,0)(0,0){1}{
\put(   0,  0){\line(1,0){115}}\put(   0,12){\line(1,0){115}}
\put( 115,  0){\line(0,1){ 12}}}
\multiput(116,0)(0,0){1}{
\put(   0,  0){\line(1,0){77}}\put(   0,12){\line(1,0){77}}
\put(   0,  0){\line(0,1){ 12}}}
\put(  55,  4){$\tilde{\seq{s}}$}
\put( 125,  4){$\seq{\alpha}^*$}
\put( 152,  4){$\seq{\beta}^*$}
}
\multiput(156,15)(0,4){9}{\line(0,1){2}}
\multiput(181,15)(0,4){9}{\line(0,1){2}}
\multiput(206,15)(0,4){9}{\line(0,1){2}}
\end{picture}
\]
\normalsize
If we again assume that $n$ is the smallest integer such that (\ref{eq: a>g(a)}) does not hold then $\seq{z}$ is the $e_i$-minimal prefix of $\seq{a}_k$. But then $\seq{z}$ is shift-bounded and we must have $\seq{\alpha}^*\geq \seq{\gamma}$ and $\seq{\beta}^*>\seq{\alpha}$, which contradicts that $\seq{z}$ is the $e_i$-minimal prefix of $\seq{a}_k$.

For $n = r|\seq{s}|$ with $1<r<k$ the inequality (\ref{eq: a>g(a)}) holds because $\seq{s}^* > \tilde{\seq{s}}$.

Let $0<j<k$. Then 
for $j|\seq{s}|<n<(j+1)|\seq{s}|-|\seq{u}|$ let
$\seq{\alpha} = g_n(\seq{a}_k)^{\infty}[n+1, (j+1)|\seq{s}|]$ and
$\seq{\gamma} = \seq{a}_k[n+1, (j+1)|\seq{s}|]$.

\footnotesize
\[
\begin{picture}(230,65)
\multiput(40,37)( 0,0){1}{
\put( -35,  0){\makebox(30,12)[r]{$\seq{a}_k=\ldots$}}
\multiput(0,0)(0,0){1}{
\put(   0,  0){\line(1,0){ 60}}\put(   0,12){\line(1,0){ 60}}
\put(  60, 0){\line(0,1){ 12}}}
\multiput(61,0)(0,0){1}{
\put(   0,  0){\line(1,0){ 80}}\put(   0,12){\line(1,0){ 80}}
\put(   0,  0){\line(0,1){ 12}}\put(  80, 0){\line(0,1){ 12}}}
\multiput(142,0)(0,0){1}{
\put(   0,  0){\line(1,0){ 51}}\put(   0,12){\line(1,0){ 51}}
\put(   0,  0){\line(0,1){ 12}}\multiput( 30,2)(0,4){3}{\line(0,1){2}}}
\put(  62,  0){$\overbrace{\makebox(78,15)[c]{}}_{}$}\put( 99,24){$\seq{s}^*$}
\put(  27,  4){$\seq{s}^*$}
\put( 112,  4){$\seq{\gamma}^*$}
\put( 155,  4){$\seq{u}^*$}
}
\multiput(40,15)( 0,0){1}{
\put(-35,  0){\makebox(30,12)[r]{$g_n(\seq{a}_k)^{\infty}=\ldots$}}
\multiput(0,0)(0,0){1}{
\put(   0,  0){\line(1,0){ 90}}\put(   0,12){\line(1,0){ 90}}
\put(  90,  0){\line(0,1){ 12}}}
\multiput(91,0)(0,0){1}{
\put(   0,  0){\line(1,0){102}}\put(   0,12){\line(1,0){102}}
\put(   0,  0){\line(0,1){ 12}}}
\put(  92, -2){$\underbrace{\makebox(78,1)[c]{}}_{}$}\put( 129,-17){$\tilde{\seq{s}}$}
\put( 112,  4){$\seq{\alpha}$}
\multiput(170,0)(0,4){3}{\line(0,1){2}}
}
\multiput(131,15)(0,4){9}{\line(0,1){2}}
\multiput(181,15)(0,4){9}{\line(0,1){2}}
\end{picture}
\]
\normalsize
By the shift-boundedness of $\seq{s}$ and the definition of $\seq{u}$ we must have $\seq{\gamma}^*>\seq{\alpha}$, which implies (\ref{eq: a>g(a)}).

For $(j+1)|\seq{s}|-|\seq{u}|\leq n <(j+1)|\seq{s}|$ we let
$\seq{\alpha} = g_n(\seq{a}_k)^{\infty}[n, (j+1)|\seq{s}|]$,
$\seq{\gamma} = \seq{a}_k[n+1, (j+1)|\seq{s}|]$ and 
$\seq{\beta} =  g_n(\seq{a}_k)^{\infty}[(j+1)|\seq{s}|+1, 2(j+1)|\seq{s}|-n]$.

\footnotesize
\[
\begin{picture}(230,65)
\multiput(40,37)( 0,0){1}{
\put( -35,  0){\makebox(30,12)[r]{$\seq{a}_k=\ldots$}}
\multiput(0,0)(0,0){1}{
\put(   0,  0){\line(1,0){ 60}}\put(   0,12){\line(1,0){ 60}}
\put(  60, 0){\line(0,1){ 12}}}
\multiput(61,0)(0,0){1}{
\put(   0,  0){\line(1,0){ 80}}\put(   0,12){\line(1,0){ 80}}
\put(   0,  0){\line(0,1){ 12}}\put(  80, 0){\line(0,1){ 12}}}
\multiput(142,0)(0,0){1}{
\put(   0,  0){\line(1,0){ 51}}\put(   0,12){\line(1,0){ 51}}
\put(   0,  0){\line(0,1){ 12}}\multiput( 30,2)(0,4){3}{\line(0,1){2}}}
\put(  62,  0){$\overbrace{\makebox(78,15)[c]{}}_{}$}\put(  99,24){$\seq{s}^*$}
\put( 142,  0){$\overbrace{\makebox(28,15)[c]{}}_{}$}\put( 152,24){$\seq{u}^*$}
\put(  27,  4){$\seq{s}^*$}
\put( 125,  4){$\seq{\gamma}^*$}
\put( 152,  4){$\seq{\alpha}^*$}
}
\multiput(40,15)( 0,0){1}{
\put(-35,  0){\makebox(30,12)[r]{$g_n(\seq{a}_k)^{\infty}=\ldots$}}
\multiput(0,0)(0,0){1}{
\put(   0,  0){\line(1,0){115}}\put(   0,12){\line(1,0){115}}
\put( 115,  0){\line(0,1){ 12}}}
\multiput(116,0)(0,0){1}{
\put(   0,  0){\line(1,0){77}}\put(   0,12){\line(1,0){77}}
\put(   0,  0){\line(0,1){ 12}}}
\put( 125,  4){$\seq{\alpha}$}
\put( 152,  4){$\seq{\beta}$}
}
\multiput(156,15)(0,4){9}{\line(0,1){2}}
\multiput(181,15)(0,4){9}{\line(0,1){2}}
\multiput(206,15)(0,4){9}{\line(0,1){2}}
\end{picture}
\]
\normalsize
Again, as $\seq{s}$ is shift-bounded we have $\seq{\gamma}^*\seq{\alpha}^* > \seq{\alpha}\seq{\beta}$ since $\seq{\gamma}^*\geq \seq{\alpha}$ and $\seq{\alpha}^*> \seq{\beta}$, which gives (\ref{eq: a>g(a)}).

For $k|\seq{s}|<n<|\seq{a}_k|$ let
$\seq{\gamma} = \seq{a}_k[n+1,|\seq{a}_k|]$ and
$\seq{\alpha} = g_n(\seq{a}_k)^{\infty}[n+1,|\seq{a}_k|]$.

\footnotesize
\[
\begin{picture}(230,65)
\multiput(40,37)( 0,0){1}{
\put( -35,  0){\makebox(30,12)[r]{$\seq{a}_k=\ldots$}}
\multiput(0,0)(0,0){1}{
\put(   0,  0){\line(1,0){ 60}}\put(   0,12){\line(1,0){ 60}}
\put(  60, 0){\line(0,1){ 12}}}
\multiput(61,0)(0,0){1}{
\put(   0,  0){\line(1,0){ 80}}\put(   0,12){\line(1,0){ 80}}
\put(   0,  0){\line(0,1){ 12}}\put(  80, 0){\line(0,1){ 12}}}
\multiput(142,0)(0,0){1}{
\put(   0,  0){\line(1,0){ 30}}\put(   0,12){\line(1,0){ 30}}
\put(   0,  0){\line(0,1){ 12}}\put(  30, 0){\line(0,1){ 12}}}
\put( 142,  0){$\overbrace{\makebox(28,15)[c]{}}_{}$}\put( 152,24){$\seq{u}^*$}
\put(  27,  4){$\seq{s}^*$}
\put(  98,  4){$\seq{s}^*$}
\put( 158,  4){$\seq{\gamma}$}
}
\multiput(40,15)( 0,0){1}{
\put(-35,  0){\makebox(30,12)[r]{$g_n(\seq{a}_k)^{\infty}=\ldots$}}
\multiput(0,0)(0,0){1}{
\put(   0,  0){\line(1,0){150}}\put(   0,12){\line(1,0){150}}
\put( 150, 0){\line(0,1){ 12}}}
\multiput(151,0)(0,0){1}{
\put(   0,  0){\line(1,0){42}}\put(   0,12){\line(1,0){42}}
\put(   0,  0){\line(0,1){ 12}}}
\put( 158,  4){$\seq{\alpha}$}
}
\multiput(191,15)(0,4){9}{\line(0,1){2}}
\multiput(212,15)(0,4){9}{\line(0,1){2}}
\end{picture}
\]
\normalsize
Then as $\seq{s}$ is shift-bounded and having $\seq{u}^*$ as a suffix we must have $\seq{\gamma}>\seq{\alpha}$, which again gives (\ref{eq: a>g(a)}) and completing the proof.
\end{proof}

\begin{lemma}
\label{lemma: e-b_k}
Let $\seq{s}$ be a finite $e_i$-minimal sequence. Define $\seq{b}_k(\seq{s}) = \seq{s}^k p(\seq{s})$ for $k\geq 1$. Then 
\begin{enumerate}
\item the following chain of inequalities holds $\seq{e}_{i-1} \leq \seq{s} < \seq{b}_k < \seq{e}_i[1,|\seq{b}_k|] < \seq{e}_i$ for all $k\geq 1$, (in particular $\seq{b}_k$ is an $e_i$-sequence),  

\item the inequality $e(g_n(\seq{b}_k))[1,|\seq{b}_k|] > \seq{b}_k$ holds for $2^i\leq n \leq k|\seq{s}|$ and the inequality $e(g_n(\seq{b}_k))[1,|\seq{b}_k|] \geq \seq{b}_k$ holds for $k|\seq{s}|< n \leq [\seq{b}_k|$.

\item the $\seq{b}_k$'s are $e_i$-minimal and $\seq{b}_k \searrow \seq{s}^{\infty}$ when $k$ tends to infinity.
\end{enumerate}
\end{lemma}

\begin{proof}
\textit{(1.)} 
The inequalities $\seq{e}_{i-1} \leq \seq{s} < \seq{b}_k$ are clear by definition. The only inequality we have to prove is $\seq{b}_k < \seq{e}_i[1,|\seq{b}_k|]$. 
As $\seq{s}$ is a finite $e_i$-sequence we have that $\seq{s} \leq \seq{e}_i[1,|\seq{s}|]$ and since $\seq{s}$ is $e_i$-minimal we must have $2^i<|\seq{s}|$. We only have to consider the case when $\seq{s} = \seq{e}_i[1,|\seq{s}|]$. To do so, let $\seq{w} = f^{i}(1)$ and let $t$ be the smallest integer such that $t|\seq{w}|>|\seq{s}|$.

If $t|\seq{w}|-|\seq{s}| > \frac{1}{2}|\seq{w}|$ then let
$\seq{\alpha}^*= \seq{b}_k[(t-1)|\seq{w}|+1, |\seq{s}|]$ and
$\seq{\beta}^* = \seq{e}_i[|\seq{s}|+1, 2|\seq{s}|-(t-1)|\seq{w}|]$.

\footnotesize
\[
\begin{picture}(230,65)
\multiput(40,37)( 0,0){1}{
\put(-35,  0){\makebox(30,12)[r]{$\seq{b}_k=\ldots$}}
\multiput(0,0)(0,0){1}{
\put(   0,  0){\line(1,0){ 90}}\put(   0,12){\line(1,0){ 90}}
\put(  90,  0){\line(0,1){ 12}}}
\multiput(91,0)(0,0){1}{
\put(   0,  0){\line(1,0){102}}\put(   0,12){\line(1,0){102}}
\put(   0,  0){\line(0,1){ 12}}}
\put(  92,  0){$\overbrace{\makebox(78,15)[c]{}}_{}$}\put( 129, 24){$\tilde{\seq{w}}$}
\put(  45,  4){$\seq{s}$}
\put( 103,  4){$\seq{\alpha}$}
\multiput(170,2)(0,4){3}{\line(0,1){2}}
}
\multiput(40,15)( 0,0){1}{
\put( -35,  0){\makebox(30,12)[r]{$\seq{e}_i=\ldots$}}
\multiput(0,0)(0,0){1}{
\put(   0,  0){\line(1,0){ 60}}\put(   0,12){\line(1,0){ 60}}
\put(  60, 0){\line(0,1){ 12}}}
\multiput(61,0)(0,0){1}{
\put(   0,  0){\line(1,0){ 80}}\put(   0,12){\line(1,0){ 80}}
\put(   0,  0){\line(0,1){ 12}}\put(  80, 0){\line(0,1){ 12}}}
\multiput(142,0)(0,0){1}{
\put(   0,  0){\line(1,0){ 51}}\put(   0,12){\line(1,0){ 51}}
\put(   0,  0){\line(0,1){ 12}}}
\put(  62, -2){$\underbrace{\makebox(78,1)[c]{}}_{}$}\put( 99, -17){$\seq{w}^*$}
\put(  27,  4){$\seq{w}^*$}
\put(  73,  4){$\seq{\alpha}^*$}
\put( 103,  4){$\seq{\beta}^*$}
\put( 165,  4){$\seq{w}^*$}
}
\multiput(131,15)(0,4){9}{\line(0,1){2}}
\multiput(161,15)(0,4){9}{\line(0,1){2}}
\end{picture}
\]
\normalsize
Since $\seq{w}$ is shift-bounded we must have $\seq{\alpha}<\seq{\beta}^*$, which implies $\seq{b}_k<\seq{e}_i[1,|\seq{b}_k|]$.

If $t|\seq{w}|-|\seq{s}|\leq \frac{1}{2}|\seq{w}|$ then let 
$\seq{\alpha}   = \seq{b}_k[|\seq{s}|+1, t|\seq{w}|]$,
$\seq{\beta}    = \seq{b}_k[t|\seq{w}|+1, 2t|\seq{w}|-|\seq{s}|]$ and 
$\seq{\gamma}^* = \seq{e}_i[|\seq{s}|+1, t|\seq{w}|]$, 

\footnotesize
\[
\begin{picture}(230,65)
\multiput(40,37)( 0,0){1}{
\put(-35,  0){\makebox(30,12)[r]{$\seq{b}_k=\ldots$}}
\multiput(0,0)(0,0){1}{
\put(   0,  0){\line(1,0){115}}\put(   0,12){\line(1,0){115}}
\put( 115,  0){\line(0,1){ 12}}}
\multiput(116,0)(0,0){1}{
\put(   0,  0){\line(1,0){77}}\put(   0,12){\line(1,0){77}}
\put(   0,  0){\line(0,1){ 12}}}
\put(  55,  4){$\seq{s}$}
\put( 125,  4){$\seq{\alpha}$}
\put( 152,  4){$\seq{\beta}$}
}
\multiput(40,15)( 0,0){1}{
\put( -35,  0){\makebox(30,12)[r]{$\seq{e}_i=\ldots$}}
\multiput(0,0)(0,0){1}{
\put(   0,  0){\line(1,0){ 60}}\put(   0,12){\line(1,0){ 60}}
\put(  60, 0){\line(0,1){ 12}}}
\multiput(61,0)(0,0){1}{
\put(   0,  0){\line(1,0){ 80}}\put(   0,12){\line(1,0){ 80}}
\put(   0,  0){\line(0,1){ 12}}\put(  80, 0){\line(0,1){ 12}}}
\multiput(142,0)(0,0){1}{
\put(   0,  0){\line(1,0){ 51}}\put(   0,12){\line(1,0){ 51}}
\put(   0,  0){\line(0,1){ 12}}}
\put(  62, -2){$\underbrace{\makebox(78,15)[c]{}}_{}$}\put(  99,-17){$\seq{w}^*$}
\put(  27,  4){$\seq{w}^*$}
\put( 125,  4){$\seq{\gamma}^*$}
\put( 152,  4){$\seq{\alpha}^*$}
}
\multiput(156,15)(0,4){9}{\line(0,1){2}}
\multiput(181,15)(0,4){9}{\line(0,1){2}}
\multiput(206,15)(0,4){9}{\line(0,1){2}}
\end{picture}
\]
\normalsize
The shift-boundedness of $\seq{w}$ gives that $\seq{\alpha}\leq \seq{\gamma}^*$ and $\seq{\beta}<\seq{\alpha}^*$ and therefore $\seq{b}_k < \seq{e}_i[1,|\seq{b}_k|]$.

\vspace{1.5ex}
\textit{(2.)} Let us consider the inequality 
\begin{equation}
\label{eq: e(b)>b}
e(g_n(\seq{b}_k)) > \seq{b}_k,
\end{equation}
for $0<n<|\seq{b}_k|$. The inequality (\ref{eq: e(b)>b}) fails whenever $n$ is such that $(\seq{b}_k)_n=1$, hence we may assume that $(\seq{b}_k)_n=0$. For $2^i\leq n < |\seq{s}|$ we have as $\seq{s}$ is an $e_i$-minimal sequence that $\seq{b}_k[1,|\seq{s}|] \leq e(g_n(\seq{b}_k))[1,|\seq{s}|]$.
If the inequality is strict we are done, hence we only have to consider the case when having equality, $\seq{b}_k[1,|\seq{s}|] = e(g_n(\seq{b}_k))[1,|\seq{s}|]$. Let $r$ be the smallest integer such that $nr>|\seq{s}|$.

If $nr-|\seq{s}|>\frac{1}{2}n$ then let 
$\seq{z} = g_n(\seq{b}_k)$,
$\seq{\alpha} = \seq{b}_k[|\seq{s}|+1, 2|\seq{s}|-n(r-1)]$ and
$\seq{\beta}  = e(g_n(\seq{b}_k))[|\seq{s}|+1, 2|\seq{s}|-n(r-1)]$.

\footnotesize
\[
\begin{picture}(230,65)
\multiput(40,37)( 0,0){1}{
\put( -35,  0){\makebox(30,12)[r]{$e(g_n(\seq{b}_k))=\ldots$}}
\multiput(0,0)(0,0){1}{
\put(   0,  0){\line(1,0){ 60}}\put(   0,12){\line(1,0){ 60}}
\put(  60, 0){\line(0,1){ 12}}}
\multiput(61,0)(0,0){1}{
\put(   0,  0){\line(1,0){ 80}}\put(   0,12){\line(1,0){ 80}}
\put(   0,  0){\line(0,1){ 12}}\put(  80, 0){\line(0,1){ 12}}}
\multiput(142,0)(0,0){1}{
\put(   0,  0){\line(1,0){ 51}}\put(   0,12){\line(1,0){ 51}}
\put(   0,  0){\line(0,1){ 12}}}
\put(  62,  0){$\overbrace{\makebox(78,15)[c]{}}_{}$}\put( 99,24){$\seq{z}^*$}
\put(  27,  4){$\seq{z}^*$}
\put(  73,  4){$\seq{\alpha}^*$}
\put( 103,  4){$\seq{\beta}^*$}
\put( 165,  4){$\seq{z}^*$}
}
\multiput(40,15)( 0,0){1}{
\put(-35,  0){\makebox(30,12)[r]{$\seq{b}_k=\ldots$}}
\multiput(0,0)(0,0){1}{
\put(   0,  0){\line(1,0){ 90}}\put(   0,12){\line(1,0){ 90}}
\put(  90,  0){\line(0,1){ 12}}}
\multiput(91,0)(0,0){1}{
\put(   0,  0){\line(1,0){102}}\put(   0,12){\line(1,0){102}}
\put(   0,  0){\line(0,1){ 12}}}
\put(  92, -2){$\underbrace{\makebox(78,1)[c]{}}_{}$}\put( 129,-17){$\seq{z}$}
\put(  45,  4){$\seq{s}$}
\put( 103,  4){$\seq{\alpha}$}
\multiput(170,0)(0,4){3}{\line(0,1){2}}
}
\multiput(131,15)(0,4){9}{\line(0,1){2}}
\multiput(161,15)(0,4){9}{\line(0,1){2}}
\end{picture}
\]
\normalsize
If we assume that $n$ is the smallest integer such that (\ref{eq: e(b)>b}) does not hold then $\seq{z}$ is the $e_i$-minimal prefix of $\seq{b}_k$. But then $\seq{z}$ is shift-bounded and we must have $\seq{\beta}^*>\seq{\alpha}$, which contradicts that $\seq{z}$ is the $e_i$-minimal prefix of $\seq{b}_k$. Therefore we get $e(g_n(\seq{b}_k))[1,\seq{b}_k] > \seq{b}_k$.

If $nr - |\seq{s}| \leq \frac{1}{2}n$ then let
$\seq{z}        = g_n(\seq{b}_k)$,
$\seq{\alpha}   = \seq{b}_k[|\seq{s}|+1, nr]$,
$\seq{\gamma}^* = e(g_n(\seq{b}_k))[|\seq{s}|+1, nr]$ and
$\seq{\beta}    = \seq{b}_k[nr+1,2nr-|\seq{s}| ]$.

\footnotesize
\[
\begin{picture}(230,65)
\multiput(40,37)( 0,0){1}{
\put( -35,  0){\makebox(30,12)[r]{$e(g_n(\seq{b}_k))=\ldots$}}
\multiput(0,0)(0,0){1}{
\put(   0,  0){\line(1,0){ 60}}\put(   0,12){\line(1,0){ 60}}
\put(  60, 0){\line(0,1){ 12}}}
\multiput(61,0)(0,0){1}{
\put(   0,  0){\line(1,0){ 80}}\put(   0,12){\line(1,0){ 80}}
\put(   0,  0){\line(0,1){ 12}}\put(  80, 0){\line(0,1){ 12}}}
\multiput(142,0)(0,0){1}{
\put(   0,  0){\line(1,0){ 51}}\put(   0,12){\line(1,0){ 51}}
\put(   0,  0){\line(0,1){ 12}}}
\put(  62,  0){$\overbrace{\makebox(78,15)[c]{}}_{}$}\put(  99,24){$\seq{z}^*$}
\put(  27,  4){$\seq{z}^*$}
\put( 125,  4){$\seq{\gamma}^*$}
\put( 152,  4){$\seq{\alpha}^*$}
}
\multiput(40,15)( 0,0){1}{
\put(-35,  0){\makebox(30,12)[r]{$\seq{b}_k=\ldots$}}
\multiput(0,0)(0,0){1}{
\put(   0,  0){\line(1,0){115}}\put(   0,12){\line(1,0){115}}
\put( 115,  0){\line(0,1){ 12}}}
\multiput(116,0)(0,0){1}{
\put(   0,  0){\line(1,0){77}}\put(   0,12){\line(1,0){77}}
\put(   0,  0){\line(0,1){ 12}}}
\put(  55,  4){$\seq{s}$}
\put( 125,  4){$\seq{\alpha}$}
\put( 152,  4){$\seq{\beta}$}
}
\multiput(156,15)(0,4){9}{\line(0,1){2}}
\multiput(181,15)(0,4){9}{\line(0,1){2}}
\multiput(206,15)(0,4){9}{\line(0,1){2}}
\end{picture}
\]
\normalsize
If we again assume that $n$ is the smallest integer such that (\ref{eq: e(b)>b}) does not hold then $\seq{z}$ is the $e_i$-minimal prefix of $\seq{b}_k$. But then $\seq{z}$ is shift-bounded and we must have $\seq{\gamma}^*\geq \seq{\alpha}$ and $\seq{\alpha}^* > \seq{\beta}$, which contradicts that $\seq{z}$ is the $e_i$-minimal prefix of $\seq{b}_k$ and again $e(g_n(\seq{b}_k))[1,|\seq{b}_k|] > \seq{b}_k$.

For $n = r|\seq{s}|$ with $1<r<k$ we have that $e(g_n(\seq{b}_k))[1,|\seq{b}_k|] > \seq{b}_k$ holds as $(\tilde{\seq{s}})^*>\seq{s} \geq p(\seq{s})$. Note that we only have to consider those $n$ such that $(\seq{b}_k)_n = 0$.

For $j|\seq{s}| < n < j|\seq{s}|+|p(\seq{s})|$ where $0<j<k$ let
$\seq{\alpha}^* = e(g_n(\seq{b}_k))[n+1, (j+1)|\seq{s}|]$ and
$\seq{\gamma} = \seq{b}_k[n+1,(j+1)|\seq{s}|]$.

\footnotesize
\[
\begin{picture}(230,65)
\multiput(40,37)( 0,0){1}{
\put(-35,  0){\makebox(30,12)[r]{$e(g_n(\seq{b}_k))=\ldots$}}
\multiput(0,0)(0,0){1}{
\put(   0,  0){\line(1,0){ 90}}\put(   0,12){\line(1,0){90}}
\put(  90,  0){\line(0,1){ 12}}}
\multiput( 91,0)(0,0){1}{
\put(   0,  0){\line(1,0){102}}\put(   0,12){\line(1,0){102}}
\put(   0,  0){\line(0,1){ 12}}}
\put(  92,  0){$\overbrace{\makebox(78,15)[c]{}}_{}$}\put( 129, 24){$\seq{s}^*$}
\put( 112,  4){$\seq{\alpha}^*$}
\multiput(170,2)(0,4){3}{\line(0,1){2}}
}
\multiput(40,15)( 0,0){1}{
\put( -35,  0){\makebox(30,12)[r]{$\seq{b}_k=\ldots$}}
\multiput(0,0)(0,0){1}{
\put(   0,  0){\line(1,0){ 60}}\put(   0,12){\line(1,0){ 60}}
\put(  60, 0){\line(0,1){ 12}}}
\multiput(61,0)(0,0){1}{
\put(   0,  0){\line(1,0){ 80}}\put(   0,12){\line(1,0){ 80}}
\put(   0,  0){\line(0,1){ 12}}\put(  80, 0){\line(0,1){ 12}}}
\multiput(142,0)(0,0){1}{
\put(   0,  0){\line(1,0){ 51}}\put(   0,12){\line(1,0){ 51}}
\put(   0,  0){\line(0,1){ 12}}}
\put(  62, -2){$\underbrace{\makebox(78,1)[c]{}}_{}$}\put( 97,-17){$\seq{s}$}
\put(  27,  4){$\seq{s}$}
\put( 112,  4){$\seq{\gamma}$}
}
\multiput(131,15)(0,4){9}{\line(0,1){2}}
\multiput(181,15)(0,4){9}{\line(0,1){2}}
\end{picture}
\]
\normalsize
From the shift-boundedness of $\seq{s}$ and the definition of $p(\seq{s})$ we have that $\seq{\alpha}^*>\seq{\gamma}$, which implies $e(g_n(\seq{b}_k))[1,|\seq{b}_k|] > \seq{b}_k$.

For $j|\seq{s}| +|p(\seq{s})| \leq n < (j+1)|\seq{s}|$ where $0<j<k$ let
$\seq{\alpha}^* = e(g_n(\seq{b}_k))[n+1, (j+1)|\seq{s}|]$,
$\seq{\gamma} = \seq{b}_k[n+1,(j+1)|\seq{s}|]$ and
$\seq{\beta}^* = e(g_n(\seq{b}_k))[(j+1)|\seq{s}|+1,(j+1)|\seq{s}|-n]$.

\footnotesize
\[
\begin{picture}(230,65)
\multiput(40,37)( 0,0){1}{
\put(-35,  0){\makebox(30,12)[r]{$e(g_n(\seq{b}_k))=\ldots$}}
\multiput(0,0)(0,0){1}{
\put(   0,  0){\line(1,0){115}}\put(   0,12){\line(1,0){115}}
\put( 115,  0){\line(0,1){ 12}}}
\multiput(116,0)(0,0){1}{
\put(   0,  0){\line(1,0){77}}\put(   0,12){\line(1,0){77}}
\put(   0,  0){\line(0,1){ 12}}}
\put( 125,  4){$\seq{\alpha}^*$}
\put( 152,  4){$\seq{\beta}^*$}
}
\multiput(40,15)( 0,0){1}{
\put( -35,  0){\makebox(30,12)[r]{$\seq{b}_k=\ldots$}}
\multiput(0,0)(0,0){1}{
\put(   0,  0){\line(1,0){ 60}}\put(   0,12){\line(1,0){ 60}}
\put(  60, 0){\line(0,1){ 12}}}
\multiput(61,0)(0,0){1}{
\put(   0,  0){\line(1,0){ 80}}\put(   0,12){\line(1,0){ 80}}
\put(   0,  0){\line(0,1){ 12}}\put(  80, 0){\line(0,1){ 12}}}
\multiput(142,0)(0,0){1}{
\put(   0,  0){\line(1,0){ 51}}\put(   0,12){\line(1,0){ 51}}
\put(   0,  0){\line(0,1){ 12}}}
\put(  62, -2){$\underbrace{\makebox(78,15)[c]{}}_{}$}\put(  99,-17){$\seq{s}$}
\put(  27,  4){$\seq{s}$}
\put( 125,  4){$\seq{\gamma}$}
\put( 152,  4){$\seq{\alpha}$}
}
\multiput(156,15)(0,4){9}{\line(0,1){2}}
\multiput(181,15)(0,4){9}{\line(0,1){2}}
\multiput(206,15)(0,4){9}{\line(0,1){2}}
\end{picture}
\]
\normalsize
Again, as $\seq{s}$ is shift-bounded we have $\seq{\alpha}^*\seq{\beta}^* > \seq{\gamma}\seq{\alpha}$ since $\seq{\alpha}^*\geq \seq{\gamma}$ and $\seq{\beta}^*> \seq{\alpha}$, which gives $e(g_n(\seq{b}_k))[1,|\seq{b}_k|] > \seq{b}_k$.

For $k|\seq{s}|<n<|\seq{b}_k|$ let
$\seq{\alpha}^* = e(g_n(\seq{b}_k))[n+1, |\seq{b}_k|]$ and
$\seq{\gamma}   = \seq{b}_k[n+1,|\seq{b}_k|]$.

\footnotesize
\[
\begin{picture}(230,65)
\multiput(40,37)( 0,0){1}{
\put( -35,  0){\makebox(30,12)[r]{$e(g_n(\seq{b}_k))=\hdots$}}
\multiput(0,0)(0,0){1}{
\put(   0,  0){\line(1,0){130}}\put(  0,12){\line(1,0){130}}
\put( 130,  0){\line(0,1){ 12}}}
\multiput(131,0)(0,0){1}{
\put(   0,  0){\line(1,0){62}}\put(  0,12){\line(1,0){62}}
\put(   0,  0){\line(0,1){ 12}}}
\put( 148,  4){$\seq{\alpha}^*$}
}
\multiput(40,15)( 0,0){1}{
\put(-35,  0){\makebox(30,12)[r]{$\seq{b}_k=\hdots$}}
\multiput(0,0)(0,0){1}{
\put(   0,  0){\line(1,0){ 50}}\put(   0,12){\line(1,0){ 50}}
\put(  50,  0){\line(0,1){ 12}}}
\multiput(51,0)(0,0){1}{
\put(   0,  0){\line(1,0){ 70}}\put(   0,12){\line(1,0){ 70}}
\put(   0,  0){\line(0,1){ 12}}\put(  70, 0){\line(0,1){ 12}}}
\multiput(122,0)(0,0){1}{
\put(   0,  0){\line(1,0){ 50}}\put(   0,12){\line(1,0){ 50}}
\put(   0,  0){\line(0,1){ 12}}\put(  50, 0){\line(0,1){ 12}}}
\put( 122, -2){$\underbrace{\makebox(48,1)[c]{}}_{}$}\put( 137,-17){$p(\seq{s})$}
\put(  22,  4){$\seq{s}$}
\put(  83,  4){$\seq{s}$}
\put( 148,  4){$\seq{\gamma}$}
}
\multiput(171,15)(0,4){9}{\line(0,1){2}}
\multiput(212,15)(0,4){9}{\line(0,1){2}}
\end{picture}
\]
\normalsize
Since $p(\seq{s})$ is shift-bounded we have that $\seq{\alpha}^*\geq \seq{\gamma}$, and therefore we obtain $e(g_n(\seq{b}_k))[1,|\seq{b}_k|] \geq \seq{b}_k$.

\vspace{1.5ex}
\textit{(3.)} 
The shift-boundedness of $\seq{s}$ gives that $g_n(\seq{b}_k)^{\infty} \leq \seq{s}^{\infty} < \seq{b}_k$ for all $2^i\leq n < |\seq{b}_k|$ such that $(\seq{b}_k)_n=1$. Hence the $e_i$-minimality of $\seq{b}_k$ follows by combining this with the previous statements of this lemma.
\end{proof}

\section{The Set $\mathcal{A}$}

Let $\seq{u}$ be a finite sequence ending with a 1 and let $\Ac{\seq{u}}\subset\Sigma_2$ be the set of infinite sequences created from the finite sequences $\tilde{\seq{u}}$, $\seq{u}$, $\seq{u}^*$, and $\seq{u}'$ following the transition matrix
\[
A = \left(\begin{array}{cccc}0 & 0 & 1 & 1 \\1 & 1 & 0 & 0 \\0 & 0 & 1 & 1 \\1 & 1 & 0 & 0\end{array}\right),
\]
where the rows and columns are ordered in the order $\tilde{\seq{u}},\seq{u},\seq{u}^*,\seq{u}'$. The elements of $\Ac{\seq{u}}$ are sequences similar to the suffix given in (\ref{eq: uuu}) but where the upper bounding $k$ has been removed. Note that the transition matrix $A$ is primitive and has the spectral radius $\rho(A) = 2$. By Proposition \ref{prop: holder-dim} we have

\begin{lemma}
\label{lemma: dim A}
Let $\seq{u}$ be a finite sequence ending with a 1. Then $\dim_H \Ac{\seq{u}} = \frac{1}{|\seq{u}|}$.
\end{lemma}

For the special case when having $\seq{u}=f^i(1)$ we have $\Ac{1} = \Sigma_2^{}$ for $i=0$ and for $i\geq0$ we have the nested inclusions
\begin{equation}	
\label{eq. A(i+1) < Ai}
\Ac{f^{i+1}(1)}\subset \Ac{f^{i}(1)}.
\end{equation}

\begin{definition}
Let $\seq{u}$ be a finite sequence ending with a 1 and let 
\[
\mu_{\seq{u}}: \Ac{\seq{u}}\to  \Sigma_2^{}
\]
be the map $(\tilde{\seq{u}},\seq{u},\seq{u}^*, \seq{u}') \mapsto (0,1,0,1)$. Let $\mu_{\seq{u}}^{-1}$ map the first 0 in  each block of zeros to $\tilde{\seq{u}}$ otherwise 0 is mapped to $\seq{u}^*$ and let the first 1 in each block of ones be mapped to $\seq{u}'$ otherwise 1 mapped to $\seq{u}$.
\end{definition}

Note that $\mu$ could equally have been defined as a function between sets of finite sequences, that is, $\mu_{\seq{u}}:\{ \seq{x} [1,n|\seq{u}|] : \seq{x} \in \Ac{\seq{u}}\} \to \{ \seq{x} [1,n] : \seq{x} \in \Sigma_2^{} \}$. 

\vspace{1.5ex}
A function $T$ similar to our $\mu_{\seq{u}}$ is defined by Allouche in \cite{allouche83:1}. The function $T$ is there used to show that the set $\Gamma$, (see (\ref{eq: allouche set})), is self similar. 

\vspace{1.5ex}

The function $\mu_{\seq{u}}: \Ac{\seq{u}}\to\Sigma_2^{}$ is not bijective as for $U_1 = [\tilde{\seq{u}}] \cap \Ac{\seq{u}}$ and $U_2 = [\seq{u}'] \cap \Ac{\seq{u}}$ we have $\mu_{\seq{u}}^{-1}(\Sigma_2^{}) = U_1 \cup U_2$, where the right-hand-side is a proper subset of $\Ac{\seq{u}}$ if $|\seq{u}|>1$. The violation of the bijectivity is however only in the first positions, so by shifting these out we have
\[
\sigma^{|\seq{u}|}\left(\mu_{\seq{u}}^{-1}(\Sigma_2^{})\right) = \Ac{\seq{u}}.
\]
If we restrict $\mu$ to map sequences from $[\tilde{\seq{u}}] \cap \Ac{\seq{u}}$ into $[0]$ we obtain a bijection as the sequences causing a collisions due to the definition of the inverse of $\mu$ have been removed. 

\begin{example}
Let $\seq{u} = 01$.  Then $\seq{c}_1 = 0011$ and  $\seq{c}_2 = 1011$ are prefixes of sequences in $\Ac{\seq{u}}$. We have $\mu_{01}(0011) = 01$ and $\mu_{01}(1011) = 01$, but $\mu_{01}^{-1}(01) = 0011$.
\end{example}

\begin{lemma}
\label{lemma: q preserves order}
Let $\seq{u}$ be a finite sequence not ending with $0$ and such that $\seq{u}\leq\seq{u}'$. Put $U = [\tilde{\seq{u}}]\cap \Ac{\seq{u}}$ and $V = [0]$, (or $U = ([\tilde{\seq{u}}]\cap \Ac{\seq{u}}) [1,|\seq{u}|n]$ and $V = ([0])[1,n]$ in the case of finite sequeces). Then $\mu_{\seq{u}}:U\to V$ is bijective and order-preserving.
\end{lemma}

\begin{proof}
The bijectivity is clear from the just above reasoning of the definition of the inverse of $\mu$. For the order preservation let $\seq{c}_1<\seq{c}_2$ be two sequences in $U$ and let $\seq{s}_1 = \mu(\seq{c}_1)$ and $\seq{s}_2 = \mu(\seq{c}_2)$. Assume for contradiction that $\seq{s}_1> \seq{s}_2$. There is a smallest $n$ such that $(\seq{s}_1)_n = 1$ and $(\seq{s}_2)_n = 0$. Let $\seq{w}_1 = \seq{c}_1[n|\seq{u}|+1,(n+1)|\seq{u}|]$ and $\seq{w}_2 = \seq{c}_2[n|\seq{u}|+1,(n+1)|\seq{u}|]$. That is, $\seq{w}_1$ is the subsequence in $\seq{c}_1$ mapped into $(\seq{s}_1)_n$ by $\mu_{\seq{u}}$, and similarly for $\seq{w}_2$. If $(\seq{s}_1)_{n-1} = 1$ then $\seq{w}_1 = \seq{u}$ and $\seq{w}_2 = \tilde{\seq{u}}$, contradicting $\seq{c}_1 < \seq{c_2}$. If $(\seq{s}_1)_{n-1} = 0$ then $\seq{w}_1 = \seq{u}'$ and $\seq{w}_2 = \seq{u}^*$, again contradicting $\seq{c}_1 < \seq{c_2}$.  Finally, if $n = 1$ then $\seq{w}_1 = \seq{u}'$ and $\seq{w}_2 = \tilde{\seq{u}}$, then similarly this would imply $\seq{c}_1 < \seq{c_2}$.
\end{proof}

\begin{lemma}
\label{lemma: length of SB c}
For $\seq{u} = f^{i-1}(1)$ where $i\geq1$ let $U = [\tilde{\seq{u}}]\cap \Ac{\seq{u}}$. If $\seq{c}$ is an infinite  shift-bounded $e_i$-sequence then $\seq{c} \in U$. If $\seq{c}$ is a  finite shift-bounded $e_i$-sequence then $\seq{c}$ is a prefix of a sequences in $U$ and $|\seq{c}| = k\cdot 2^{i-1}$ for some $k\geq 3$. 
\end{lemma}

\begin{proof}
There is a maximal $N$ and a sequence $\seq{v}$ such that $\seq{c} = \tilde{\seq{u}}(\seq{u}^*)^N \seq{v}$ with $\seq{v}>\seq{u}^*$ as an $e_i$-sequence must start with $\tilde{\seq{u}}(\seq{u}^*)$. By shifting $n = (1+N)|\seq{u}|$ times we obtain from $\seq{c}'= \seq{u}'\seq{u}^N\seq{v}' > \sigma^n(\seq{c}) = \seq{v}'$ that $\seq{u}'$ must be a prefix of $\seq{v}$. Hence Lemma \ref{lemma: uuu} gives that $\seq{c}\in U$ if $\seq{c}$ is infinite or that $\seq{c}$ is the prefix of a sequence in $U$ if $\seq{c}$ is finite. Moreover, since $\seq{u}'$ is a prefix of $\seq{v}$ we have that $|\seq{c}|\geq 3 |\seq{u}|$. 

For the length of $\seq{c}$ in the finite case we have to show that $\seq{u}$ and $\seq{u}'$ are the only allowed suffixes of $\seq{c}$ of length $|\seq{u}|$ and moreover we may not find $\seq{u}$ or $\seq{u}'$ by cutting an ending $\tilde{\seq{u}}$, $\seq{u}$, $\seq{u}^*$ or $\seq{u}'$ off.

The sequence $\seq{c}$ cannot end with $\tilde{\seq{u}}$ or $\seq{u}^*$ as it then would end with a zero, contradicting $\seq{c}$ being shift-bounded. 

If $\seq{c}$ ends with a prefix $\seq{v}$ of $\tilde{\seq{u}}$ then there is an $n$ such that $\sigma^n(\seq{c})=\seq{v} \leq \tilde{\seq{u}} < \seq{c}$, contradicting the shift-boundedness of $\seq{c}$. The same procedure holds for a proper prefix of $\seq{u}$.

If $\seq{c}$ ends with a proper prefix $\seq{v}$ of $\seq{u}^*$ then $\seq{c}$ must end with $\tilde{\seq{u}}(\seq{u}^*)^m \seq{v}$ for some $0\leq m$, as $\seq{c}$ is prefix of a sequence in $U$. If $\seq{c} = \tilde{\seq{u}}(\seq{u}^*)^m \seq{v}$ then it would not be an $e_i$-sequence. Hence $\seq{c}$ must end with 
\[
\tilde{\seq{u}} (\seq{u}^*)^r  \seq{u}' (\seq{u})^s  \tilde{\seq{u}} (\seq{u}^*)^m \seq{v}.
\]
But then for $n = (2+r+s)|\seq{u}|$ we have $\sigma^n(\seq{c})=\tilde{\seq{u}}(\seq{u}^*)^m \seq{v}  < \seq{c}$, contradicting the shift-boundedness of $\seq{c}$. 
\end{proof}

\begin{lemma}
\label{lemma: e21-min/SB <=> e2i-min/SB}
For $\seq{u} = f^{i-1}(1)$ where $i\geq1$ let $U = [\tilde{\seq{u}}]\cap \Ac{\seq{u}}$ and $V = [0]$, (or $U = ([\tilde{\seq{u}}]\cap \Ac{\seq{u}}) [1,|\seq{u}|n]$ and $V = ([0])[1,n]$ in the case of finite sequences). Then $\mu_{\seq{u}}:U\to V$ is a bijection between $e_i$-minimal sequences and $e_1$-minimal sequences.
\end{lemma}

\begin{proof}
It is clear that an $e_1$-minimal sequence is a prefix of a sequence in $V$ and by Lemma \ref{lemma: length of SB c} an $e_i$-minimal sequence is a prefix of a sequences in $U$.

Let $\seq{c}$ be a prefix of a sequence in $U$ such that $|\seq{c}| = k|\seq{u}|$, for some $k\geq 3$, and let $\seq{s} = \mu_{\seq{u}}(\seq{c})$. 
Since $g_{n|\seq{u}|}(\seq{c})$ ends with either $\seq{u}$ or $\seq{u}'$ and begins with $\tilde{\seq{u}}$ it follows that $g_{n|\seq{u}|}(\seq{c})^{\infty}$ is an element in $U$. 
Similarly we have that $e(g_{n|\seq{u}|}(\seq{c}))[1,n|\seq{u}|]$ ends with either $\tilde{\seq{u}}$ or $\seq{u}^*$ and  since $(g_{n|\seq{u}|}(\seq{c}))^*$ begins with $\seq{u}'$ we have that also $e(g_{n|\seq{u}|}(\seq{c}))$ is an element of $U$.

Lemma \ref{lemma: length of SB c} gives that we only have to check for minimality of $\seq{c}$ in prefixes of length $k|\seq{u}|$ for $k\geq 3$. Assume there is an $n\geq 2 |\seq{u}| = 2^i$ such that
\[	
e(g_{n|\seq{u}|}(\seq{c})) \leq \seq{c} \leq g_{n|\seq{u}|}(\seq{c})^{\infty}
\]
does not hold. Then the order-preservation of $\mu$ gives that 
\[	
e(g_{n}(\seq{s})) \leq \seq{s} \leq g_{n}(\seq{s})^{\infty}
\]
cannot hold either. 
\end{proof}

\begin{theorem}
\label{thm: self-similar}
For $\seq{u} = f^{i-1}(1)$ for $i\geq 1$ put $U = \mathcal{A}(\seq{u})$ and $V = \Sigma_2$, (or put $U = \mathcal{A}(\seq{u})[1,|\seq{u}|n]$ and $V=\Sigma_2[1,n]$ for the finite case). Let $\mu_{\seq{u}}:U \to V$ and let $\seq{c}$ be an $e_i$-sequence such that $\mu_{\seq{u}}(c)$ is well defined. Then 
\[
\dim_H F(\seq{c}) =  \frac{1}{2^{i-1}} \dim_H F\big(\mu_{\seq{u}}(\seq{c})\big).
\]
\end{theorem}

\begin{proof}
Let $S = ([\tilde{\seq{u}}]\cup[\seq{u}'])\cap F(\seq{c})$. By the order-preservation of $\mu$ and Lemma \ref{lemma: uuu} we have $\mu_{\seq{u}}(S) = F(\mu_{\seq{u}}(\seq{c}))$. Moreover $\mu_{\seq{u}}^{-1}(\mu_{\seq{u}}(S)) \subset F(\seq{c})$. Hence 
\[
\frac{1}{2^{i-1}}\dim_H F(\mu_{\seq{u}}(\seq{c})) \leq \dim_H F(\seq{c}).
\]

For the reversed inequality, let $\seq{x} \in F(\seq{c})$. If $\seq{x}$ does not contain $00$ nor $11$ then $\seq{x}$ is either of the sequences $(01)^{\infty}$ or $(10)^{\infty}$. If $\seq{x}$ does contain two consecutive zeros or ones then Lemma \ref{lemma: uuu} gives that $\seq{x}$ ends with a sequence which is an element in 
$\Ac{f^1(1)}$. Therefore, and by the use of the nested inclusion (\ref{eq. A(i+1) < Ai}), we have that 
\[
\bigcup_{n=1}^{\infty}\left\{ \seq{v}[1,n] \seq{w} : \seq{v}\in F(f^i(1)),\, \seq{w}\in \bigcup_{k=|\seq{u}|}^{2|\seq{u}|-1} \sigma^k\Big(\mu_{\seq{u}}^{-1}\big(F(\mu_{\seq{u}}(\seq{c}))\big)\Big)\right\},
\]
contains $F(\seq{c})$, which implies the desired inequality.
\end{proof}

\begin{corollary}
\label{cor: dim e_i}
Let $i\geq1$. Then $\dim_H F(\seq{e}_i) = \frac{1}{2^{i}}$.
\end{corollary}

From Corollary \ref{cor: dim e_i} we can directly derive Moreira's Theorem \ref{thm: moreira}.

\section{Results and Proofs}

Let us define the interval $I(\seq{c})$ to be the set of sequences 
\[
I(\seq{c}) = \left\{ \seq{x} \in \Sigma_2^{}: e(g_{m_{\seq{c}}}(\seq{c})) \leq \seq{x} \leq  (g_{m_{\seq{c}}}(\seq{c}))^{\infty} \right\}.
\]
Note that if $\seq{c}$ is an infinite $e_i$-minimal sequence the interval $I(\seq{c})$ will only contain one element, $I(\seq{c}) = \{\seq{c}\}$. We have to show that the definition of the interval $I(\seq{c})$ is independent of the choice of the representative $\seq{c}$.

\begin{lemma}
\label{lemma: c_ma=1}
Let $\seq{c}$ be a finite $e_i$-minimal sequence and $\seq{a}\in I(\seq{c})$. If $m_{\seq{a}}\leq m_{\seq{c}}$ then $(\seq{c})_{m_{\seq{a}}}=1$.
\end{lemma}

\begin{proof}
Assume for contradiction that $(\seq{c})_{m_{\seq{a}}}=0$. Form the assumption we have $m_{\seq{a}} < m_{\seq{c}}$, as $\seq{c}$ ends with a $1$. If $e(g_{m_{\seq{a}}}(\seq{a}))< \seq{c}$ then we have $e(g_{m_{\seq{a}}}(\seq{a}))< \seq{c} < g_{m_{\seq{a}}}(\seq{a})^{\infty}$, which contradicts the $e_i$-minimality of $\seq{c}$. Hence we must have $\seq{c} < e(g_{m_{\seq{a}}}(\seq{a}))$. If $\seq{c}$ is not a proper prefix of $e(g_{m_{\seq{a}}}(\seq{a}))$ then $\seq{c}^{\infty} < e(g_{m_{\seq{a}}}(\seq{a}))<\seq{a}$, which contradicts that $\seq{a}\in I(\seq{c})$. 

If $\seq{c}$ is a proper prefix of $e(g_{m_{\seq{a}}}(\seq{a}))$ let $k\geq 1 $ be the largest integer such that $k\,m_{\seq{a}} \leq m_{\seq{c}}$. If $k\,m_{\seq{a}} = m_{\seq{c}}$ then we have that $\seq{c}^{\infty} < e(g_{m_{\seq{a}}}(\seq{a}))$ as 
$\seq{c}[1,m_{\seq{a}}] < (g_{m_{\seq{a}}}(\seq{a}))^*$, which implies $\seq{a}\notin I(\seq{c})$, a contradiction. 

\footnotesize
\[
\begin{picture}(230,65)
\multiput(40,37)( 0,0){1}{
\put(-35,  0){\makebox(30,12)[r]{$e(g_{m_{\seq{a}}}(\seq{a}))=\hdots$}}
\multiput(0,0)(0,0){1}{
\put(   0,  0){\line(1,0){ 24}}\put(   0,12){\line(1,0){ 24}}
\put(  24,  0){\line(0,1){ 12}}}
\multiput(25,0)(66,0){2}{
\put(   0,  0){\line(1,0){ 65}}\put(   0,12){\line(1,0){ 65}}
\put(   0,  0){\line(0,1){ 12}}\put(  65, 0){\line(0,1){ 12}}
\put(  20,  4){$g_{m_{\seq{a}}}(\seq{a})^*$}}
\multiput(157,0)(0,0){1}{
\put(   0,  0){\line(1,0){ 36}}\put(   0,12){\line(1,0){ 36}}
\put(   0,  0){\line(0,1){ 12}}}
}
\multiput(40,15)( 0,0){1}{
\put( -35,  0){\makebox(30,12)[r]{$\seq{c}^{\infty}=\hdots$}}
\multiput(0,0)(0,0){1}{
\put(   0,  0){\line(1,0){ 90}}\put(  0,12){\line(1,0){ 90}}
\put(  90,  0){\line(0,1){ 12}}
\put(  45,  4){$\seq{c}$}}
\multiput( 91,0)(0,0){1}{
\put(   0,  0){\line(1,0){102}}\put(  0,12){\line(1,0){102}}
\put(   0,  0){\line(0,1){ 12}}
\put(  45,  4){$\seq{c}$}}
}
\multiput(131,15)(0,4){11}{\line(0,1){2}}
\put(125,60){$k\,m_{\seq{a}}$}
\multiput(197,15)(0,4){9}{\line(0,1){2}}
\end{picture}
\]
\normalsize

If $ m_{\seq{c}} - k\,m_{\seq{a}} < \frac{1}{2}  m_{\seq{a}}$ let 
$\seq{\alpha}  = \seq{c}^{\infty}[m_{\seq{c}}+1,2m_{\seq{c}}-k\,m_{\seq{a}}]$ and
$\seq{\beta}^* = e(g_{m_{\seq{a}}}(\seq{a}))[m_{\seq{c}}+1,2m_{\seq{c}}-k\,m_{\seq{a}}]$.

\footnotesize
\[
\begin{picture}(230,65)
\multiput(40,37)( 0,0){1}{
\put(-35,  0){\makebox(30,12)[r]{$e(g_{m_{\seq{a}}}(\seq{a}))=\hdots$}}
\multiput(0,0)(0,0){1}{
\put(   0,  0){\line(1,0){ 14}}\put(   0,12){\line(1,0){ 14}}
\put(  14,  0){\line(0,1){ 12}}}
\multiput(15,0)(66,0){2}{
\put(   0,  0){\line(1,0){ 65}}\put(   0,12){\line(1,0){ 65}}
\put(   0,  0){\line(0,1){ 12}}\put(  65, 0){\line(0,1){ 12}}}
\put(  33,  4){$g_{m_{\seq{a}}}(\seq{a})^*$}
\put(  88,  4){$\seq{\alpha}^*$}
\put( 108,  4){$\seq{\beta}^*$}
\put(  82,  0){$\overbrace{\makebox(63,15)[c]{}}_{}$}\put( 104, 24){$g_{m_{\seq{a}}}(\seq{a})^*$}
\multiput(147,0)(0,0){1}{
\put(   0,  0){\line(1,0){ 46}}\put(   0,12){\line(1,0){ 46}}
\put(   0,  0){\line(0,1){ 12}}}
}
\multiput(40,15)( 0,0){1}{
\put( -35,  0){\makebox(30,12)[r]{$\seq{c}^{\infty}=\hdots$}}
\multiput(0,0)(0,0){1}{
\put(   0,  0){\line(1,0){100}}\put(  0,12){\line(1,0){100}}
\put( 100,  0){\line(0,1){ 12}}
\put(  50,  4){$\seq{c}$}}
\multiput(101,0)(0,0){1}{
\put(   0,  0){\line(1,0){ 92}}\put(  0,12){\line(1,0){ 92}}
\put(   0,  0){\line(0,1){ 12}}
\put(  40,  4){$\seq{c}$}}
\put( 108,  4){$\seq{\alpha}$}
}
\multiput(141,15)(0,4){9}{\line(0,1){2}}
\multiput(161,15)(0,4){9}{\line(0,1){2}}
\end{picture}
\]
\normalsize
The shift-boundedness of $g_{m_{\seq{a}}}(\seq{a})$ gives that $\seq{\alpha}<\seq{\beta}^*$ and therefore $\seq{c}^{\infty} < e(g_{m_{\seq{a}}}(\seq{a}))$, which implies $\seq{a}\notin I(\seq{c})$, a contradiction.

If $ m_{\seq{c}} - k\,m_{\seq{a}} \geq \frac{1}{2}  m_{\seq{a}}$ let 
$\seq{\alpha}  = \seq{c}^{\infty}[m_{\seq{c}}+1, (k+1)m_{\seq{a}}]$, 
$\seq{\beta} = \seq{c}^{\infty}[(k+1)m_{\seq{a}}+1, 2(k+1)m_{\seq{a}}-m_{\seq{c}}]$ and
$\seq{\gamma}^* = e(g_{m_{\seq{a}}}(\seq{a}))[m_{\seq{c}}+1, (k+1)m_{\seq{a}}]$.

\footnotesize
\[
\begin{picture}(230,65)
\multiput(40,37)( 0,0){1}{
\put(-35,  0){\makebox(30,12)[r]{$e(g_{m_{\seq{a}}}(\seq{a}))=\hdots$}}
\multiput(0,0)(0,0){1}{
\put(   0,  0){\line(1,0){ 14}}\put(   0,12){\line(1,0){ 14}}
\put(  14,  0){\line(0,1){ 12}}}
\multiput(15,0)(66,0){2}{
\put(   0,  0){\line(1,0){ 65}}\put(   0,12){\line(1,0){ 65}}
\put(   0,  0){\line(0,1){ 12}}\put(  65, 0){\line(0,1){ 12}}}
\put(  33,  4){$g_{m_{\seq{a}}}(\seq{a})^*$}
\put( 133,  4){$\seq{\gamma}^*$}
\put( 153,  4){$\seq{\alpha}^*$}
\put(  82,  0){$\overbrace{\makebox(63,15)[c]{}}_{}$}\put( 104, 24){$g_{m_{\seq{a}}}(\seq{a})^*$}
\multiput(147,0)(0,0){1}{
\put(   0,  0){\line(1,0){ 46}}\put(   0,12){\line(1,0){ 46}}
\put(   0,  0){\line(0,1){ 12}}}
}
\multiput(40,15)( 0,0){1}{
\put( -35,  0){\makebox(30,12)[r]{$\seq{c}^{\infty}=\hdots$}}
\multiput(0,0)(0,0){1}{
\put(   0,  0){\line(1,0){125}}\put(  0,12){\line(1,0){125}}
\put( 125,  0){\line(0,1){ 12}}
\put(  50,  4){$\seq{c}$}}
\multiput(126,0)(0,0){1}{
\put(   0,  0){\line(1,0){ 67}}\put(  0,12){\line(1,0){ 67}}
\put(   0,  0){\line(0,1){ 12}}}
\put( 133,  4){$\seq{\alpha}$}
\put( 153,  4){$\seq{\beta}$}
}
\multiput(166,15)(0,4){9}{\line(0,1){2}}
\multiput(186,15)(0,4){9}{\line(0,1){2}}
\multiput(206,15)(0,4){9}{\line(0,1){2}}
\end{picture}
\]
\normalsize
The shift-boundedness of $g_{m_{\seq{a}}}(\seq{a})$ and $\seq{c}$ gives that $\seq{\gamma}^*\seq{\alpha}^* > \seq{\alpha}\seq{\beta}$ and therefore $\seq{c}^{\infty} < e(g_{m_{\seq{a}}}(\seq{a}))$, which implies $\seq{a}\notin I(\seq{c})$ again a contradiction. 
\end{proof}

\begin{theorem}
For any $\seq{a} \in I(\seq{c})$ we have $I(\seq{a}) = I(\seq{c})$.
\end{theorem}

\begin{proof}
We may assume that $\seq{c}$ is a finite $e_i$-minimal sequence.
If $m_{\seq{c}} > m_{\seq{a}}$ then Lemma \ref{lemma: c_ma=1} gives that $(\seq{c})_{m_{\seq{a}}} = 1$ and then Lemma \ref{lemma: prefix is minimal} gives $g_{m_{\seq{a}}}(\seq{a})^{\infty}<\seq{c}$. 
Thus we must have $g_{m_{\seq{a}}}(\seq{a})^{\infty}[1, m_{\seq{c}}] \leq e(\seq{c})[1,m_{\seq{c}}]$. If the inequality is strict then $\seq{a}\notin I(\seq{c})$, which contradicts our assumption. Hence we only have to consider the case when having equality. 

Let $k\geq1$ be that largest integer such that $k\,m_{\seq{a}} \leq m_{\seq{c}}$. If $k\,m_{\seq{a}} = m_{\seq{c}}$ then $g_{m_{\seq{a}}}(\seq{a})^{\infty}<e(\seq{c})$, since $g_{m_{\seq{a}}}(\seq{a})< c^*[1,{m_{\seq{a}}}]$, which contradicts $\seq{a}\in I(\seq{c})$. 

\footnotesize
\[
\begin{picture}(230,65)
\multiput(40,37)( 0,0){1}{
\put( -35,  0){\makebox(30,12)[r]{$e(\seq{c})=\hdots$}}
\multiput(0,0)(0,0){1}{
\put(   0,  0){\line(1,0){ 90}}\put(  0,12){\line(1,0){ 90}}
\put(  90,  0){\line(0,1){ 12}}
\put(  45,  4){$\tilde{\seq{c}}$}}
\multiput( 91,0)(0,0){1}{
\put(   0,  0){\line(1,0){102}}\put(  0,12){\line(1,0){102}}
\put(   0,  0){\line(0,1){ 12}}
\put(  45,  4){$\seq{c}^*$}}
}
\multiput(40,15)( 0,0){1}{
\put(-35,  0){\makebox(30,12)[r]{$g_{m_{\seq{a}}}(\seq{a})^{\infty}=\hdots$}}
\multiput(0,0)(0,0){1}{
\put(   0,  0){\line(1,0){ 24}}\put(   0,12){\line(1,0){ 24}}
\put(  24,  0){\line(0,1){ 12}}}
\multiput(25,0)(66,0){2}{
\put(   0,  0){\line(1,0){ 65}}\put(   0,12){\line(1,0){ 65}}
\put(   0,  0){\line(0,1){ 12}}\put(  65, 0){\line(0,1){ 12}}
\put(  20,  4){$g_{m_{\seq{a}}}(\seq{a})$}}
\multiput(157,0)(0,0){1}{
\put(   0,  0){\line(1,0){ 36}}\put(   0,12){\line(1,0){ 36}}
\put(   0,  0){\line(0,1){ 12}}}
}
\multiput(131,15)(0,4){11}{\line(0,1){2}}
\put(125,60){$k\,m_{\seq{a}}$}
\multiput(197,15)(0,4){9}{\line(0,1){2}}
\end{picture}
\]
\normalsize

If $ m_{\seq{c}} - k\,m_{\seq{a}} < \frac{1}{2} m_{\seq{a}}$ let 
$\seq{\alpha}^* = e(\seq{c})[m_{\seq{c}}+1,  2m_{\seq{c}}-k\,m_{\seq{a}}]$ and 
$\seq{\beta} = g_{m_{\seq{a}}}(\seq{a})^{\infty}[m_{\seq{c}}+1,  2m_{\seq{c}}-k\,m_{\seq{a}}]$.

\footnotesize
\[
\begin{picture}(230,65)
\multiput(40,37)( 0,0){1}{
\put( -35,  0){\makebox(30,12)[r]{$e(\seq{c})=\hdots$}}
\multiput(0,0)(0,0){1}{
\put(   0,  0){\line(1,0){100}}\put(  0,12){\line(1,0){100}}
\put( 100,  0){\line(0,1){ 12}}
\put(  50,  4){$\tilde{\seq{c}}$}}
\multiput(101,0)(0,0){1}{
\put(   0,  0){\line(1,0){ 92}}\put(  0,12){\line(1,0){ 92}}
\put(   0,  0){\line(0,1){ 12}}
\put(  40,  4){$\seq{c}^*$}}
\put( 108,  4){$\seq{\alpha}^*$}
}
\multiput(40,15)( 0,0){1}{
\put(-35,  0){\makebox(30,12)[r]{$g_{m_{\seq{a}}}(\seq{a})^{\infty}=\hdots$}}
\multiput(0,0)(0,0){1}{
\put(   0,  0){\line(1,0){ 14}}\put(   0,12){\line(1,0){ 14}}
\put(  14,  0){\line(0,1){ 12}}}
\multiput(15,0)(66,0){2}{
\put(   0,  0){\line(1,0){ 65}}\put(   0,12){\line(1,0){ 65}}
\put(   0,  0){\line(0,1){ 12}}\put(  65, 0){\line(0,1){ 12}}}
\put(  33,  4){$g_{m_{\seq{a}}}(\seq{a})$}
\put(  88,  4){$\seq{\alpha}$}
\put( 108,  4){$\seq{\beta}$}
\put(  82, -2){$\underbrace{\makebox(63,1)[c]{}}_{}$}\put( 104,-14){$g_{m_{\seq{a}}}(\seq{a})$}
\multiput(147,0)(0,0){1}{
\put(   0,  0){\line(1,0){ 46}}\put(   0,12){\line(1,0){ 46}}
\put(   0,  0){\line(0,1){ 12}}}
}
\multiput(141,15)(0,4){9}{\line(0,1){2}}
\multiput(161,15)(0,4){9}{\line(0,1){2}}
\end{picture}
\]
\normalsize
The shift-boundedness of $g_{m_{\seq{a}}}(\seq{a})$ gives $\seq{\alpha}^* > \seq{\beta}$ and therefore $g_{m_{\seq{a}}}(\seq{a})^{\infty}<e(\seq{c})$, a contradiction to $\seq{a}\in I(\seq{c})$.

If $ m_{\seq{c}} - k\,m_{\seq{a}} \geq \frac{1}{2} m_{\seq{a}}$ let 
$\seq{\alpha}^* = e(\seq{c})[m_{\seq{c}}+1, (k+1)m_{\seq{a}}]$ and
$\seq{\beta}^* = e(\seq{c})[(k+1)m_{\seq{a}}+1, 2(k+1)m_{\seq{a}}-m_{\seq{c}}]$. 

\footnotesize
\[
\begin{picture}(230,65)
\multiput(40,37)( 0,0){1}{
\put( -35,  0){\makebox(30,12)[r]{$e(\seq{c})=\hdots$}}
\multiput(0,0)(0,0){1}{
\put(   0,  0){\line(1,0){125}}\put(  0,12){\line(1,0){125}}
\put( 125,  0){\line(0,1){ 12}}
\put(  50,  4){$\seq{c}$}}
\multiput(126,0)(0,0){1}{
\put(   0,  0){\line(1,0){ 67}}\put(  0,12){\line(1,0){ 67}}
\put(   0,  0){\line(0,1){ 12}}}
\put( 133,  4){$\seq{\alpha}^*$}
\put( 153,  4){$\seq{\beta}^*$}
}
\multiput(40,15)( 0,0){1}{
\put(-35,  0){\makebox(30,12)[r]{$g_{m_{\seq{a}}}(\seq{a})^{\infty}=\hdots$}}
\multiput(0,0)(0,0){1}{
\put(   0,  0){\line(1,0){ 14}}\put(   0,12){\line(1,0){ 14}}
\put(  14,  0){\line(0,1){ 12}}}
\multiput(15,0)(66,0){2}{
\put(   0,  0){\line(1,0){ 65}}\put(   0,12){\line(1,0){ 65}}
\put(   0,  0){\line(0,1){ 12}}\put(  65, 0){\line(0,1){ 12}}}
\put(  33,  4){$g_{m_{\seq{a}}}(\seq{a})$}
\put( 133,  4){$\seq{\gamma}$}
\put( 153,  4){$\seq{\alpha}$}
\put(  82, -2){$\underbrace{\makebox(63,1)[c]{}}_{}$}\put( 104, -14){$g_{m_{\seq{a}}}(\seq{a})$}
\multiput(147,0)(0,0){1}{
\put(   0,  0){\line(1,0){ 46}}\put(   0,12){\line(1,0){ 46}}
\put(   0,  0){\line(0,1){ 12}}}
}
\multiput(166,15)(0,4){9}{\line(0,1){2}}
\multiput(186,15)(0,4){9}{\line(0,1){2}}
\multiput(206,15)(0,4){9}{\line(0,1){2}}
\end{picture}
\]
\normalsize
The shift-boundedness of $\seq{c}$ gives $\seq{\alpha}^*\seq{\beta}^* > \seq{\gamma}\seq{\alpha}$ and therefore $g_{m_{\seq{a}}}(\seq{a})^{\infty}<e(\seq{c})$, a contradiction to $\seq{a}\in I(\seq{c})$. 

\vspace{1.5ex}
Now, assume that $m_{\seq{c}} < m_{\seq{a}}$. We must have $\seq{c}\leq \seq{a} \leq \seq{c}^{\infty}$ as $e(\seq{c}) \leq \seq{a} < \seq{c}$ contradicts the $e_i$-minimality of $\seq{a}$. Hence $\seq{c}$ is a proper prefix of $\seq{a}$. Lemma \ref{lemma: prefix is minimal} gives $\seq{c}^{\infty} = g_{m_{\seq{c}}}(\seq{a})^{\infty} <  \seq{a}$, a contradiction to $\seq{a}\in I(\seq{c})$. 
\end{proof}

\begin{lemma}
\label{lemma: dim Fc = dim Fa}
For all $\seq{a}\in I(\seq{c})$ we have $\dim_H F(\seq{c}) = \dim_H F(\seq{a})$.
\end{lemma}

\begin{proof}
The statement of the lemma is clear if $\seq{c}$ is an infinite $e_i$-minimal sequence. Hence we may assume that  $\seq{c}$ is a finite $e_{i}$-minimal sequence. 
By Lemma \ref{lemma: Fc = Fcoo} we have that $F(\seq{c}) = F(\seq{c}^{\infty})$ and therefore 
$\dim_H F(\seq{c}) = \dim_H F(\seq{c}^{\infty})$.
Lemma \ref{lemma: dim Fc = dim Ffc} gives that $\dim_H F(\seq{c}) = \dim_H F(f(\seq{c}))$. 
Let $\seq{x}$ be an element in $F(e(\seq{c}))\setminus F(f(\seq{c}))$. As $f(\seq{c}) = \tilde{\seq{c}} \, \seq{c}'$ we have from Lemma \ref{lemma: uuu} that $\seq{x}$ must end with a sequence in $\Ac{\seq{c}}$. Moreover, since $\seq{c}$ is a finite $e_{i}$-sequence we have
\begin{equation}
\label{eq: dim A < dim F}
\dim_H \Ac{\seq{c}}
= \frac{1}{|\seq{c}|}
\leq \frac{1}{2^i}
= \dim_H F(\seq{e}_{i})
\leq \dim_H F(\seq{c}).
\end{equation}
Combining (\ref{eq: dim A < dim F}) with Lemma \ref{lemma: uuu} shows that $\dim_H F(\seq{c}) = \dim_H F(e(\seq{c}))$. 
\end{proof}

\begin{definition}
Let $\textit{IM}(i)$ be the set of infinite $e_i$-minimal sequences and define $\textit{IM} = \cup_{i=1}^{\infty}\textit{IM}(i)$.
\end{definition}

\begin{theorem}
\label{thm: derivative zero}
The derivative of $\phi$ is zero Lebesgue \mbox{a.e.}
\end{theorem}

\begin{proof}
It is clear that the derivative of $\phi$ is zero on an interval $I(\seq{c})$, where $\seq{c}$ is a finite $e_i$-minimal sequence. Hence we only have to show that the set $\textit{IM}$ has Lebesgue measure zero. From Lemma \ref{lemma: min is SB} we have $\textit{IM} \subset\textit{ISB}$, where $\textit{ISB}$ is the set of infinite shift-bounded sequences. Therefore by Lemma \ref{lemma: S measure 0} $\textit{IM}$ has Lebesgue measure 0.
\end{proof}

\begin{theorem}
\label{thm: maximality}
The interval $I(\seq{c})$ is the largest interval $I$ on which we have $\dim_H F(\seq{c}) = \dim_H F(\seq{a})$ for $\seq{a} \in I$.
\end{theorem}

\begin{proof}
By Lemma \ref{lemma: dim Fc = dim Fa} we have that $\dim_H F(\seq{c}) = \dim_H F(\seq{a})$ for all $\seq{a}\in I(\seq{c})$.
For the maximallity, assume first that $\seq{c}$ is a finite $e_1$-minimal sequence. Let $A_{\seq{c}}$ be a transition matrix corresponding to $F(\seq{c})$. Lemma \ref{lemma: e-a_k} gives that there is a sequence $\{\seq{a}_k\}$ of finite $e_1$-minimal sequences growing to $e(\seq{c})$. Let $A_k$ be the transition matrix corresponding to $F(\seq{a}_k)$. From Corollary \ref{cor: A primitiv} it follows that $A_{\seq{c}}$ and $A_k$ are primitive matrices. As $(\seq{a}_k)^{\infty}\in F(\seq{a}_k)\setminus F(\seq{c})$ we have that  $A_k \geq A_{\seq{c}}$, entry by entry, (we may rescale the matrices to have the same size), and where the inequality is strict for at least one pair of indices. As $A_k$ is primitive it follows from the Perron-Frobenius Theorem \ref{thm: perron-frobenius} and Theorem \ref{thm: pesin} that 
\[
\dim_H F(\seq{a}_k)>\dim_H F(\seq{c}),
\]
and therefore the interval $I(\seq{c})$ cannot be extended leftward. Similarly, we use the sequence $\{\seq{b}_k\}$ from Lemma \ref{lemma: e-b_k} to show that $\seq{c}^{\infty}$ is the right endpoint of the interval $I$.

For the case when $\seq{c}$ is an infinite $e_1$-minimal sequence let $\seq{a}$ be a finite $e_1$ minimal sequence. Then by our general assumptions the endpoint of $I(\seq{a})$ are uniquely coded and not infinite $e_1$-minimal. The result now follows by the fact that the intervals constructed from infinite $e_1$-minimal sequence have zero Lebesgue measure.

By Lemma \ref{lemma: e21-min/SB <=> e2i-min/SB} and Theorem \ref{thm: self-similar} we may now extend the result to be valid in any interval $[\seq{e}_{i-1},\seq{e}_i)$ for $i\geq2$.
\end{proof}

\begin{lemma}
\label{lemma: entropy estimate}
Let $\seq{u}$ be a finite shift-bounded sequence and $\seq{u}^{\infty}<\seq{v}$ a sequence such that $F(\seq{v})$ is a sub-shift of finite type and that $\seq{u}$ is a prefix of $\seq{v}$. Let $j$ be the first position such that $(\seq{u}^{\infty})[1,j] \neq \seq{v}[1,j]$ and let $m$ be the largest integer such that $|\seq{u}^m|<j$. Then 
\begin{equation}
\label{eq: entropy estimate}
	\lambda_{\seq{u}} \leq  \lambda_{\seq{v}} \left(1+\frac{2}{\lambda_{\seq{v}}^{m|\seq{u}|-m+1} }\right)
\end{equation}
where $\log \lambda_{\seq{u}}$ is the topological entropy of $F(\seq{u})$ and $\log \lambda_{\seq{v}}$ is the topological entropy of $F(\seq{v})$. 
\end{lemma}

\begin{proof}
A sequence $\seq{x}$ in $F(\seq{u})[1,n]\setminus F(\seq{v})[1,n]$ must contain the sub-sequence $\seq{u}^m$ or $(\seq{u}^m)^*$. Since $\seq{u}$ is a prefix of $\seq{v}$ and $\seq{u}^{\infty}<\seq{v}$ the number $m$ is well defined and $m\geq1$. The number of sequence $\seq{x}$ of length $n$ containing $\seq{u}^m$ or $(\seq{u}^m)^*$ precisely $r$ times is bounded by 
\[	
	2^r \binom{n-r(m|\seq{u}|+m-1)}{r} |F(\seq{v})[1,n-r(m|\seq{u}|+m-1)]|,
\]
as the shift-boundedness of $\seq{u}$ gives that $\seq{u}^m$ or $(\seq{u}^m)^*$ may overlap in at most $m-1$ positions. As $F(\seq{v})$ is a subshift of finite type there is a constant $C_{\seq{v}}$ such that $|F(\seq{v})[1,n]|\leq C_{\seq{v}} \lambda_{\seq{v}}^n$ for $n$ large enough. If summing up we get
\begin{align*}
|F(\seq{u})&[1,n]| \leq \\
&\leq  \sum_{r\geq0} 2^r \binom{n-r(m|\seq{u}|+m-1)}{r} |F(\seq{v})[1,n-r(m|\seq{u}|+m-1)]| \\
&\leq  C_{\seq{v}} \lambda_{\seq{v}}^{n} \sum_{r\geq0} 2^r \binom{n}{r} \frac{1}{ \lambda_{\seq{v}}^{r(m|\seq{u}|-m+1)} }\\
&\leq  C_{\seq{v}} \lambda_{\seq{v}}^{n} \left(1+ \frac{2}{ \lambda_{\seq{v}}^{m|\seq{u}|-m+1} }\right)^n.
\end{align*}
By taking the logarithm on both sides in the equation above, divide by $n$ and then letting $n$ tend to infinity we obtain (\ref{eq: entropy estimate}).
\end{proof}

\begin{theorem}
The map $\phi$ is continuous.
\end{theorem}

\begin{proof} 
By Theorem \ref{thm: pesin} we just have to show that the entropy of $F(\seq{c})$ depends continuously on $\seq{c}$. Let us first consider the case when $\seq{c}$ is a finite $e_i$-minimal sequence. By Lemma \ref{lemma: entropy estimate} we have 
\begin{equation}
\label{eq: ak-ec}
\left|h_{\textnormal{top}}(F(\seq{a}_k)) - h_{\textnormal{top}}(F(e(\seq{c})))\right|\leq \log\left(1+2\lambda_{\seq{c}}^{-k|\seq{c}|} \right),
\end{equation}
where $\seq{a}_k$ is the sequences defined in Lemma \ref{lemma: e-a_k}. Hence, when letting $k$ tend to infinity we have that $\seq{a}_k\to e(\seq{c})$ and that the right hand side of (\ref{eq: ak-ec}) tends to zero, implying the left-continuity in the left endpoint of the interval $I(\seq{c})$. The right-continuity in the left-endpoint of $I(\seq{c})$ follows trivially as the entropy is constant in a neighbourhood to the right of this point. Similarly the right-continuity in the right endpoint of $I(\seq{c})$ is clear. By the equality $F(\seq{c}^k) = F(\seq{c}^\infty)$, the sequence $\seq{w}_k = \seq{c}^k 1^{|\seq{c}|}$ and Lemma \ref{lemma: entropy estimate} we have that 
\[
\left|h_{\textnormal{top}}(F(\seq{c}^{\infty})) - h_{\textnormal{top}}(F(\seq{w}_k))\right|\leq \log\left(1+2\lambda_{\seq{w}_k}^{-k|\seq{c}|} \right),
\]
which implies the left-continuity in the right endpoint of $I(\seq{c})$.

Now assume that $\seq{c}$ is an infinite $e_i$-minimal sequence. Then Lemma \ref{lemma: min is SB} and Lemma \ref{lemma: SB prefix of infinite SB} implies that there is a sequence $\{\seq{u}_k\}$ of finite shift-bounded sequence tending to $\seq{c}$. There is an $m$ such that for $\seq{v}_k = \seq{u}_k 1^{m}$ we have $\seq{u}_k < \seq{c} < \seq{v}_k$, for all $k\geq1$. Again Lemma \ref{lemma: entropy estimate} gives that  
\[
\left|h_{\textnormal{top}}(F(\seq{u}_k)) - h_{\textnormal{top}}(F(\seq{v}_k))\right|
\leq \log\left(1+2\lambda_{\seq{v}_k}^{-|\seq{u}_k|} \right),
\]
which implies the continuity of the entropy in the point $\seq{c}$.

Finally, Corollary \ref{cor: dim e_i} implies the left-continuity in the point $\sigma(\seq{t}')$, where $\seq{t}$ is Thue-Morse sequence. 
\end{proof}

\vspace{1.5ex}

For the rest of the section we turn our interest to the set $\textit{IM}$ of infinite minimal sequences. We define the function $\psi : \Sigma_2^{} \to [0,1]$ by $\psi(\seq{c}) = \dim_H \textit{IM} \cap [\seq{c},1^{\infty}] $. Note that we equally could have defined the function $\psi$ as a function on the real interval $[0,1]$. In comparison to $\phi$ the function $\psi$ is defined on the parameter-space while $\phi$ is a function on the phase-space.

\begin{theorem}
For any sequence $\seq{c}$ we have $\psi(\seq{c}) = \phi(\seq{c})$.	
\end{theorem}

\begin{proof}
Since $\textit{IM}\cap[\seq{c},1^{\infty}]\subset F(\seq{c})$ we have $\psi(\seq{c}) \leq \phi(\seq{c})$. Let us turn to the reversed inequality. Assume that $\seq{c}$ is a finte $e_i$-minimal sequence. From Lemma \ref{lemma: e-b_k} we know that there is a sequence $\{\seq{b}_k\}$ of $e_i$-minimal sequences tending to $\seq{c}^{\infty}$. Define
\[
N_k(\seq{c}) = \left\{ \seq{c}\, \seq{u} : \seq{u}\in [\seq{b}_k]\cap
F(\seq{b}_k)\right\},
\]
where $[\cdot]$ denotes the cylinder-set. Let $\seq{x}\in N_k$. Note that \seq{x} has the prefix $\seq{b}_{k+1}$. 
By Lemma \ref{lemma: e-b_k} we have that $\seq{e}_{i-1} < \seq{c} < \seq{x} < \seq{b}_k < \seq{e}_i$, so $\seq{x}$ is an $e_i$-sequence. The $e_i$-minimality of $\seq{b}_{k+1}$ and that $\seq{x}[n+1,n+|\seq{b}_{k+1}|]> \seq{b}_{k+1}$ for $n\geq 2^i$ gives that $g_n(\seq{x})^{\infty} < \seq{x}$ for all $n\geq 2^i$ such that $x_n = 1$.

For the cases when $n$ is such that $x_n =0$ Lemma \ref{lemma: e-b_k} gives that $e(g_n(\seq{x})) > \seq{x}$ for $2^i\leq n \leq (k+1)|\seq{c}|$. 
For $(k+1)|\seq{c}| < n$ let 
$\seq{v}^* = e(g_n(\seq{x}))[n+1,n+(k+1)|\seq{c}|]$ and
$\seq{w}   = \seq{x}[n+1,n+(k+1)|\seq{c}|]$. 
Then $\seq{v}^* = (\seq{c}^{k+1})^*$ and $(\seq{b}_k^*)^{\infty}\geq  \seq{w}$ by the definition of $N_k(\seq{c})$ and Lemma \ref{lemma: Fc = Fcoo}. This implies
\[
\seq{v}^* = (\seq{c}^{k+1})^* >\left( (\seq{c}^{k}p(\seq{c}))^*\right)^{\infty} = (\seq{b}_k^*)^{\infty} \geq  \seq{w}
\]
and therefore $e(g_n(\seq{x}))> \seq{x}$ for $n\geq 2^i$ such that $x_n = 0$. Hence $\seq{x}$ is an infinite $e_i$-minimal sequence. We have  $\textit{IM}\cap[\seq{c},1^{\infty}]\supset N_k(\seq{c})$ and $\dim_H N_k(\seq{c}) = \dim_H F(\seq{b}_k)$. By choosing $k$ sufficiently large we have $\dim_H N_k(\seq{c})$ arbitrarily close to $\phi(\seq{c})$.
\end{proof}

\begin{corollary}
\label{cor: dim Gamma = dim F} 
Let $\textit{ISB}$ be the set of all infinite shift-bounded sequences. Then 
\[	\dim_H \textit{IM}\cap[\seq{c},1^{\infty}]
	= \dim_H \textit{ISB}\cap[\seq{c},1^{\infty}] 
	= \dim_H \Gamma\cap[0,\seq{c}'] 
	= \dim_H F(\seq{c}),
\]
where $\Gamma$ is the set defined in (\ref{eq: allouche set}).
\end{corollary}

\section{Numerics}

By characterising the dimension of $F(c)$ via the spectral radius of a primitive transition matrix the problem of numerically calculate an approximative value of $\phi$ reduces to calculate the eigenvalues of the transition matrix.

\begin{figure}
\[
\begin{picture}(230,165)
\put(   0,   0){\includegraphics[width=225pt]{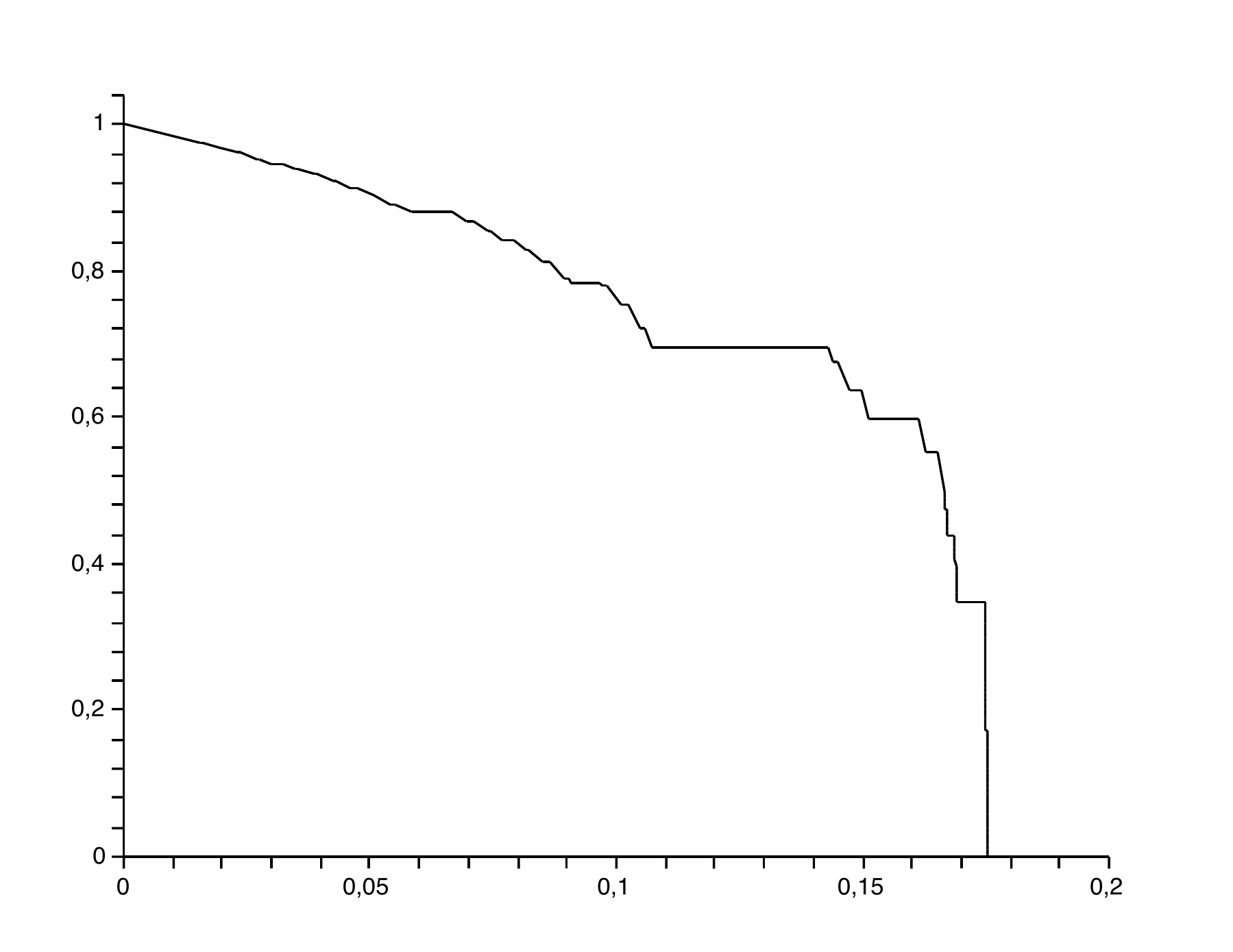}}
\put( 214,  15){\footnotesize$c$}
\put(   6, 164){\footnotesize$\dim_H F(c)$}
\end{picture}
\]
\caption{The graph of $\phi(c) = \dim_H  F(c)$.}
\label{graph}
\end{figure}

The graph of $\phi$, (see figure \ref{graph}) was calculate by considering $e_1$-minimal sequence of length at most 8, which gives transition matrices of size $128\times128$, and then using Theorem \ref{thm: self-similar} to obtain the values of $\phi$ for $e_i$-minimal sequences with $i>1$.
A finer subdivision of the interval [0,1] would require harder calculation as the runtime complexity of the computation is exponential in the length of the minimal sequences.

\section{Acknowledgements}

The author would like to express his gratitude to J. Schmeling and T. Persson for their appreciated comments and remarks.

\end{document}